\documentclass[12pt]{article}
\usepackage{amsmath}
\usepackage{graphicx,psfrag,epsf}
\usepackage{enumerate}
\usepackage{natbib}
\usepackage{url} 

\newcommand{\blind}{0}

\usepackage[english]{babel}
\usepackage{amsthm}
\newtheoremstyle{colon}%
{}
{}
{\itshape}
{}
{\bfseries}
{:}
{ }
{}

\theoremstyle{colon}
\newtheorem{ass}{Assumption}
\newtheorem{theorem}{Theorem}
\newtheorem{example}{Example}
\newtheorem{definition}{Definition}

\newtheorem{rem}{Remark}
\newtheorem{lemma}{Lemma}
\newtheorem{cor}{Corollary}
\newtheorem{fact}{Fact}

\newtheorem{innercustomgeneric}{\customgenericname}
\providecommand{\customgenericname}{}
\newcommand{\newcustomtheorem}[2]{%
  \newenvironment{#1}[1]
  {%
   \renewcommand\customgenericname{#2}%
   \renewcommand\theinnercustomgeneric{##1}%
   \innercustomgeneric
  }
  {\endinnercustomgeneric}
}

\newcustomtheorem{customex}{Example}
\newcustomtheorem{customprop}{Proposition}
\newcustomtheorem{customlem}{Lemma}
\newcustomtheorem{customrem}{Remark}
\newcustomtheorem{customdef}{Definition}
\newcustomtheorem{customcor}{Corollary}

\newcounter{bean}

\usepackage{scalerel,stackengine}			
\usepackage{amssymb}
\usepackage{authblk}
\usepackage{amsmath, amsfonts, amssymb, bm, bbm}
\usepackage{listings}
\usepackage{float}
\usepackage[titletoc]{appendix}
\usepackage{hvfloat}
\usepackage{longtable}
\usepackage{mathtools}
\usepackage{booktabs,xcolor,colortbl}
\usepackage{xparse}
\usepackage[bottom]{footmisc}
\usepackage[font=small,skip=0pt]{caption}
\usepackage{mathrsfs}
\usepackage{enumitem}
\usepackage{rotating}
\usepackage{comment}
\usepackage{accents}
\usepackage{csquotes}
\usepackage{booktabs}
\usepackage{threeparttable}
\usepackage{setspace}
\excludecomment{confidential}
\newcommand\sbullet[1][.5]{\mathbin{\vcenter{\hbox{\scalebox{#1}{$\bullet$}}}}}

\newcommand{\aggregate}[2]{\underset{#2}{\operatornamewithlimits{#1\ }}}

\newcommand{\zenodoLink}{\url{https://zenodo.org/records/15089775}}
\let\originalleft\left
\let\originalright\right
\renewcommand{\left}{\mathopen{}\mathclose\bgroup\originalleft}
\renewcommand{\right}{\aftergroup\egroup\originalright}

\begin{document}

\def\spacingset#1{\renewcommand{\baselinestretch}%
{#1}\small\normalsize} \spacingset{1}

\sloppy 


\if0\blind
{\title{\bf ARMAr-LASSO: Mitigating the Impact of Predictor
Serial Correlation on the LASSO}
  \author{\small Simone Tonini\thanks{Corresponding author: simone.tonini@santannapisa.it} \hspace{.2cm}\\
Institute of Economics \& L'EMbeDS, Sant'Anna School\\
 of Advanced Studies, Pisa, Italy\\
    and \\
    Francesca Chiaromonte\hspace{.2cm}\\
    Institute of Economics \& L'EMbeDS,\\
    Sant'Anna School of Advanced Studies, Pisa, Italy\\
Department of Statistics, \\
Penn State University, USA\\
  and  \\
    Alessandro Giovannelli\hspace{.2cm}\\
    Department of Information Engineering,
 Computer Science and Mathematics, \\
 University of L'Aquila, L'Aquila, Italy.}

\date{} 

  \maketitle
} \fi

\if1\blind
{ \bigskip
  \bigskip
  \bigskip
  \begin{center}
    {\LARGE\bf On the Impact of Predictor Serial Correlation 
on the LASSO}
\end{center}
  \medskip
} \fi
\bigskip
\begin{abstract}
We explore estimation and prediction accuracy for sparse linear models, focusing on scenarios where both predictors and errors carry serial correlations. We establish a clear link between predictor serial correlation and the performance of the LASSO, showing that even orthogonal or weakly correlated stationary AR processes can lead to significant spurious correlations due to their serial correlations. To address this challenge, we propose a novel approach named ARMAr-LASSO ({\em ARMA residuals LASSO}), which applies the LASSO to predictors that have been pre-whitened with ARMA filters and lags of dependent variable. We establish oracle inequalities for the feasible estimator and assess its estimation, selection, and forecasting performance empirically. Our findings are supported by extensive simulations and an application to real-world macroeconomic data, which highlight the superior performance of the ARMAr-LASSO for handling sparse linear models in the context of time series.
\end{abstract}

\noindent%
\vfill

\newpage
\spacingset{1.8} 

\section{Introduction}\label{sec:intro}
The LASSO (\citealt{tibshirani96}) is perhaps the most commonly employed approach to handle regressions with a large number of predictors. From a theoretical standpoint, its effectiveness in terms of estimation, prediction, and feature selection is contingent upon either orthogonality or reasonably weak correlation among predictors (see~\citealt{zhaoyu2006, Bickel2009, Negahban2012}). This hinders the use of the LASSO for the analysis of economic time series data, which are notoriously characterized by intrinsic multicollinearity; that is, by predictor correlations at the population level (\citealt{Lippi2000, SeW2002a, demol2008}). A common procedure to address this issue is to model multicollinearity and remove it, as proposed, e.g., by~\cite{Fan2020}, who filter time series using common factors and then apply the LASSO to the filtered residuals. 
However, mitigating or even eliminating multicollinearity is not the end of the story, as effectiveness of the LASSO can also be affected by spurious correlations. These occur when predictors are orthogonal or weakly correlated at the population level, but a lack of sufficient independent replication (lack of degrees of freedom) introduces correlations at the sample level, potentially leading to false scientific discoveries and incorrect statistical inferences (\citealp{Fan16}). This issue has been broadly explored in ultra-high dimensional settings, where the number of predictors can vastly exceed the available sample size (\citealp{Fan13}).

This article introduces two elements of novelty. First, we establish an explicit link between predictor serial correlation and LASSO performance. Building on~\cite{Bartlett35,Mcgregor1965}, we show how predictor serial correlation amplifies spurious sample correlations in an AR(1) benchmark. Monte Carlo simulations illustrate this phenomenon for Gaussian AR(1) processes and show that it persists in more general settings, including cross-correlated and non-Gaussian ARMA processes. We then investigate the implications for sparse regression. Monte Carlo experiments show that increasing predictor serial correlation can deteriorate LASSO performance, even when predictors are mutually independent in the population.  
Second, we propose an approach that rescues LASSO performance in the presence of serially correlated predictors. Our approach, which we name ARMAr-LASSO ({\em ARMA residuals LASSO}), relies upon a working model where, instead of the observed predictor time series, we use as regressors the residuals of ARMA processes fitted on such series, augmented with lags of the dependent variable. We employ the mixingale and near-epoch dependence framework (\citealp{davidson1994, ADAMEK2023}) to prove oracle inequalities for the estimation and prediction error of the ARMAr-LASSO, while simultaneously addressing the issue of estimating ARMA residuals. To complete the analysis, we use simulations to validate and generalize theoretical results. We illustrate the method by forecasting Euro Area inflation from a high-dimensional predictor panel. Simulations and empirical analyses show ARMAr-LASSO produces more parsimonious models, better coefficient estimates, and more accurate forecasts than LASSO-based benchmarks. 

On the serial correlation front, most of the theoretical econometric literature has focused on its impact in the error terms, particularly regarding post-LASSO inference (see, e.g., \citealp{Chernozhukov2021,chronopoulos2023,ADAMEK2023}). The present study shows that serial correlation in the predictors deserves similar attention, as it can adversely affect both the estimation and prediction accuracy of the LASSO. Our work complements the vast literature on error bounds for LASSO-based methods in time series analysis, which addresses estimation and prediction consistency in scenarios with autocorrelated errors and autoregressive processes (see, e.g., \citealp{nardi2011,Uematsu2019}). Such scenarios are ubiquitous, e.g., they are easily found in US and Euro Area monthly macroeconomic data (see~\citealp{McCracken2016} and~\citealp{Ale2021}). Moreover, our methodology is consistent with the existing literature on pre-whitening filters, which aim to mitigate autocorrelation and multicollinearity by applying LASSO or related methods to filtered residuals (see, e.g., \citealp{
medeiros2017ADA,Hansen2019,Fan2020,chronopoulos2023}). 
 
The remainder of the article is organized as follows. Section~\ref{sec:ProbSet} introduces the problem setup and our analysis concerning the link between serial correlations and LASSO performance. Section~\ref{sec:armar_lasso} introduces the ARMAr-LASSO and explores its theoretical properties. Section~\ref{sec:Exp} presents simulations and real data analyses to evaluate the proposed methodology. Section~\ref{ConcRem} provides some final remarks. The Supplement encompasses the proofs of theoretical results, additional studies, and simulations excluded from the main manuscript. \sloppy
Code and replicability 
materials are at \zenodoLink

We summarize here some notation that will be used throughout. Bold letters denote vectors, for example $\mathbf{a}=(a_1,\dots,a_p)'$. Supp$(\mathbf{a})$ denotes the support of a vector, that is, $\{i\in\{1,\dots,p\}:a_i\neq0\}$, and $|\text{Supp}(\mathbf{a})|$ the support cardinality. The $\ell_q$ norm of a vector is $||\mathbf{a}||_q\coloneqq\left(\sum_{j=1}^p|a_j|^q\right)^{1/q}$ for $0<q<\infty$, with $||\mathbf{a}||_q^k\coloneqq\left(\sum_{j=1}^p|a_j|^q\right)^{k/q}$, and with the usual extension $||\mathbf{a}||_0\coloneqq|\text{Supp}(\mathbf{a})|$. Bold capital letters denote matrices, for example $\mathbf{A}$, where $\left(\mathbf{A}\right)_{ij}=a_{ij}$ is the $i$-row $j$-column element. Furthermore, $\pmb{0}_{p}$ denotes a $p$-length vector of zeros, $\pmb{I}_{p}$ the $p\times p$ identity matrix. 
To simplify the presentation, we frequently use $C$ to indicate arbitrary positive finite constants.

\section{Problem Setup}\label{sec:ProbSet}
\sloppy
Consider the linear regression model
\begin{equation}
y_t=\mathbf{x}_t'\pmb{\alpha}^*+\varepsilon_t\ \ , \label{DGP_y}\ \ \ \ \ \ t=1, \ldots, T\ \ ,
\end{equation}
where $\mathbf{x}_t=(x_{1,t},\ldots,x_{n,t})'$ is a $n\times1$ vector of predictors,  $\pmb{\alpha}^*$ is a $n\times1$ unknown $s$-sparse vector of regression coefficients, i.e. $||\pmb{\alpha}^*||_0=s<n$, and $\varepsilon_t$ is an error term. We impose the following assumptions on the processes $\left\{\mathbf{x}_t\right\}$ and $\left\{\varepsilon_t\right\}$.

\begin{ass}\label{ass:CovStat}
(a) $\left\{\mathbf{x}_t\right\}$ and $\left\{\varepsilon_t\right\}$ are 
causal, invertible stationary ARMA processes
\begin{eqnarray}
 x_{i,t}&=&\sum_{\ell=1}^{p_i}\phi_{i,\ell}x_{i,t-\ell}+\sum_{k=1}^{q_i}\theta_{i,k}u_{i,t-k}+u_{i,t}\ \ ,\label{DGP_x}
\ \ \ \  i=1,\dots,n\ \ ,\ \ \ \ p_i,q_i<\infty\ \ ,  \\
 \varepsilon_t &=&\sum_{\ell=1}^{p_{\varepsilon}}\phi_{\varepsilon,\ell}\varepsilon_{t-\ell}+\sum_{k=1}^{q_{\varepsilon}}\theta_{\varepsilon,k}\omega_{t-k}+\omega_t\ \ ,\ \ \ \ \ \ p_{\varepsilon},q_{\varepsilon}<\infty\ \ .\label{DGP_varep}
\end{eqnarray}
(b) The innovation processes $u_{i,t}\sim i.i.d.(0,\sigma_i^2)$, $\omega_{t}\sim i.i.d.(0,\sigma_{\omega}^2)$, where $u_{i,t}\perp u_{j,t-\ell}$ for any $i\neq j,\ t$ and $\ell\neq0$; and $u_{i,t-\ell}\perp\omega_t$ for any $i,t$ and $\ell$. Moreover, the innovation vectors
$(u_{1,t},\ldots,u_{n,t},\omega_t)'$
are i.i.d. over time.\end{ass}

\noindent There are several approaches to estimate a sparse $\pmb{\alpha}^*$ 
(\citealp{zhang2012,James2013}); here we focus on the LASSO estimator (\citealp{tibshirani96}) given by $\widehat{\pmb{\alpha}}= \aggregate{argmin}{\pmb{\alpha}\in \mathbb{R}^n}\:\left\{\frac{1}{2T}||\mathbf{y}-\mathbf{X}'\pmb{\alpha}||_2^2 + \ddot\lambda||\pmb{\alpha}||_1 \right \},$ where $\mathbf{y}=(y_1,\dots,y_T)'$ is the $T\times1$ response vector, $\mathbf{X}=(\mathbf{x}_1,\dots,\mathbf{x}_T)$ is the $n\times T$ design matrix, and $\ddot\lambda>0$ is the weight of the $\ell_1$ penalty and must be ``tuned'' to guarantee that regression coefficient estimates are effectively shrunk to zero -- thus ensuring predictor, or feature, selection.

However, linear associations among predictors are well known to affect LASSO performance.~\cite{Bickel2009,Buhlmann2011} and~\cite{Negahban2012} have shown that the LASSO estimation and prediction accuracy are inversely proportional to the {\em minimum eigenvalue} of the predictor sample covariance matrix. 
Thus, highly correlated predictors deteriorate estimation and prediction performance.
Moreover, ~\cite{zhaoyu2006} proved 
that the LASSO struggles to differentiate between \textit{relevant} (i.e., $\{i\in\{1,\dots,n\}:\alpha_i^*\neq0\}$) and \textit{irrelevant} (i.e., $\{i\in\{1,\dots,n\}:\alpha_i^*=0\}$) predictors when they are closely correlated, leading to false positives. Thus, highly correlated predictors may also deteriorate feature selection performance. The \textit{irrepresentable condition} addresses this issue ensuring both estimation and feature selection consistency through bounds on the sample correlations between relevant and irrelevant predictors \citep{zhaoyu2006,Buhlmann2011}. Nevertheless, orthogonality or weak correlation seldom hold in the context of economic and financial data. For instance, decades of literature provide evidence for co-movements of macroeconomic variables (\citealp{Lippi2000,Lippi2005,SeW2002a,SeW2002b}). Special methods have been proposed to mitigate the negative effects of these linear associations, such as Factor-Adjusted Regularized Model Selection (FarmSelect) (\citealp{Fan2020}), which applies the LASSO to the idiosyncratic components of economic variables, obtained by filtering the variables through a factor model. Although approaches such as FarmSelect can be very effective in addressing multicollinearity, strong spurious correlations can emerge at the sample level and affect the LASSO even when predictors are orthogonal or weakly correlated at the population level. 
Sample-level spurious correlations can be particularly prominent in regressions with many predictors, especially if the sample sizes are relatively small, and the problem can be yet more serious for time series data, where independent replication can be further hindered by serial correlations.
This is precisely the focus of this article. In the next section, we present a theoretical result, building on established time-series literature~\citep{Bartlett35,Mcgregor1965}, that links serial correlation within time series to sample correlation between them.

\subsection{Serial Correlation and Spurious Sample Correlations}
\label{sec:spurious_correlations}
Consider a first order $n$-variate autoregressive process $\mathbf{x}_t=\pmb{\phi}\mathbf{x}_{t-1}+\mathbf{u}_t$, $t=1,\dots,T$, where $\pmb{\phi}$ is the $n\times n$ diagonal matrix with $\text{diag}(\pmb{\phi})=(\phi_1,\dots,\phi_n)$. 
Let $\mathbf{u}_t$ satisfy Assumption~\ref{ass:CovStat}(b), with mutually independent coordinates, $\sigma_i^2>0$, and $|\phi_i|<1$ for all $i$. Assume that $\mathbf{x}_0$ follows the joint stationary distribution, independently of subsequent innovations. Then $\mathbf{x}_t$ is stationary with $E(\mathbf{x}_t)=\mathbf{0}_n$ and covariance matrix $\mathbf{C}_x$, where $(\mathbf{C}_x)_{ii}=\sigma_i^2/(1-\phi_i^2)$ and $(\mathbf{C}_x)_{ij}=0$ for $i\neq j$.
Let $\widehat{\mathbf{C}}_x=T^{-1}\mathbf{X}\mathbf{X}'$ denote the uncentered empirical Gram matrix, with eigenvalues $\widehat{\psi}_{\max}^x\geq\cdots\geq\widehat{\psi}_{\min}^x$. Finally, denote Pearson's sample correlation by $\widehat{\rho}_{ij}^x$.

\begin{lemma}
\label{lem:sample_correlations}
Under the stationary AR(1) model above, for $i\neq j$, define $z_{k,t}=\sigma_k^{-1}\sqrt{1-\phi_k^2}x_{k,t}$, $\ddot\phi=\phi_i\phi_j$, and $A_{ij,T}=T^{-1}\sum_{t=1}^Tz_{i,t}z_{j,t}$. Then: (a) For every $T\geq1$, $\operatorname{Var}(A_{ij,T})=T^{-1}[1+2\sum_{\ell=1}^{T-1}(1-\ell/T)\ddot\phi^\ell]$. (b) As $T\to\infty$, $\sqrt{T}\widehat{\rho}_{ij}^x\xrightarrow{d}\mathcal{N}\left(0,(1+\ddot\phi)/(1-\ddot\phi)\right)$
\end{lemma}

\begin{figure}
\graphicspath{{images/}}
\centering
{\includegraphics[width=5.4cm]{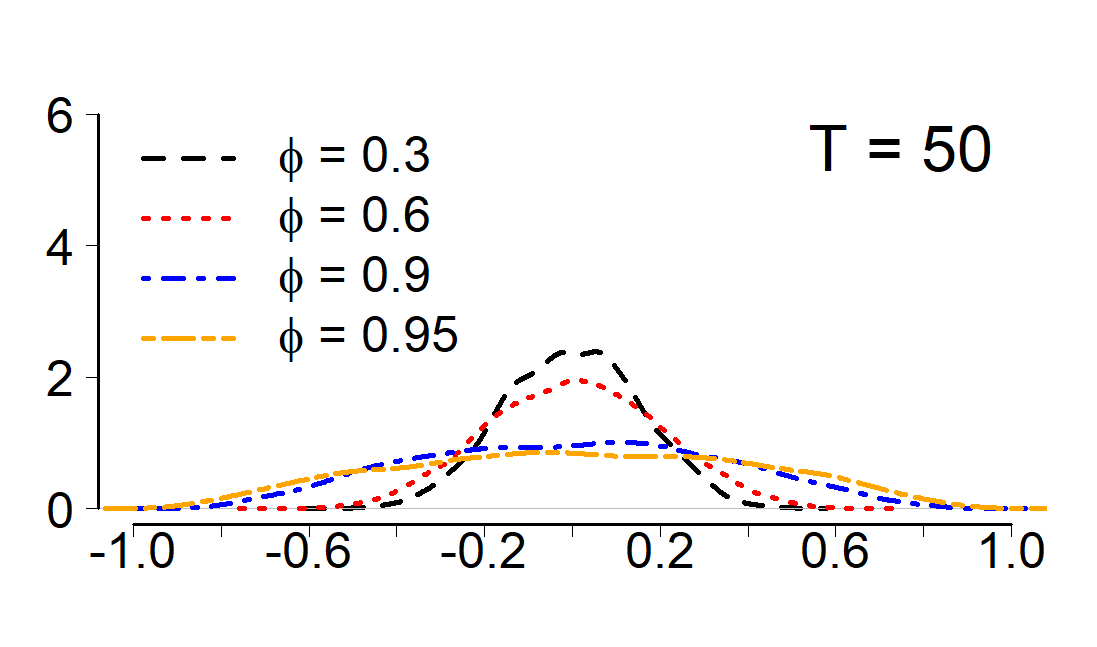}}\hfil
\subfloat{\includegraphics[width=5.4cm]{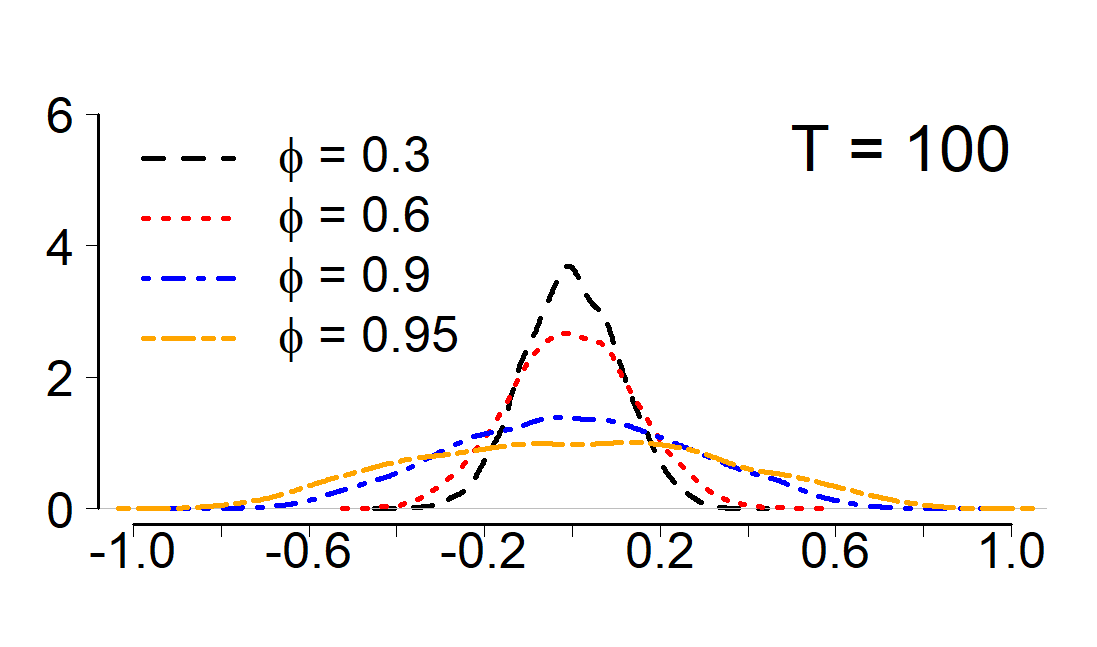}}\hfil
\subfloat{\includegraphics[width=5.4cm]{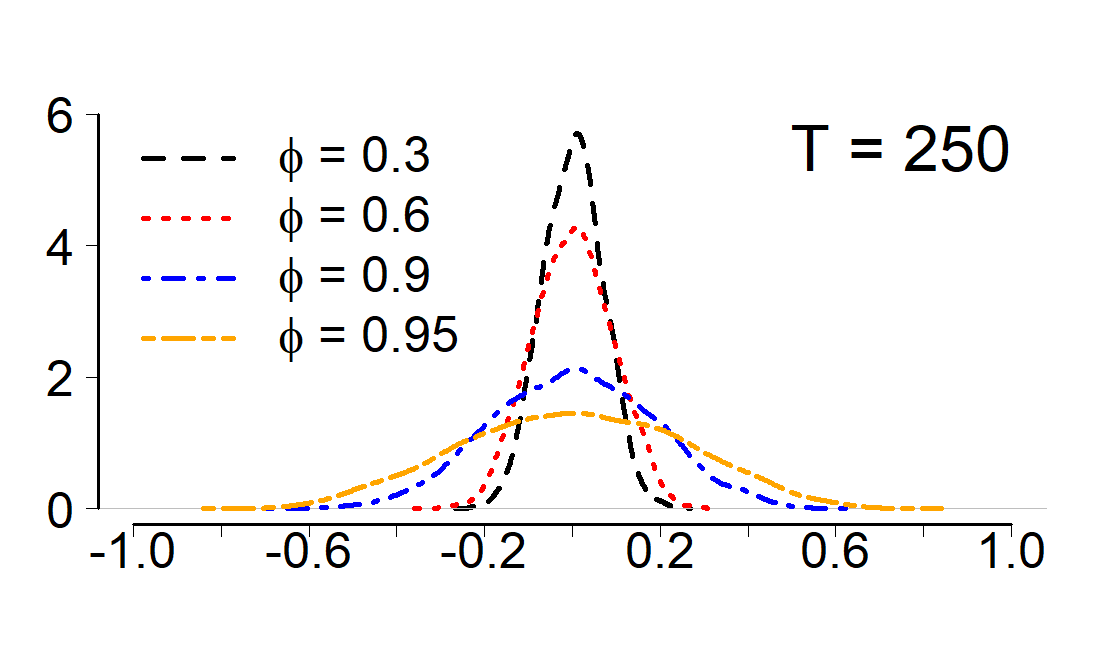}}
\caption{\footnotesize Densities of $\widehat{\rho}_{12}^x$ based on 5000 Monte Carlo replications for different values of $T$ and $\phi$.}\label{fig:DistributionsMC}
\end{figure}

Lemma~\ref{lem:sample_correlations} establishes the exact variance of standardized cross-products and the asymptotic variance of Pearson's correlation. Positive $\ddot\phi$ increases the latter relative to independent observations, whereas negative $\ddot\phi$ reduces it. Figure~\ref{fig:DistributionsMC} reports Monte Carlo densities of $\widehat{\rho}_{12}^x$ from 5000 replications in the case where $n=2$ and $\mathbf{u}_t\overset{\mathrm{i.i.d.}}{\sim}\mathcal{N}(\mathbf{0}_2,\mathbf{I}_2)$, with $\phi_1=\phi_2\in\{0.3,0.6,0.9,0.95\}$ and $T\in\{50,100,250\}$. The simulations show increasing dispersion and more frequent large absolute correlations as persistence increases. Supplement~\ref{MonteCarlo} explores additional settings, including correlated, non-Gaussian, and ARMA processes.

\subsection{Implications for Sparse Regression}
\label{sec:sparse_regression}
The performance of the LASSO depends on the properties of the empirical design matrix, particularly on suitable compatibility or restricted eigenvalue conditions~\citep{Bickel2009,Buhlmann2011}. 
Let $S=\{j:\alpha_j^*\neq0\}$, with $s=|S|\geq1$. The compatibility constant associated with $S$ is $\gamma_x^2= \underset{\|\pmb{\alpha}_{S^c}\|_1\leq3\|\pmb{\alpha}_{S}\|_1;\ \pmb{\alpha}\in\mathbb R^{n}\backslash\{0\}}{\text{min}}\frac{s\,\pmb{\alpha}'\mathbf{X}\mathbf{X}'\pmb{\alpha}}{T\|\pmb{\alpha}_{S}\|_1^2}.$ The compatibility condition (\citealp{Bickel2009}) requires $\gamma_x^2>0$, ensuring positive definiteness restricted to the cone $||\pmb{\alpha}_{S^c}||_1\leq3||\pmb{\alpha}_S||_1$. Assume that each row of $\mathbf{X}$ is standardized to have mean 0 and variance 1. If $T^{-1}\mathbf{XX}'$ is nonsingular, then 
$\gamma_x^2\geq\widehat{\psi}_{\min}^x>0$. Thus, the minimum eigenvalue of $\frac{1}{T}\mathbf{XX}'$ is a lower bound on the compatibility constant, so the compatibility condition is considerably weaker than assuming $\frac{1}{T}\mathbf{XX}'$ to be positive definite. The following remark highlights how strong predictor correlations can adversely affect compatibility conditions, as shown in previous studies.

\begin{rem}
For a standardized design, large sample correlations can reduce the minimum eigenvalue, since $\Pr\bigl(\widehat{\psi}_{\min}^{x}\leq 1-\tau\bigr)\geq\Pr\bigl(|\widehat{\rho}_{ij}^{x}|\geq\tau\bigr)$, for $\tau\in[0,1)$ (see Supplement~\ref{UppMinPsy}). Although $\gamma_x^2\geq\widehat{\psi}_{\min}^{x}$ does not imply a monotonic relationship, \citet{BUHLMANN2013} highlight that strong sample correlations involving active predictors can also reduce the compatibility constant. Moreover, \citet[Ch.~6, p.~106]{Buhlmann2011} relate the compatibility condition to a condition on the minimum eigenvalue. The following simulations investigate these relationships. 
\end{rem}

\begin{figure}[h]
\graphicspath{{images/}}
\centering
\includegraphics[width=15cm]{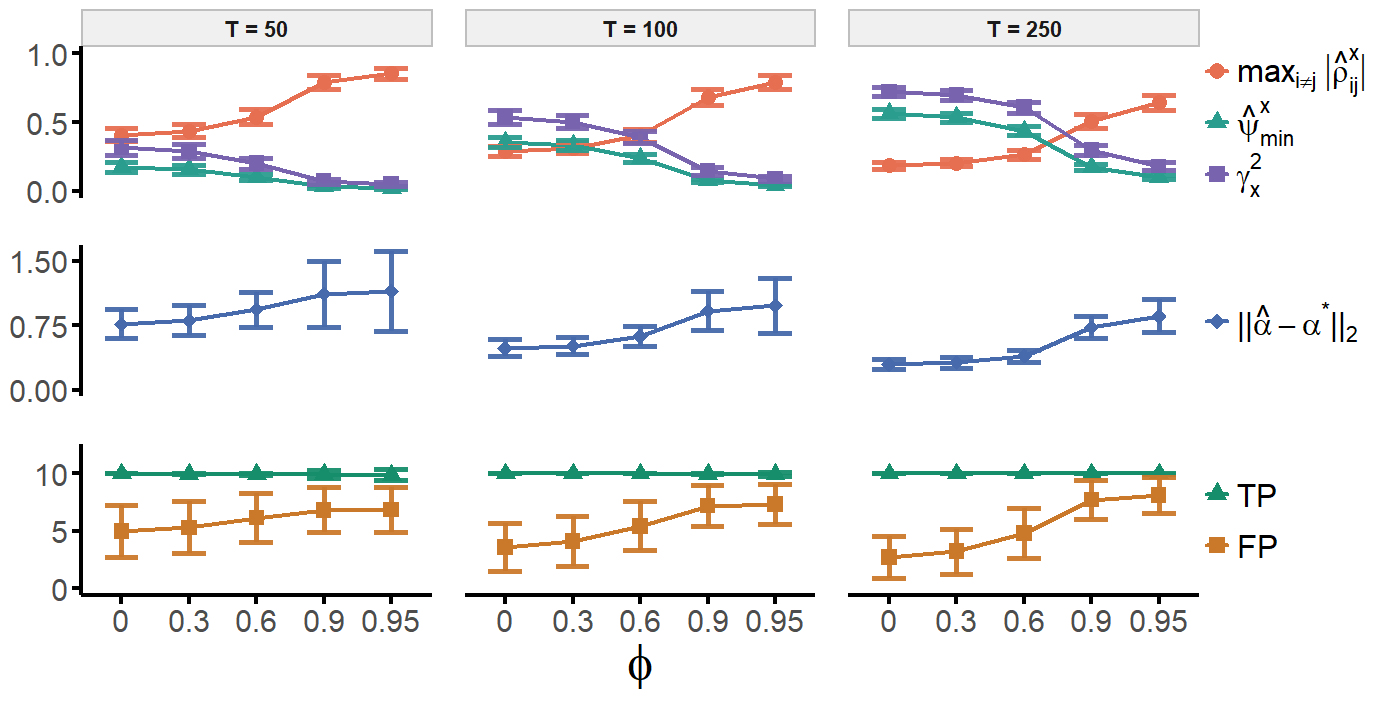}
\caption{\footnotesize Numerical ``toy example'' based on
5000 Monte Carlo replications.
Columns correspond to $T=50$, $100$, and $250$.
From top to bottom, rows report:
(i) maximum absolute sample correlation, minimum eigenvalue,
and compatibility constant;
(ii) $\ell_2$ coefficient estimation error;
and (iii) the numbers of true and false positives.
Points represent Monte Carlo averages, with error bars
indicating standard deviations.}
\label{Fig:toy}
\end{figure}

We conduct a Monte Carlo experiment based on the 20-dimensional AR(1) process $\mathbf{x}_t=\phi\mathbf{x}_{t-1}+\mathbf{u}_t$, where $\mathbf{u}_t\overset{\mathrm{i.i.d.}}{\sim}\mathcal{N}(\mathbf{0}_{20},\mathbf{I}_{20})$. To assess LASSO performance, we generate $y_t=\mathbf{x}_{t-1}'\pmb{\alpha}^*+\varepsilon_t$, where $\pmb{\alpha}^*=(\pmb1_{10}',\pmb0_{10}')'$, and $\varepsilon_t=\phi\varepsilon_{t-1}+\omega_t$, with $\omega_t\overset{\mathrm{i.i.d.}}{\sim}\mathcal{N}(0,\sigma_{\omega}^2)$ and $\sigma_{\omega}^2$ is set to fix the population signal-to-noise ratio at 10. We consider $\phi\in\{0,0.3,0.6,0.9,0.95\}$ and $T\in\{50,100,250\}$, with 5000 replications per configuration. We deliberately focus on the low-dimensional setting ($n<T$) to isolate the effects of serial correlation on the design properties and LASSO performance, while highlighting that substantial spurious sample correlations can arise even in the absence of high dimensionality. High-dimensional settings are investigated extensively in the Monte Carlo experiments of Section~\ref{sec:Exp}. Fixing $S=\{1,\ldots,10\}$, Figure~\ref{Fig:toy} shows that increasing $\phi$ amplifies spurious sample correlations and reduces the minimum eigenvalue and compatibility constant of the sample-standardized design. These effects persist under strong serial correlation as $T$ increases, despite predictors being independent in the population. Moreover, increasing $\phi$ raises the LASSO coefficient estimation error and false positives, while true positives remain nearly unchanged. 

\section{The ARMAr-LASSO}\label{sec:armar_lasso}
We introduce ARMAr-LASSO, designed to improve sparse regression with serially correlated predictors. ARMAr-LASSO is formulated as a two-step procedure. In the first step we estimate a univariate ARMA model on each predictor. In the second step, we run the LASSO using, instead of the original predictors, the residuals from the ARMA model, i.e. estimates of the $u$'s in equation~\eqref{DGP_x}, plus lags of the response. 
We start by introducing the ``working model'' on which our proposal relies; that is, the model that contains the true, 
non-observable ARMA residuals (their estimation will be addressed later)
\begin{equation}\label{WM}
y_t= \mathbf{w}_t'\pmb\beta^* + v_t \ \ .
\end{equation}
Model~\eqref{WM} is the linear projection of $y_t$ on $\mathbf{w}_t=(u_{1,t},\dots,u_{n,t},y_{t-1},\dots,y_{t-p_y})'$. $\pmb\beta^*=(\pmb\alpha^{*'},\phi_{y,1},\dots,\phi_{y,p_y})'$ represents the corresponding best linear projection coefficients and $v_t$ is the error term, which is unlikely to be $i.i.d.$. The model includes $p_y$ lags of $y_t$, where $p_y$ is arbitrary, and some lags may be relevant while others may not. The relevant lags will be directly selected using LASSO. Model~\eqref{WM} is misspecified, in the sense that it does not correspond to the true data generating process (DGP) for the response. 
The idea behind model~\eqref{WM} is to leverage the serial independence of the $u$ terms, thereby avoiding the risk of sizeable spurious correlation. However, the $u$ terms alone may explain only a small portion of the variance of $y_t$, particularly in situations with high persistence. Response lags are included to capture persistence and improve forecasting. 

\begin{rem}\label{rem:p_y}
When $p_y=0$, model~\eqref{WM} relies solely on $u_t$, leaving serial dependence in $v_t$ and reducing the signal-to-noise ratio. Including lags of $y_t$ allows this dependence to be exploited. Under the common AR($p$) restriction (see Remark~\ref{rem1} below), $p$ lags yield white-noise $v_t$, whereas $p_y<p$ impairs forecasting. Although not formally established, simulations show improved performance over benchmarks when $p_y>p$ under common AR($p$) (Supplement~\ref{sec:FiltersExp}).
\end{rem}

\begin{fact}\label{fact1}
By linear projection, $E(\mathbf{w}_tv_t)=\pmb0$. Assumption~\ref{ass:CovStat} and causality imply $E(\mathbf{u}_ty_{t-\ell})=\pmb0$ for $\ell\geq1$ and $E(\mathbf{u}_ty_t)=\mathbf{C}_u\pmb{\alpha}^*$. The coefficients on the innovations in the population linear projection equal $\pmb{\alpha}^*$ whenever $\mathbf{C}_u$ is nonsingular, despite working-model misspecification.
\end{fact}

\begin{fact}\label{fact2}
By linear projection, $E(y_{t-\ell}v_t)=0$ for $1\leq\ell\leq p_y$, although $E(v_t\mid y_{t-1},y_{t-2},\ldots)$ may be nonzero. Moreover, $E(y_{t-\ell}\mid y_{t-\ell-1},y_{t-\ell-2},\ldots)$ may be nonzero, reflecting dependence among response lags, whose exact selection is not our primary objective.
\end{fact}

Next, we provide three illustrative examples. In the first, and simplest, predictors and error terms have an AR(1) representation with a common coefficient; we refer to this as the \textit{common AR(1) restriction} case. In the second, the AR(1) processes have different autoregressive coefficients. In the third, predictors admit a common factor representation with AR(1) idiosyncratic components. Note that in all the examples $p_y=1$.
\begin{example}\label{ex:commonAR1}
Suppose both predictors and error terms in model~\eqref{DGP_y} admit a common AR(1) representation; that is,  $x_{i,t}=\phi x_{i,t-1}+u_{i,t}$ and $\varepsilon_t=\phi \varepsilon_{t-1}+\omega_t$. In this case $y_t=\sum_{i=1}^n\alpha_i^*x_{i,t}+\varepsilon_t=\sum_{i=1}^n\alpha_i^*(\phi x_{i,t-1}+u_{i,t})+\phi\varepsilon_{t-1}+\omega_t=\sum_{i=1}^n\alpha_i^*u_{i,t}+\phi y_{t-1} + \omega_t.$ Thus, under the common AR(1) restriction (or common factor restriction,~ \citealp{Mizon1995}), the working model~\eqref{WM} is equivalent to the true model~\eqref{DGP_y}. 
\end{example}

\begin{rem}\label{rem1}
    The working model~\eqref{WM} coincides with the true model~\eqref{DGP_y}
    under a common AR($p$) restriction;
    that is, when $x_{i,t}=\sum_{\ell=1}^p\phi_\ell x_{i,t-\ell}+u_{i,t}$ and $\varepsilon_t=\sum_{\ell=1}^p\phi_\ell\varepsilon_{t-\ell}+\omega_t$.
    In fact, it is easy to show that $y_t=\sum_{i=1}^n\alpha_i^*x_{i,t}+\varepsilon_t=\sum_{i=1}^n\alpha_i^*u_{i,t}+\sum_{\ell=1}^p\phi_\ell y_{t-\ell} + \omega_t$ for any autoregressive order $p$. Moreover, $v_t=\omega_t$ and $E(v_t\mid\mathbf{w}_t)=0$ so the working model is correctly specified, including the response-lag coefficients.  
\end{rem}

\begin{example}\label{ex:NOcommonAR1}
Suppose $x_{i,t}=\phi_i x_{i,t-1}+u_{i,t}$ and $\varepsilon_t=\phi_\varepsilon\varepsilon_{t-1}+\omega_t$, with mutually independent, zero-mean, unit-variance innovations. Then $v_t=\sum_{i=1}^n\alpha_i^*(\phi_i-\phi_y)x_{i,t-1}+(\phi_\varepsilon-\phi_y)\varepsilon_{t-1}+\omega_t$, where $\phi_y=\frac{E(y_ty_{t-1})}{E(y_t^2)}=\left(\sum_{i=1}^n\frac{\phi_i(\alpha_i^*)^2}{1-\phi_i^2}+\frac{\phi_\varepsilon}{1-\phi_\varepsilon^2}\right)\big/\left(\sum_{i=1}^n\frac{(\alpha_i^*)^2}{1-\phi_i^2}+\frac{1}{1-\phi_\varepsilon^2}\right)$, and $E(v_t\mathbf{w}_t)=\pmb0$, whereas $E(v_t\mid\mathbf{x}_{t-1},\varepsilon_{t-1})=\sum_{i=1}^n\alpha_i^*(\phi_i-\phi_y)x_{i,t-1}+(\phi_\varepsilon-\phi_y)\varepsilon_{t-1}$.
\end{example}

\begin{example}\label{ex:FactorModel}
Suppose $x_{i,t}=\lambda_i f_t+z_{i,t}$, $f_t=\phi_f f_{t-1}+\delta_t$, $z_{i,t}=\phi_i z_{i,t-1}+\eta_{i,t}$, and $\varepsilon_t=\phi_{\varepsilon}\varepsilon_{t-1}+\omega_t$, where $0<|\phi_f|,|\phi_i|,|\phi_{\varepsilon}|<1$, $\phi_f\neq\phi_i$, and $\lambda_i\neq0$. The innovations $\delta_t,\eta_{i,t},\omega_t$ are mutually independent i.i.d. $\mathcal{N}(0,1)$. Predictors generally admit marginal ARMA(2,1) representations.
\end{example}

\noindent 
The LASSO estimator of $\pmb\beta^*$ of working model~\eqref{WM} 
is
\begin{equation}\label{wmLASSO}
\widehat{\pmb{\beta}}= \aggregate{argmin}{\pmb{\beta}\in\mathbb{R}^{n+p_y}}\:\left \{\:\frac{1}{2T}\:\left|\left|\mathbf{y}-\mathbf{W}'\pmb{\beta} \right|\right|_2^2 + \lambda||\pmb{\beta}||_1\:\right \}\ \ ,
\end{equation}
where $\lambda>0$ is a tuning parameter. In Section~\ref{sec:feasibleARMAr_LAS}, we address the estimation of the innovations and derive oracle inequalities for ARMAr-LASSO estimation and prediction errors, allowing $n$ to grow with the sample size ($n=n_T$). 

Let $\mathbf{q}_t=(\mathbf{w}_t',v_t)'$. Under Assumption~\ref{ass:CovStat}, for fixed $n$, the causal ARMA representations imply that $\mathbf{q}_t$ is geometrically L$^2$-near-epoch dependent (NED) on the underlying innovation process $(\mathbf{u}_t',\omega_t)'$ (see \citealp[Ch.~17]{davidson1994} and \citealp{ADAMEK2023}). In particular, the variables in Example~\ref{ex:NOcommonAR1} admit linear representations with geometrically decaying coefficients. The same NED property holds in Example~\ref{ex:FactorModel} with respect to the underlying shocks $(\delta_t,\pmb{\eta}_t',\omega_t)'$, although the marginal ARMA innovations need not satisfy the joint independence requirements of Assumption~\ref{ass:CovStat}. Under suitable moment and dependence conditions, NED provides a framework for establishing mixingale properties without requiring the process itself to be mixing (see Supplement~\ref{theo:PrAB_LEMMAS_Proof}).

\subsection{Feasible ARMAr-LASSO}\label{sec:feasibleARMAr_LAS}
Let $\widehat{\mathbf{W}}$ be the feasible design, whose regressors include estimated ARMA innovations $\widehat{\mathbf{w}}_t=(\widehat{u}_{1,t},\dots,\widehat{u}_{n_T,t},y_{t-1},\dots,y_{t-p_y})'$ and response lags.
The use of estimated innovations also changes the working-model residual evaluated at $\pmb{\beta}^*$. Let $\boldsymbol{\Delta}=\widehat{\mathbf W}-\mathbf W$, where only the innovation rows differ. Since $\mathbf y=\mathbf W'\pmb{\beta}^*+\mathbf v$, the corresponding residual for the feasible design is
$
\widehat{\mathbf v}^{\,*}=\mathbf y-\widehat{\mathbf W}'\pmb{\beta}^*=\mathbf v-\boldsymbol{\Delta}'\pmb{\beta}^*=\mathbf v-(\widehat{\mathbf U}-\mathbf U)'\pmb{\alpha}^*.
$
Thus, the relevant ARMAr-LASSO score is $T^{-1}\widehat{\mathbf W}\widehat{\mathbf v}^{\,*}$ rather than $T^{-1}\widehat{\mathbf W}\mathbf v$.

We now address the estimation of ARMA residuals. For each predictor $x_{i,t}$ in~\eqref{DGP_x}, we select the ARMA orders using standard criteria (e.g., BIC) and estimate the parameters using a consistent method, such as Yule--Walker for AR models or maximum likelihood for ARMA models. Under correct specification and standard regularity conditions, the estimated innovations are consistent for each fixed predictor (\citealp{hannan1980,brockwell2016}). However, pointwise consistency does not ensure uniform accuracy as $n=n_T$ grows with $T$. We therefore impose the following assumption.

\begin{ass}
\label{ass:ARMAuniform}
Define
$
\delta_T=\max_{1\leq i\leq n_T}\left(\frac{1}{T}\sum_{t=1}^T(\widehat u_{i,t}-u_{i,t})^2\right)^{1/2}.
$
For some deterministic sequence $r_{T,n_T}>0$ with $r_{T,n_T}\to0$, $\delta_T=O_p(r_{T,n_T})$.
\end{ass}

\begin{lemma}
\label{lem:ARMAscore}
Let Assumption~\ref{ass:ARMAuniform} hold. Define $\|\mathbf a\|_T=(T^{-1}\sum_{t=1}^T a_t^2)^{1/2}$, $M_T=\max_{j\leq n_T+p_y}(T^{-1}\sum_{t=1}^T\widehat w_{j,t}^2)^{1/2}$ and $R_T=\max_{1\leq\ell\leq T}\left\|\frac{1}{T}\sum_{t=1}^{\ell}(\widehat{\mathbf w}_t\widehat v_t^{\,*}-\mathbf w_tv_t)\right\|_\infty$. Then
$
R_T\leq\delta_T\left(\|\mathbf v\|_T+M_T\|\pmb{\alpha}^*\|_1\right).
$
If $\|\mathbf v\|_T=O_p(1)$, $M_T=O_p(1)$, and $r_{T,n_T}(1+\|\pmb{\alpha}^*\|_1)=o(\lambda_T)$, then $R_T=o_p(\lambda_T)$. Consequently, if for some $a>0$,
$
P\left(\max_{1\leq\ell\leq T}\left\|\frac{1}{T}\sum_{t=1}^{\ell}\mathbf w_tv_t\right\|_\infty\leq a\lambda_T\right)\to1,
$
then, for every $b>a$,
$
P\left(\max_{1\leq\ell\leq T}\left\|\frac{1}{T}\sum_{t=1}^{\ell}\widehat{\mathbf w}_t\widehat v_t^{\,*}\right\|_\infty\leq b\lambda_T\right)\to1.
$
\end{lemma}

We derive oracle inequalities for ARMAr-LASSO, allowing $n=n_T$ to grow with $T$, subject to uniform ARMA accuracy, score concentration, estimated-design compatibility, and the rate conditions specified in this section. 
\begin{ass}\label{ass:boundXandQ}
   Consider $\mathbf{q}_t=(\widehat{\mathbf{w}}_t',v_t)'$.
        There exist constants $c_2>c_1>2$ such that $\underset{i\leq n_T+p_y+1,\ t\leq T}{\text{max}}E(\left|q_{i,t}\right|^{2c_2})\leq C$ and $\underset{i\leq n_T,t\leq T}{\text{max}} E(\left| x_{i,t} \right|^{2c_{2}})\leq C$. 

\end{ass}

\begin{ass}
\label{ass:ScoreMix}
Let $\xi_{j,t}=w_{j,t}v_t$, $1\leq j\leq n_T+p_y$, be centered and adapted to a common filtration $(\mathcal{F}_t)$.
Assume that
$
\sup_{T,j,t}\|E(\xi_{j,t}\mid\mathcal{F}_{t-m})\|_{c_1}\leq a_m,\ \sum_{m=0}^{\infty}a_m<\infty,
$
where $c_1>2$ is as in Assumption~\ref{ass:boundXandQ} and $(a_m)$ does not depend on $T$ or $j$.
\end{ass}

Assumption~\ref{ass:boundXandQ} imposes uniform moment bounds and implies $||\mathbf v||_T=O_p(1)$, whereas Assumption~\ref{ass:ScoreMix} requires summable conditional dependence for the centered oracle scores. 

\noindent To bound the ARMAr-LASSO estimation error, we follow the typical procedure in technical textbooks (see, e.g.,~ \citealp[Ch.~6]{Buhlmann2011}). We need $\lambda$ to be sufficiently large as to exceed the empirical process $\underset{i\leq n_T+p_y,t\leq T}{\text{max}}\left| \sum_{t=1}^T\widehat{w}_{i,t}\widehat{v}^*_t \right|$ with high probability.

\begin{theorem}
\label{theo:PrA}
Let Assumptions~\ref{ass:CovStat}-~\ref{ass:ScoreMix} hold. Suppose that $M_T=O_p(1)$, and $r_{T,n_T}(1+\|\pmb{\alpha}^*\|_1)=o(\lambda_T)$. Define $\mathcal{A}_T=\left\{\max_{1\leq j\leq n_T+p_y,\,1\leq\ell\leq T}\left|\sum_{t=1}^{\ell}\widehat w_{j,t}\widehat v_t^{\,*}\right|\leq\frac{T\lambda_T}{4}\right\}$. If, for a sufficiently large constant $C_\lambda>0$, $\lambda_T\geq C_\lambda\left(\sqrt{\log T}\right)^{1/c_1}\frac{(2n_T+p_y)^{1/c_1}}{\sqrt T},$ then, for some constant $C>0$ independent of $T$, $P(\mathcal{A}_T)\geq1-C(\sqrt{\log T})^{-1}-o(1)\longrightarrow1.$
\end{theorem}

\begin{rem}\label{rem:lambdaARMA}
Let $L_T=\frac{(2n_T+p_y)^{1/c_1}(\sqrt{\log T})^{1/c_1}}{\sqrt{T}}$. Choose $\lambda_T=C_\lambda L_T$, with $C_\lambda>0$ sufficiently large. Under the remaining conditions of Lemma~2, $r_{T,n_T}(1+\|\pmb{\alpha}^*\|_1)=o(L_T)$ implies $R_T=o_p(\lambda_T)$.
\end{rem}

\noindent Theorem~\ref{theo:PrA} establishes that the inequalities we need for the error bound of the proposed ARMAr-LASSO estimator hold with high probability (see Remark~\ref{rem:NED} in Supplement~\ref{theo:PrAB_LEMMAS_Proof}). 
We impose the compatibility condition directly on the feasible design $\widehat{\mathbf{W}}$, which allows us to generalize subsequent results to the high-dimensional framework.

\begin{ass}\label{ass:REC} 
Let $\widetilde S=\{j:\beta_j^*\neq0\}$, with $\widetilde s=|\widetilde S|\geq1$, and define the compatibility constant $\gamma_{\widehat w}^2= \underset{\|\pmb{\beta}_{\widetilde S^c}\|_1\leq3\|\pmb{\beta}_{\widetilde S}\|_1;\ \pmb{\beta}\in\mathbb R^{n_T+p_y}\backslash\{0\}}{\text{min}}\frac{\widetilde s\,\pmb{\beta}'\widehat{\mathbf W}\widehat{\mathbf W}'\pmb{\beta}}{T\|\pmb{\beta}_{\widetilde S}\|_1^2}.$ Assume that $\gamma_{\widehat w}^2>0$. This implies that
$\left |\left |\pmb{\beta}_{\tilde{S}}\right |\right |_1^2\leq \tilde{s}\ \frac{\pmb{\beta}'\widehat{\mathbf{W}}\widehat{\mathbf{W}}'\pmb{\beta}}{T\gamma_{\widehat{w}}^2}$ for all $\pmb\beta$ satisfying $||\pmb{\beta}_{{\widetilde S}^c}||_1\leq 3||\pmb{\beta}_{\widetilde{S}}||_1$.
\end{ass}

\begin{rem}
Let $\gamma_w$ denote the compatibility constant of the oracle design $\mathbf W$, defined on the same support $\widetilde S$ and cone as $\gamma_{\widehat w}$. Since $\|(\widehat{\mathbf W}-\mathbf W)'\pmb{\beta}\|_T\leq\delta_T\|\pmb{\beta}\|_1\leq4\delta_T\|\pmb{\beta}_{\widetilde S}\|_1$ for every $\pmb{\beta}$ in the cone, the triangle inequality yields
$
\gamma_{\widehat w}\geq\gamma_w-4\sqrt{\widetilde s}\,\delta_T.
$
Consequently, if $\gamma_w^{-2}=O_p(1)$ and $\sqrt{\widetilde s},r_{T,n_T}\to0$, Assumption~\ref{ass:ARMAuniform} implies $\gamma_{\widehat w}^{-2}=O_p(1)$.\end{rem}

\noindent Theorem~\ref{theo:Bounds} establishes oracle inequalities for the in-sample prediction and $\ell_1$ estimation errors, along with a variable-screening guarantee under a beta-min condition.
\begin{theorem}
\label{theo:Bounds}
Let Assumption~\ref{ass:REC} and the conditions of Theorem~\ref{theo:PrA} hold. Let $\widehat{\pmb{\beta}}=\aggregate{argmin}{\pmb{\beta}\in \mathbb{R}^{n+p_y}}\:\left\{\frac{1}{2T}\|\mathbf y-\widehat{\mathbf W}'\pmb{\beta}\|_2^2+\lambda_T\|\pmb{\beta}\|_1\right\}$, $\widetilde S=\{j:\beta_j^*\neq0\}$, and put $\widetilde s=|\widetilde S|\geq1$. Then, simultaneously with probability at least $1-C(\sqrt{\log T})^{-1}-o(1)$,
$
\text{(a)}\|\widehat{\pmb{\beta}}-\pmb{\beta}^*\|_1\leq\frac{C\widetilde s\lambda_T}{\gamma_{\widehat w}^2},\
\text{(b)}\frac1T\|\widehat{\mathbf W}'(\widehat{\pmb{\beta}}-\pmb{\beta}^*)\|_2^2\leq\frac{C\widetilde s\lambda_T^2}{\gamma_{\widehat w}^2}. 
$
On the same event, if $\min_{j\in\widetilde S}|\beta_j^*|>C\widetilde s\lambda_T/\gamma_{\widehat w}^2$, the variable-screening property holds.
\end{theorem}

\begin{cor}
\label{rem4} 
Under the conditions of Theorem~\ref{theo:Bounds}, suppose that
$\lambda_T\asymp\frac{(\sqrt{\log T})^{1/c_1}(2n_T+p_y)^{1/c_1}}{\sqrt{T}},\ \gamma_{\widehat w}^{-2}=O_p(1).
$
Then the following rates hold:
$\text{(a)}\left\|\widehat{\pmb{\beta}}-\pmb{\beta}^*\right\|_1
=O_p\left(\frac{\widetilde s}{\sqrt{T}}\left[(n_T+p_y)\sqrt{\log T}\right]^{1/c_1}\right), 
\text{(b)}\frac{1}{T}\left\|\widehat{\mathbf W}'(\widehat{\pmb{\beta}}-\pmb{\beta}^*)\right\|_2^2
=O_p\left(\frac{\widetilde s}{T}\left[(n_T+p_y)\sqrt{\log T}\right]^{2/c_1}\right).$ If $\widetilde s\lambda_T\to0$, both quantities converge to zero in probability.
\end{cor}

\section{Simulations and Empirical Application}\label{sec:Exp}
We evaluate ARMAr-LASSO through simulations and a real-data application.

\subsection{Simulation Experiments}\label{ARMArLASSOMonteCarlo}
The response variable is generated using the model $y_t=\sum_{i=1}^n\alpha_i^*x_{i,t-1}+\varepsilon_t$, where predictors are lagged to reflect their availability at the forecast origin. We consider the following data generating processes (DGPs) for predictors and error terms:
\setcounter{bean}{0}
\begin{list}
{(\Alph{bean})}{\usecounter{bean}}   
\item Common AR(1) Restriction: $x_{i,t}=\phi_i x_{i,t-1}+u_{i,t},\ \varepsilon_t=\phi_{\varepsilon}\varepsilon_{t-1}+\omega_t$.
\item General AR/ARMA: $x_{i,t}= qf_t+z_{i,t}$, where $f_t=0.9 f_{t-1}+\delta_t$, ${z}_{j,t}=0.8{z}_{j,t-1}+\eta_{j,t}$; ${z}_{h,t}=0.6{z}_{h,t-1}+0.3{z}_{h,t-2}+\eta_{h,t}$; $z_{w,t}=0.5z_{w,t-1}+0.4z_{w,t-2}+\eta_{w,t}+0.3\eta_{w,t-1}$; $z_{k,t}=0.7z_{k,t-1}+\eta_{k,t}+0.4\eta_{k,t-1}$, for $t=1,\dots,T$, and $j=1,\dots,4;\ h=5,\dots,7;\ w=8,\dots,10;\ k=11,\dots,n$. The error terms are generated as $\varepsilon_t=0.7\varepsilon_{t-1}+0.2\varepsilon_{t-2}+\omega_t$. 
\end{list}

\noindent The shocks are generated as follows: $u_{i,t} \overset{\mathrm{i.i.d.}}{\sim} \mathcal{N}(0,1)$ with $\left(\mathbf{C}_{u}\right){ij}=c_{ij}^{u}=\rho^{|i-j|}$, $\delta_t,\ \eta_{i,t} \overset{\mathrm{i.i.d.}}{\sim} \mathcal{N}(0,1)$ with $\left(\mathbf{C}_{\eta}\right){ij}=c_{ij}^{\eta}=\rho^{|i-j|}$, and $\omega_t \overset{\mathrm{i.i.d.}}{\sim} \mathcal{N}(0,\sigma_{\omega}^2)$. For the DGP (A) and for the DGP (B) with $q=0$ we set $\rho=0.8$, while for the DGP (B) with $q=1$ we set $\rho=0.4$ to generate predictors primarily influenced by the common factor, with weakly correlated AR or ARMA idiosyncratic components. Finally, we vary the value of $\sigma_{\omega}^2$ to explore different signal-to-noise ratios (SNRs). For each DGP, the performance of ARMAr-LASSO and benchmarks is evaluated based on average results from 1000 independent simulations, focusing on the coefficient estimation error (CoEr) obtained as $||\widehat{\pmb\alpha}-\pmb\alpha||_2$, the root mean square forecast error (RMSFE), and the percentages of true positives (\%TP) and false positives (\%FP) in selecting relevant predictors. For all methods, the tuning parameter $\lambda$ is selected using the Bayesian Information Criterion (BIC). Regardless of the choice of $n$, $\pmb\alpha^*$ always has the first $10$ entries equal to $1$ and all others equal to $0$. In Supplement~\ref{sec:FiltersExp}, we extend the analysis of the ARMA filtering procedure by examining the effects of misspecifying the autoregressive orders of the predictors, using alternatives to BIC for model selection and including fewer lags of $y_t$.

\subsubsection{DGP(A)}\label{sec:AR1Filter}
We consider two complementary AR(1) settings with $p_y=1$. Under DGP (A1), $\phi_i \overset{\mathrm{i.i.d.}}{\sim} U(0.6,0.95)$ and $\phi_{\varepsilon}=0.8$. Heterogeneous autoregressive coefficients imply that a single response lag generally cannot fully capture the dynamics of $y_t$. Under DGP (A2), $\phi_i=\phi_{\varepsilon}=\phi$, yielding the exact representation $y_t = \mathbf{u}_{t -1}'\pmb{\alpha}^* + \phi y_{t-1} + \omega_t$. Thus, the two settings allow us to evaluate ARMAr-LASSO both with and without an exact dynamic representation.

\begin{figure}
\graphicspath{{images/}}
\centering
{\includegraphics[width=12cm]{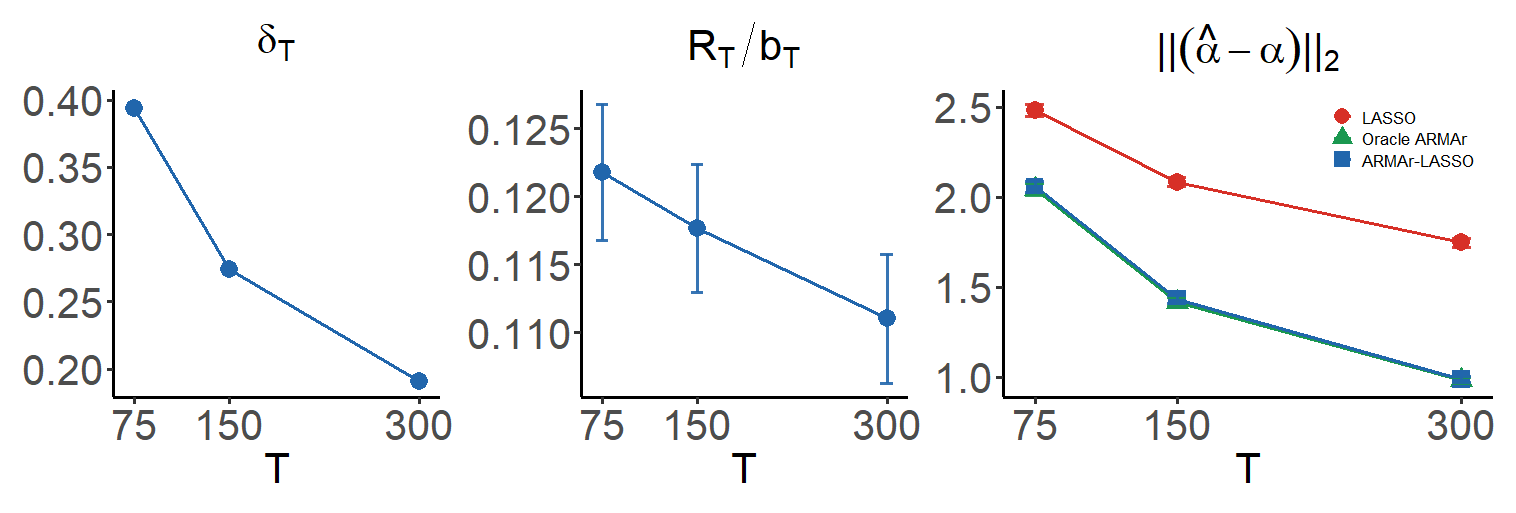}}
\caption{\footnotesize DGPs (A1). $\delta_T$, $R_T/b_T$, and $\ell_2$-estimation errors. Red circles, blue squares, and green triangles denote LASSO, ARMAr-LASSO, and Oracle ARMAr-LASSO, respectively. Points and bars represent means and 95\% confidence intervals from 1000 Monte Carlo replications.
}\label{fig:score_perturbation}
\end{figure}

\noindent{\bf DGP(A1): Filter Evaluation.} We fix $n=150$, SNR=10, and consider sample sizes $T\in\{75,150,300\}$.
For ARMAr-LASSO, innovations are estimated by fitting AR(1) models to each predictor using the Yule–Walker method. The population projection coefficient $\phi_y$ depends on the persistence and relative variances of the signal and error components. We compare ARMAr-LASSO with standard LASSO applied to the observed time series, and the oracle version of ARMAr-LASSO described in~\eqref{wmLASSO}.
Figure~\ref{fig:score_perturbation} reports the Monte Carlo averages of the uniform innovation estimation error $\delta_T$, the ratio $R_T/b_T$, where $b_T=\delta_T(\|\mathbf v\|_T+M_T\|\pmb{\alpha}^*\|_1)$, and the $\ell_2$-estimation errors of the three estimators. As $T$ increases, $\delta_T$ decreases, while $R_T/b_T$ remains well below one, consistent with Lemma~\ref{lem:ARMAscore}. Moreover, ARMAr-LASSO achieves smaller estimation errors than conventional LASSO and performs remarkably close to its oracle counterpart. These findings suggest that the additional error from estimating the innovations is limited, even in relatively small samples.

\noindent{\bf DGP(A2): Impact of $\phi$.} We compare our ARMAr-LASSO (ARMAr-LAS) with the standard LASSO (LAS). Simulations have varying numbers of predictors (dimensionality), $n\in\{75, 150, 300\}$, and a fixed sample size, $T=150$. In this way we cover low ($n=75$), intermediate ($n=150$), and high ($n=300$) dimensional scenarios. As $n$ varies, we also cover different levels of sparsity consequently to $||\pmb{\alpha}^*||_0=10$. We investigate settings with different $\phi\in\{0.3, 0.6, 0.9, 0.95\}$ and SNR$\in\{0.5,1,5,10\}$. 
Results are presented in Table~\ref{TabDGPA&B} for SNR=10 (complete results are provided in Supplement~\ref{ARMArLASSOMonteCarlo_supp}). For each SNR, CoEr and RMSFE (both expressed in relative terms to the values obtained by LAS), as well as \%TP, and \%FP are given for every $n$ and $\phi$ considered. 
Results have ARMAr-LAS as the best performer in terms of CoEr and RMSFE across values of $\phi$ and $n$ -- demonstrating superior accuracy in both estimation and forecasting compared to LASSO. ARMAr-LAS improves TP selection and sharply reduces FP under stronger persistence ($\phi\geq0.6$).

\subsubsection{DGP(B): General AR/ARMA.}\label{sec:AR/ARMA}
We compare our ARMAr-LASSO (ARMAr-LAS) with the standard LASSO (LAS), standard LASSO plus lags of $y_t$ (LASy), GLS-LASSO as proposed by~\cite{chronopoulos2023} (GLS-LAS), autoregressive distributed lag LASSO (ARDL-LAS), and FarmSelect as proposed by~\cite{Fan2020}, which employs LASSO on factor model residuals (FaSel). For GLS-LAS, we filter both response and predictors using the coefficients of an AR($p_{\varepsilon}$) model applied to $\widehat{\varepsilon}_t$, with the order $p_{\varepsilon}$ (max 2) selected with BIC, with parameter estimation by maximum likelihood. For ARDL-LAS, we consider two lags of the response and two lags of each predictor as additional regressors -- bringing the number of term undergoing selection to $n\times3+2$. For ARMAr-LAS, the $\widehat{u}$'s are obtained by filtering each series with an ARMA($p_i,q_i$) process, with the orders $p_i$ and $q_i$ (max 2) selected via BIC. We consider $p_y=3$; that is, three lags of $y_t$ as additional predictors. Simulations have varying $n\in\{75, 150, 300\}$, SNR$\in\{1,10\}$, $q\in\{0,1\}$ and a fixed sample size, $T=150$. The common-factor design allows us to assess ARMAr-LASSO beyond the innovation-independence conditions of Assumption~\ref{ass:CovStat}, examining whether ARMA filtering can still mitigate serial dependence and improve performance relative to competing methods. Results are presented in Table~\ref{TabC&D}. ARMAr-LAS outperforms all other LASSO-based methods in terms of estimation accuracy, forecasting, and feature selection, regardless of the SNR and the presence of a common factor (i.e., when $q=1$). In these DGPs, comparisons with GLS-LASSO and ARDL-LASSO are particularly informative. Unlike error-based GLS filtering, ARMAr-LASSO targets each predictor’s dynamics, reducing persistence and yielding more accurate estimates and forecasts than GLS-LASSO. Moreover, ARMAr-LASSO includes only three lags of $y_t$, making it more parsimonious than ARDL-LASSO, which becomes overparameterized as predictor lags increase. These results support ARMAr-LASSO under heterogeneous ARMA dynamics and common factors, beyond the common-AR case.

To conclude, we summarize key findings from additional simulations reported in the Supplement due to space constraints. Supplement~\ref{sec:FiltersExp} examines varying numbers of $y_t$ lags ($p_y$), AR filter orders ($p_i$), and filter selection methods. Supplements~\ref{ARMArLASSOMonteCarlo_supp},~\ref{sec:largeT} and~\ref{sec:misspecifiedAR} extend the analysis to varying AR(1) coefficients $\phi$, larger sample sizes $T$, and misspecified predictor AR models. In particular, Supplement~\ref{ARMArLASSOMonteCarlo_supp} shows that ARMAr-LASSO outperforms the LASSO-based benchmarks of Section~\ref{sec:AR/ARMA} across all four performance metrics as $\phi$ increases. Finally, Supplement~\ref{Sec:MinEig} reports larger minimum eigenvalues for ARMAr-LASSO than for LASSO and GLS-LASSO in the examined low-dimensional settings. 
\begin{table}
\centering
\caption{\footnotesize DGPs (A2). CoEr, RMSFE (relative to LAS), \%TP and \%FP for LASSO and ARMAr-LASSO.}
\label{TabDGPA&B}
\makebox[\linewidth]{
\scalebox{0.5}{
\begin{tabular}{
@{}l@{\hspace{0.3cm}}rrrr
@{\hspace{0.6cm}}rrrr
@{\hspace{0.6cm}}rrrr@{}
}
\hline
$n$ & \multicolumn{4}{c}{75}
    & \multicolumn{4}{c}{150}
    & \multicolumn{4}{c}{300}\\\hline
$\phi$
    & 0.3 & 0.6 & 0.9 & 0.95
    & 0.3 & 0.6 & 0.9 & 0.95
    & 0.3 & 0.6 & 0.9 & 0.95\\\hline

CoEr &&&&&&&&&&&&\\\cline{1-1}
ARMAr-LAS
    & 0.94 & 0.77 & 0.51 & 0.50
    & 0.95 & 0.77 & 0.59 & 0.57
    & 0.95 & 0.79 & 0.63 & 0.60\\\hline

RMSFE &&&&&&&&&&&&\\\cline{1-1}
ARMAr-LAS
    & 0.97 & 0.83 & 0.70 & 0.69
    & 0.99 & 0.83 & 0.70 & 0.65
    & 0.98 & 0.84 & 0.74 & 0.69\\\hline

\% TP &&&&&&&&&&&&\\\cline{1-1}
LASSO
    & 97.20 & 93.10 & 88.00 & 87.70
    & 97.10 & 93.60 & 87.30 & 85.30
    & 96.70 & 93.00 & 88.30 & 84.90\\
ARMAr-LAS
    & 98.20 & 97.80 & 97.20 & 96.60
    & 97.70 & 97.90 & 96.90 & 96.30
    & 97.50 & 97.60 & 96.60 & 95.90\\\hline

\% FP &&&&&&&&&&&&\\\cline{1-1}
LASSO
    & 2.60 & 7.10 & 30.20 & 31.60
    & 1.30 & 4.30 & 13.70 & 12.30
    & 0.70 & 2.60 & 8.00 & 7.30\\
ARMAr-LAS
    & 2.90 & 3.20 & 3.30 & 3.10
    & 1.40 & 1.50 & 1.50 & 1.60
    & 0.70 & 0.80 & 0.80 & 0.80\\
\hline
\end{tabular}
}
}
\end{table}
\begin{table}[t]
\centering
\caption{\footnotesize DGP (B). CoEr, RMSFE (relative to LAS), \%TP, and \%FP for LASSO-based benchmarks and ARMAr-LASSO. For each configuration, the best CoEr and RMSFE are in bold.}
\label{TabC&D}
\begingroup
\setlength{\tabcolsep}{2pt}
\renewcommand{\arraystretch}{1.03}
\resizebox{\textwidth}{!}{%
\begin{tabular}{@{}ccc!{\hspace{20pt}}rrrrr!{\hspace{20pt}}rrrrr!{\hspace{20pt}}rrrrrr!{\hspace{20pt}}rrrrrr@{}}
\toprule
\multicolumn{3}{c}{DGP settings} & \multicolumn{5}{c}{CoEr} & \multicolumn{5}{c}{RMSFE} & \multicolumn{6}{c}{\%TP} & \multicolumn{6}{c}{\%FP} \\
\cmidrule(lr){4-8}\cmidrule(lr){9-13}\cmidrule(lr){14-19}\cmidrule(l){20-25}
$q$ & SNR & $n$ & \rotatebox{65}{LASSOy} & \rotatebox{65}{GLS-LAS} & \rotatebox{65}{ARDL-LAS} & \rotatebox{65}{FaSel} & \rotatebox{65}{ARMAr-LAS} & \rotatebox{65}{LASSOy} & \rotatebox{65}{GLS-LAS} & \rotatebox{65}{ARDL-LAS} & \rotatebox{65}{FaSel} & \rotatebox{65}{ARMAr-LAS} & \rotatebox{65}{LASSO} & \rotatebox{65}{LASSOy} & \rotatebox{65}{GLS-LAS} & \rotatebox{65}{ARDL-LAS} & \rotatebox{65}{FaSel} & \rotatebox{65}{ARMAr-LAS} & \rotatebox{65}{LASSO} & \rotatebox{65}{LASSOy} & \rotatebox{65}{GLS-LAS} & \rotatebox{65}{ARDL-LAS} & \rotatebox{65}{FaSel} & \rotatebox{65}{ARMAr-LAS} \\
\midrule
0 & 1 & 75 & 0.49 & 0.65 & 0.45 & 1.04 & \textbf{0.41} & 0.80 & 0.79 & 0.80 & 0.96 & \textbf{0.67} & 57.10 & 46.50 & 60.10 & 46.70 & 47.10 & 65.60 & 32.30 & 10.90 & 16.00 & 4.10 & 29.80 & 4.10 \\
0 & 1 & 150 & 0.59 & 0.87 & 0.54 & 1.03 & \textbf{0.50} & 0.88 & 0.88 & 0.89 & 0.95 & \textbf{0.74} & 52.00 & 43.80 & 56.90 & 42.20 & 45.60 & 62.40 & 18.40 & 6.50 & 13.90 & 2.30 & 17.70 & 2.00 \\
0 & 1 & 300 & 0.66 & 0.93 & 0.61 & 1.06 & \textbf{0.55} & 0.90 & 0.93 & 0.90 & 0.96 & \textbf{0.79} & 52.60 & 44.10 & 55.00 & 42.20 & 49.70 & 60.80 & 11.30 & 4.40 & 8.90 & 1.70 & 11.60 & 0.90 \\
\addlinespace[3pt]
0 & 10 & 75 & 0.92 & 0.77 & 0.81 & 1.17 & \textbf{0.53} & 0.97 & 0.79 & 1.00 & 0.96 & \textbf{0.71} & 87.60 & 87.20 & 93.30 & 85.90 & 73.60 & 97.20 & 29.50 & 27.60 & 18.50 & 8.60 & 32.20 & 5.10 \\
0 & 10 & 150 & 0.95 & 0.90 & 0.89 & 1.07 & \textbf{0.60} & 0.97 & 0.88 & 0.99 & 0.94 & \textbf{0.75} & 87.00 & 86.50 & 91.40 & 86.50 & 81.80 & 96.60 & 14.60 & 13.50 & 10.70 & 5.70 & 16.50 & 2.40 \\
0 & 10 & 300 & 0.97 & 0.96 & 0.91 & 1.01 & \textbf{0.64} & 0.98 & 0.94 & 0.97 & 0.94 & \textbf{0.80} & 87.50 & 87.40 & 91.40 & 87.20 & 87.10 & 96.40 & 9.10 & 8.70 & 7.80 & 3.60 & 9.60 & 1.20 \\
\addlinespace[3pt]
1 & 1 & 75 & 0.50 & 0.64 & 0.46 & 1.02 & \textbf{0.43} & 0.81 & 0.77 & 0.78 & 0.99 & \textbf{0.70} & 48.10 & 31.90 & 52.80 & 35.00 & 45.60 & 54.80 & 31.80 & 12.20 & 19.60 & 5.20 & 32.30 & 10.50 \\
1 & 1 & 150 & 0.56 & 0.88 & 0.52 & 1.04 & \textbf{0.48} & 0.83 & 0.93 & 0.80 & 0.97 & \textbf{0.74} & 39.40 & 27.60 & 43.30 & 29.70 & 37.90 & 50.00 & 19.70 & 7.80 & 16.40 & 3.10 & 21.10 & 6.20 \\
1 & 1 & 300 & 0.63 & 0.93 & 0.58 & 1.04 & \textbf{0.53} & 0.87 & 0.94 & 0.85 & 0.98 & \textbf{0.79} & 33.40 & 22.90 & 35.60 & 24.00 & 32.10 & 43.40 & 11.30 & 5.00 & 10.00 & 1.80 & 12.70 & 3.60 \\
\addlinespace[3pt]
1 & 10 & 75 & 0.93 & 0.77 & 0.80 & 1.06 & \textbf{0.52} & 0.94 & 0.79 & 0.92 & 0.98 & \textbf{0.70} & 82.10 & 81.70 & 91.60 & 81.10 & 78.60 & 96.40 & 31.50 & 29.50 & 25.40 & 8.70 & 35.90 & 13.30 \\
1 & 10 & 150 & 0.95 & 0.92 & 0.87 & 1.06 & \textbf{0.59} & 0.97 & 0.90 & 0.93 & 0.94 & \textbf{0.74} & 79.30 & 78.40 & 85.70 & 78.60 & 79.00 & 95.30 & 16.50 & 15.30 & 14.50 & 5.40 & 22.40 & 8.30 \\
1 & 10 & 300 & 0.97 & 0.96 & 0.91 & 1.03 & \textbf{0.62} & 0.98 & 0.94 & 0.98 & 0.96 & \textbf{0.79} & 76.60 & 76.00 & 83.20 & 75.20 & 78.20 & 94.30 & 9.70 & 9.20 & 9.30 & 3.20 & 13.50 & 5.10 \\
\bottomrule
\end{tabular}%
}
\endgroup
\end{table}

\smallskip

\subsection{Empirical Application}\label{EmpApp}
We consider Euro Area (EA) data composed of 309 monthly macroeconomic time series spanning the period between January 1997 and December 2018. 
We use a single revised-data vintage, so the exercise evaluates pseudo-out-of-sample rather than real-time
forecasting. The series are listed in Supplement~\ref{Stat_EA}, grouped according to their measurement domain: Industry \& Construction Survey (ICS), Consumer Confidence Indicators (CCI), Money \& Interest Rates (M\&IR), Industrial Production (IP), Harmonized Consumer Price Index (HCPI), Producer Price Index (PPI), Turnover \& Retail Sale (TO), Harmonized Unemployment Rate (HUR), and Service Surveys (SI). Supplement~\ref{Stat_EA} also reports transformations applied to the series to achieve stationarity (we did not attempt to identify or remove outliers), as well as 
an analysis of the autocorrelation 
functions that justifies
the use of our ARMAr-LASSO in this context. The target variable is the Overall EA Consumer Price Index (CPI), which is transformed as I(2) (i.e.~integration of order 2) following~\cite{SeW2002b}: $y_{t+h}=(1200/h)\text{log}(CPI_{t+h}/CPI_t)-1200\ \text{log}(CPI_t/CPI_{t-1})$, 
where $y_t=1200\ \text{log}(CPI_t/CPI_{t-1})-1200\ \text{log}(CPI_{t-1}/CPI_{t-2})$, and $h$ is the forecasting horizon. We compute forecasts of $y_{t+h}$ at horizons $h=12$ and $24$ using a rolling-window scheme with a fixed span of 120 months. At each forecast origin $t$, models are re-estimated using observations from $t-120$ to $t-h$, accounting for the horizon-specific availability of the response variable. The window advances by one month at each iteration. The last forecast is  December 2018. At each forecast origin $t$, model estimation uses only responses available through $t-h$, while predictor filtering uses information available through $t$, ensuring a pseudo-out-of-sample evaluation. The methods employed for  our empirical exercise are:
\setcounter{bean}{0}
\begin{list}
{(\alph{bean})}{\usecounter{bean}}
\item {\it Univariate AR($p$):} the autoregressive forecasting model based on $p$ lagged values of the target variable, i.e.~$\widehat{y}_{t+h}=\widehat{\alpha}_0+\sum_{\ell=1}^p\widehat{\phi}_\ell y_{t-\ell+1}$, which serves as a benchmark. 

\item {\it LAS} (\citealp{tibshirani96}): forecasts are obtained from the equation $\widehat{y}_{t+h}=\widehat{\alpha}_0+\sum_{\ell=0}^{11}\widehat{\phi}_{\ell}y_{t-\ell}+\sum_{i=1}^{308}\widehat{\alpha}_ix_{it},$ where $(\widehat{\phi}_{0},\dots, \widehat{\phi}_{11}, \widehat{\alpha}_{1},\dots,\widehat{\alpha}_{308})$ is the sparse vector of penalized regression coefficients estimated by the LASSO.

\item{\it GLS-LAS} (\citealp{chronopoulos2023}): forecasts are obtained from the equation
$\widehat{y}_{t+h}=\widehat{\alpha}_0+\sum_{\ell=1}^{p_{\varepsilon}}\widehat{\phi}_{\ell}y_{t-\ell+1}+\sum_{i=1}^{308}\widehat{\alpha}_i\widetilde{x}_{it},$
where $(\widehat{\alpha}_1,\dots,\widehat{\alpha}_{308})$ is the sparse vector of penalized regression coefficients estimated by the LASSO using pre-filtered response and 
predictors (the $\widetilde{x}$'s) as detailed in Supplement~\ref{COMP_ARDLandGLS}.

\item{\it ARDL-LAS}: forecasts are obtained from the equation
$\widehat{y}_{t+h}=\widehat{\alpha}_0+\sum_{\ell=0}^{11}\widehat{\phi}_{\ell}y_{t-\ell}+\sum_{i=1}^{308}\sum_{j=0}^2\widehat{\alpha}_{i,t-j}x_{i,t-j},$ where $(\widehat{\phi}_{0},\dots, \widehat{\phi}_{11}, \widehat{\alpha}_{1,t},\dots,\widehat{\alpha}_{308,t-2})$ is the sparse vector of penalized regression coefficients estimated by the LASSO, which in this case contains two lagged values for each predictor. 

\item {\it FaSel} (\citealp{Fan2020}): FarmSelector applies feature selection on factors residuals. 
Forecasts are obtained from the equation:
$\widehat{y}_{t+h}=\widehat{\alpha}_0+\widehat{\pmb{\Lambda}}\widehat{\pmb{f}}_t+\widehat{\pmb{\alpha}}'\widehat{\mathbf{z}}_t+\sum_{\ell=1}^p\widehat{\phi}_\ell y_{t-\ell+1},$ where $\widehat{\pmb{f}}_t$ is a $r$-dimensional vector of factors estimated with PCA (as in~\cite{SeW2002a, SeW2002b}), $\widehat{\mathbf{z}}_t=\widehat{\pmb{\Lambda}}\widehat{\pmb{f}}_t-\mathbf{x}_t$, $\widehat{\pmb{\Lambda}}'$ is the $n\times r$ matrix of loadings, and $\widehat{\pmb{\alpha}}$ is the sparse vector obtained applying the LASSO. The number of factors $r$ is chosen with the approach described in~\cite{AhnHorenstein2013}.

\item {\it ARMAr-LAS}: our proposal, where LASSO is applied to the estimated ARMA residuals. Forecasts are obtained from the equation $\widehat{y}_{t+h}=\widehat{\alpha}_0+\sum_{\ell=0}^{11}\widehat{\phi}_{\ell}y_{t-\ell}+\sum_{i=1}^{308}\widehat{\alpha}_i\widehat{u}_{it},$ where $(\widehat{\phi}_{0},\dots, \widehat{\phi}_{11}, \widehat{\alpha}_{1},\dots,\widehat{\alpha}_{308})$ is the sparse vector of penalized regression coefficients estimated by the LASSO.
\end{list}

\noindent For the AR($p$) benchmark and the GLS-LAS, the lag orders $p$ and $p_{\varepsilon}$ are selected by BIC within $0\leq p,p_{\varepsilon}\leq 12$.
For the ARMAr-LAS, estimated residuals (the $\widehat{u}$'s) are obtained filtering each time series with an ARMA($p_i,q_i$), where $p_i$ and $q_i$ are selected by BIC within $0 \leq p_i,q_i\leq 12$, $i=1,\dots,n$. For all the LASSO-based methods (including our ARMAr-LAS), the shrinkage parameter $\lambda$ is also selected by BIC. Forecasting accuracy is evaluated using the root mean square forecast error (RMSFE), defined as
 $RMSFE = \sqrt{\frac{1}{T_1 - T_0}\sum_{\tau=T_0}^{T_1}\Bigl( \widehat{y}_{\tau} - y_{\tau}\Bigr)^2},$
where $T_0$ and $T_1$ are the first and last time points used for the out-of-sample evaluation. We also consider the number of selected variables. 

Table~\ref{Tab:RMSFE} reports RMSFE ratios and Diebold–Mariano comparisons of origin-aligned loss differentials. Furthermore, the column ``Selected Variables (Av.)'' reports the average number of selected variables with ARMAr-LAS (on the left side of the vertical bar), and other LASSO-based
methods (on the right side of the vertical bar). ARMAr-LAS improves RMSFE over all reported benchmarks; the AR($p$) difference is significant at $h=24$. This is, in principle, consistent with the theoretical analysis we provided earlier. The sparser ARMAr-LAS output may be due to fewer false positives, as compared to other 
methods 
which suffer from the effects of spurious correlations induced by serial correlation. Notably, GLS-LAS selects fewer 
predictors than ARMAr-LAS
but provides significantly worse forecasts. However, since in this real data application we do not know the true DGP, any comment regarding accuracy in feature selection is necessarily speculative.
\begin{table}
\centering
\caption{\footnotesize 
RMSFE (ratio): ratios of RMSFE contrasting pairs of employed methods; for each ratio, we perform a Diebold-Mariano test (alternative: the second method is less accurate in forecasting) and report p-values as 0 '***'  0.001 '**' 0.01 '*' 0.05 '$^{\sbullet}$' 0.1''. Selected Variables (Av.): average of the number of variables selected by ARMAr-LAS (left side of the vertical bar) and LASSO-based benchmarks (right side).
}\label{Tab:RMSFE}
\scalebox{0.5}{
\begin{tabular}{@{}clllllllllll@{}}\hline
&Method 1 & & Method 2 & & \multicolumn{3}{c}{RMSFE (ratio)} && \multicolumn{3}{c}{Selected Variables (Av.)}  \\\cline{2-2} \cline{4-4} \cline{6-8} \cline{10-12}
&& &  & & $h$=12 &&$h=24$ && $h$=12 && $h=24$ \\\cline{6-6} \cline{8-8} \cline{10-10} \cline{12-12} 
&ARMAr-LAS && LAS        & & 0.68***  && 0.89*                 &&6.2$|$65.8  &&11.8$|$62.2 \\
&ARMAr-LAS && GLS-LAS    & & 0.69***  && 0.60***               &&6.2$|$2.5   &&11.8$|$5.6 \\
&ARMAr-LAS && ARDL-LAS   & & 0.61***  && 0.86$^{\sbullet}$     &&6.2$|$36.6  &&11.8$|$32.9 \\
&ARMAr-LAS && FaSel      & & 0.72***  && 0.83**                &&6.2$|$77.3  &&11.8$|$72.4 \\
&ARMAr-LAS && AR($p$)    & & 0.95     && 0.89* &&--&&--\\
&LAS && AR($p$)          & & 1.39     && 1.08 &&--&&-- \\
&GLS-LAS && AR($p$)      & & 1.38     && 1.44 &&--&&-- \\
&ARDL-LAS && AR($p$)     & & 1.56     && 1.07 &&--&&-- \\
&FaSel && AR($p$)        & & 1.33     && 1.21 &&--&&-- \\
\hline
\end{tabular}
}
\end{table}

\begin{figure}
\graphicspath{{images/}}
\centering
\subfloat[\footnotesize LAS, $h$=12]{\includegraphics[width=4cm]{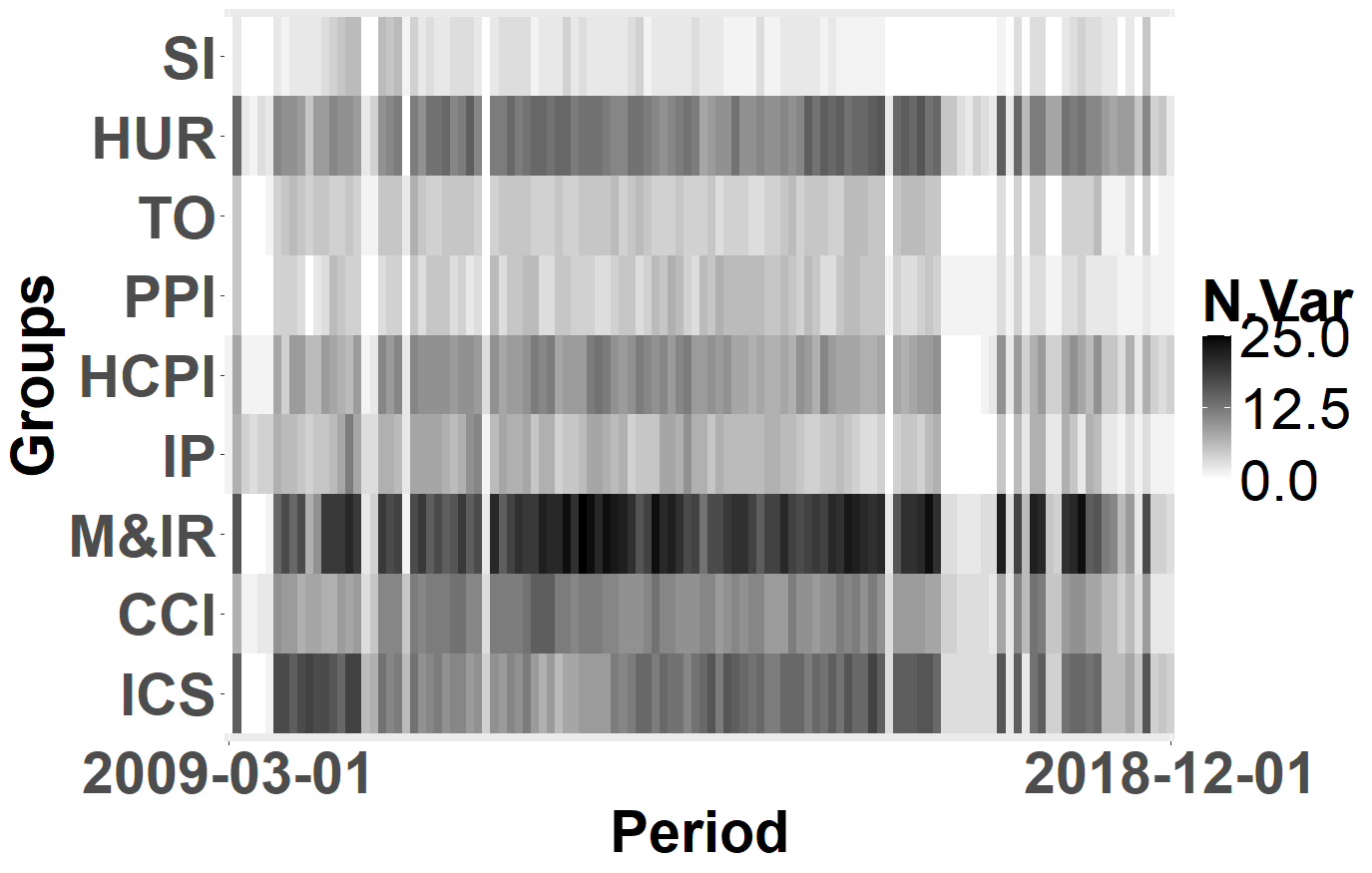}}\hfil
\subfloat[\footnotesize ARMAr-LAS, $h$=12]{\includegraphics[width=4cm]{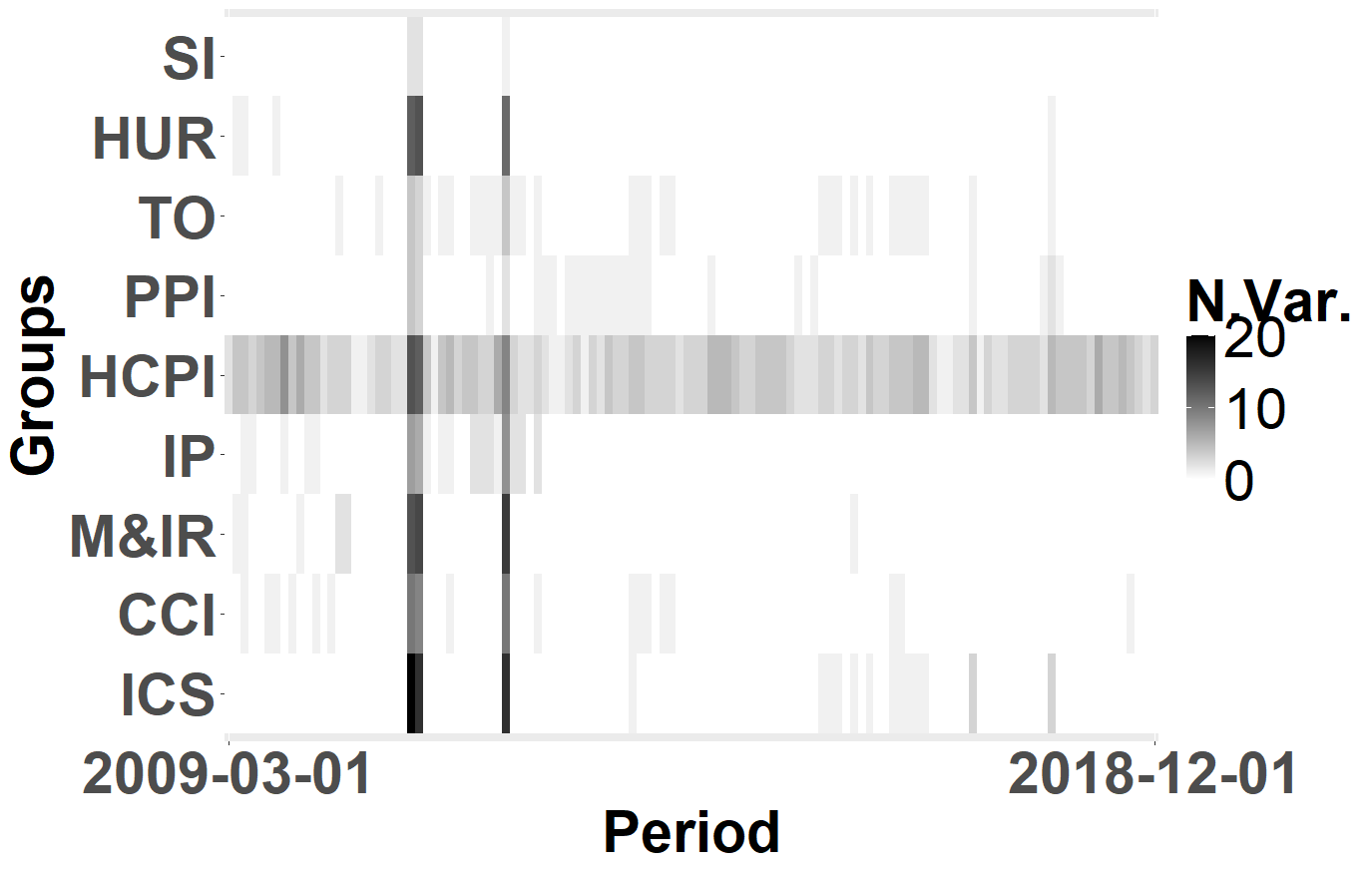}}\hfil
\subfloat[\footnotesize LAS, $h$=24]{\includegraphics[width=4cm]{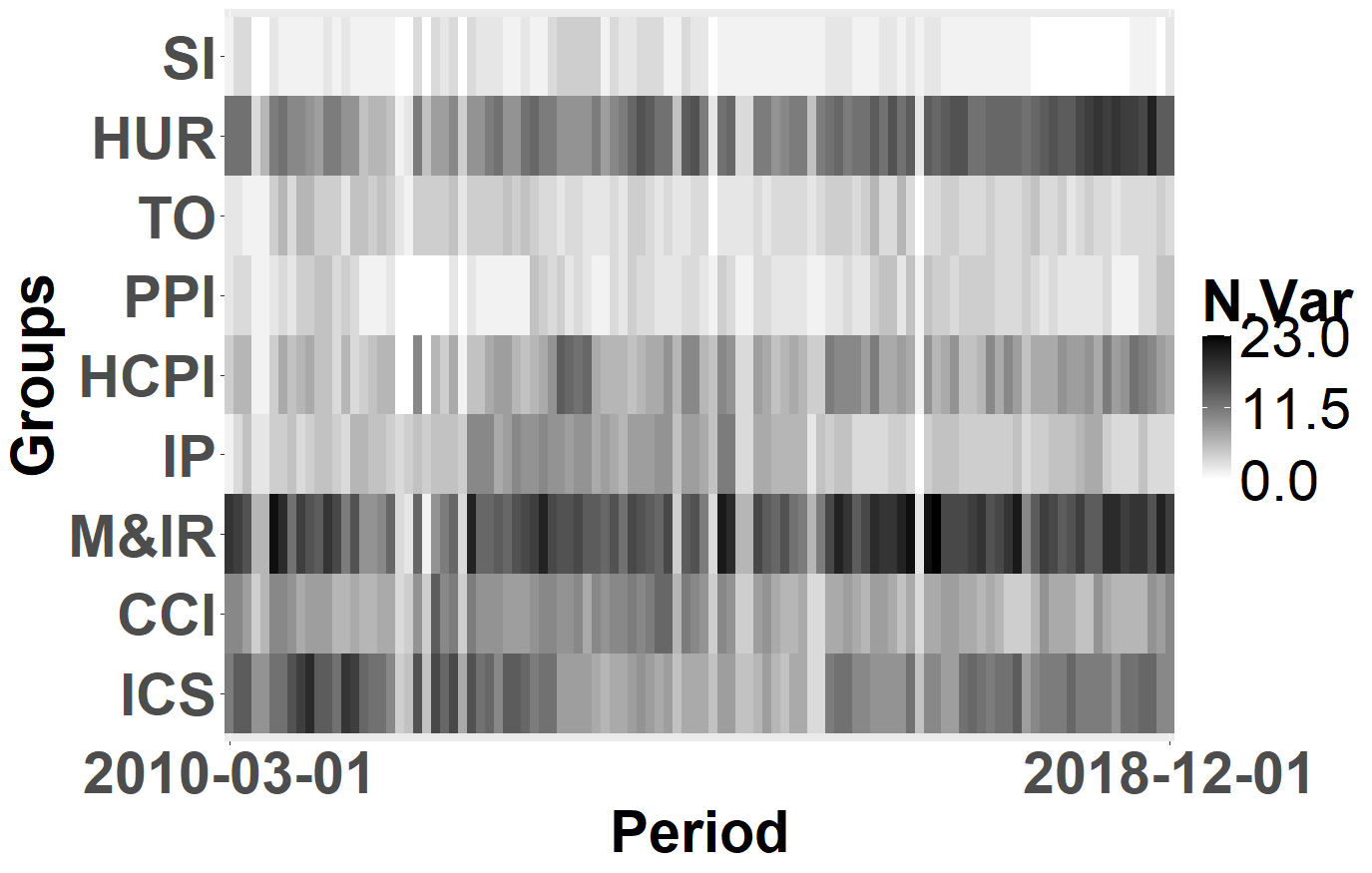}}\hfil
\subfloat[\footnotesize ARMAr-LAS, $h$=24]{\includegraphics[width=4cm]{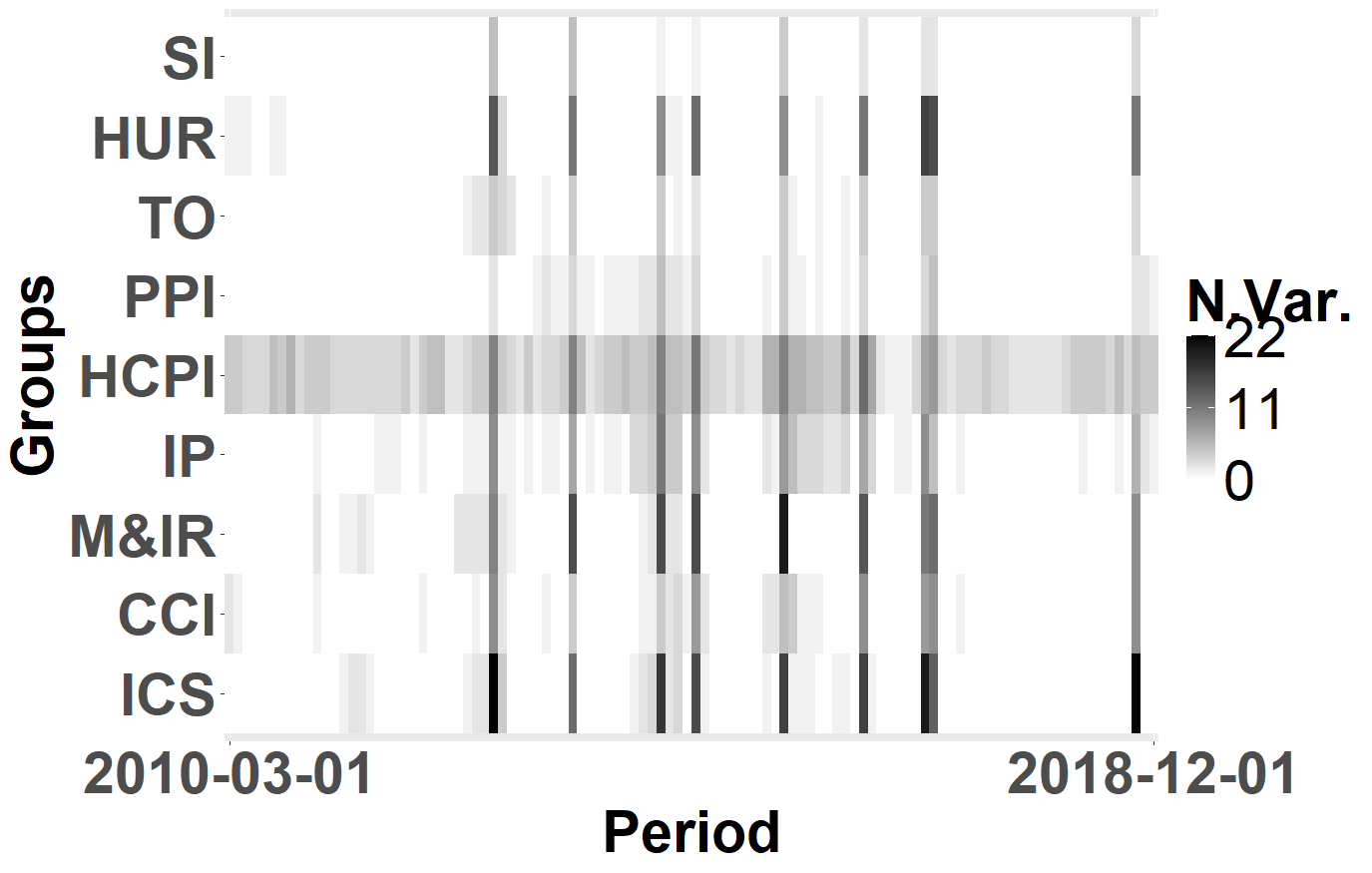}}

\caption{\footnotesize Heatmaps representing the number of 
variables selected by LAS (left)
and ARMAr-LAS (right)
in the nine main domains. 
}\label{Heatmaps}
\end{figure}

Figure~\ref{Heatmaps} summarizes patterns of selected predictors over time for LAS and ARMAr-LAS. The heatmaps represent the number of selected variables categorized according to the nine main domains (see above). LAS selects 
predictors largely, though not exclusively, from the domains ICS, M\&IR and HUR.
ARMAr-LAS is more targeted, selecting
predictors almost exclusively in the HCPI domain (in Supplement~\ref{MostSelPred}, we report the top 5 predictors in terms of selection frequency across forecasting samples). 
Note, however, that in a few instances (3 for $h=12$ and 2 for $h=24$) ARMAr-LAS does select many more predictors across multiple groups. Interestingly, these correspond to the period of the financial crisis (between the end of 2008 and the beginning of 2010), when negative shocks in some of the variables under consideration may indeed create a more complex picture in terms of feature selection. 
In summary, ARMAr-LAS exploits cross-sectional information mainly focusing on prices, and accrues a forecasting advantage -- as LAS uses many more variables to produce significantly worse forecasts. 

\section{Concluding Remarks}\label{ConcRem}
In this paper, we demonstrated that the probability of spurious correlations between stationary orthogonal or weakly correlated processes depends not only on the sample size, but also on the degree of predictors serial correlation. This motivates predictor prewhitening as a way to improve the finite-sample design used by LASSO. In order to improve the performance of LASSO in a time series context, we proposed an approach based on applying LASSO to pre-whitened (i.e., ARMA filtered) time series. In our experiments, this working model mitigates spurious sample correlations and improves estimation and
forecasting. We characterized prediction and estimation error bounds, for our proposal. Furthermore, we assessed its performance through Monte Carlo simulations and an empirical application to Euro Area macroeconomic time series. Through simulations, we observed that ARMAr-LASSO, i.e., LASSO applied to ARMA residuals, reduces the probability of large spurious correlations and outperforms other LASSO-based methods from the literature in terms of both coefficient estimation and forecasting. The empirical application confirms that ARMAr-LASSO improves forecasting performance and produces more parsimonious models.

Future research could extend ARMAr-LASSO to high-dimensional inference through a desparsified estimator incorporating Newey--West long-run covariance estimation, building on~\cite{Chernozhukov2021}. Further empirical applications could assess its performance on U.S. macroeconomic data, such as the FRED-MD dataset.

\newpage

\clearpage
\pagenumbering{arabic}

\renewcommand{\theequation}{S.\arabic{equation}}
\renewcommand{\thesection}{S\arabic{section}}
\renewcommand{\thepage}{S\arabic{page}}
\renewcommand{\thefigure}{S.\arabic{figure}}
\renewcommand{\thetable}{S.\arabic{table}}

\renewcommand{\thetheorem}{\arabic{section}.\arabic{theorem}}
\renewcommand{\theass}{\arabic{section}.\arabic{ass}}
\renewcommand{\theprop}{\arabic{section}.\arabic{prop}}
\renewcommand{\thecor}{\arabic{section}.\arabic{cor}}
\renewcommand{\thelemma}{\arabic{section}.\arabic{lemma}}
\renewcommand{\theexample}{\arabic{section}.\arabic{example}}
\renewcommand{\thedefinition}{\arabic{section}.\arabic{definition}}
\renewcommand{\therem}{\arabic{section}.\arabic{rem}}

\bigskip
\begin{center}
{\LARGE\bf Supplement - ARMAr-LASSO: Mitigating the Impact of Predictor Serial Correlation on the LASSO}
\end{center}

\begin{appendices}

\sloppy

Here we reports supplementary analyses and simulation experiments omitted from
the main manuscript. \sloppy
Code and replicability 
materials are at \zenodoLink

\section{Proofs}\label{sec:Appendix}

In this section, we provide the proof of our theoretical results.

\subsection{Proof of Lemma~\ref{lem:sample_correlations}}
\label{app:proof_lemma1}

Let $a_t=z_{i,t}z_{j,t}$ and $\ddot\phi=\phi_i\phi_j$. Since the two processes are independent and have unit variance, $E(a_t)=0$ and $\operatorname{Cov}(a_t,a_{t-\ell})=\ddot\phi^{|\ell|}$. Therefore, $\operatorname{Var}(A_{ij,T})=T^{-2}\sum_{t=1}^T\sum_{r=1}^T\operatorname{Cov}(a_t,a_r)=T^{-2}[T+2\sum_{\ell=1}^{T-1}(T-\ell)\ddot\phi^\ell]$, proving part (a).

For part (b), write $z_{k,t}=\phi_k z_{k,t-1}+\eta_{k,t}$, where $\eta_{k,t}=\sigma_k^{-1}\sqrt{1-\phi_k^2}\,u_{k,t}$, so that $E(\eta_{k,t})=0$ and $\operatorname{Var}(\eta_{k,t})=1-\phi_k^2$. Then $a_t=\ddot\phi a_{t-1}+d_t$, where $d_t=\phi_i z_{i,t-1}\eta_{j,t}+\phi_j z_{j,t-1}\eta_{i,t}+\eta_{i,t}\eta_{j,t}$. The sequence $\{d_t\}$ is a stationary martingale difference sequence with respect to the filtration $\mathcal F_t=\sigma(u_{i,s},u_{j,s}:s\leq t)$, with $E(d_t^2)=1-\ddot\phi^2$.

Moreover, $E(d_t^2\mid\mathcal F_{t-1})=\phi_i^2(1-\phi_j^2)z_{i,t-1}^2+\phi_j^2(1-\phi_i^2)z_{j,t-1}^2+(1-\phi_i^2)(1-\phi_j^2)$. By stationarity and ergodicity, $T^{-1}\sum_{t=1}^T E(d_t^2\mid\mathcal F_{t-1})\xrightarrow{p}1-\ddot\phi^2$. Since $E(d_t^2)<\infty$, the conditional Lindeberg condition also holds. Hence, the martingale central limit theorem yields $T^{-1/2}\sum_{t=1}^T d_t\xrightarrow{d}\mathcal{N}(0,1-\ddot\phi^2)$.

Since $(1-\ddot\phi)\sum_{t=1}^T a_t=\sum_{t=1}^T d_t-\ddot\phi(a_T-a_0)$ and $a_T-a_0=O_p(1)$, it follows that $\sqrt{T}\,A_{ij,T}\xrightarrow{d}\mathcal{N}(0,(1+\ddot\phi)/(1-\ddot\phi))$. Finally, the sample means of the standardized processes are $O_p(T^{-1/2})$, while their sample variances converge in probability to one. Thus $\sqrt{T}(\widehat{\rho}_{ij}^x-A_{ij,T})\xrightarrow{p}0$, and the result follows from Slutsky's theorem. 
 \hfill$\square$\\

\subsection{Proof of Lemma~\ref{lem:ARMAscore}}\label{lem:ARMAscore_proof}
Since the lagged-response regressors are unchanged, $\boldsymbol{\Delta}$ is nonzero only in its first $n_T$ rows, and
$$
\widehat v_t^{\,*}=v_t-\sum_{i=1}^{n_T}\alpha_i^*(\widehat u_{i,t}-u_{i,t}).
$$
For any $j\in\{1,\ldots,n_T+p_y\}$ and $1\leq\ell\leq T$,
$$
\frac{1}{T}\sum_{t=1}^{\ell}(\widehat w_{j,t}\widehat v_t^{\,*}-w_{j,t}v_t)=\frac{1}{T}\sum_{t=1}^{\ell}(\widehat w_{j,t}-w_{j,t})v_t-\frac{1}{T}\sum_{t=1}^{\ell}\widehat w_{j,t}\sum_{i=1}^{n_T}\alpha_i^*(\widehat u_{i,t}-u_{i,t}).
$$
By Cauchy--Schwarz and the definition of $\delta_T$,
$$
\left|\frac{1}{T}\sum_{t=1}^{\ell}(\widehat w_{j,t}-w_{j,t})v_t\right|\leq\delta_T\|\mathbf v\|_T.
$$
Similarly, the triangle and Cauchy--Schwarz inequalities yield
$$
\left|\frac{1}{T}\sum_{t=1}^{\ell}\widehat w_{j,t}\sum_{i=1}^{n_T}\alpha_i^*(\widehat u_{i,t}-u_{i,t})\right|\leq M_T\sum_{i=1}^{n_T}|\alpha_i^*|\left(\frac{1}{T}\sum_{t=1}^T(\widehat u_{i,t}-u_{i,t})^2\right)^{1/2}\leq M_T\|\pmb\alpha^*\|_1\delta_T.
$$
Combining the bounds and taking the maximum over $j$ and $\ell$ proves the first claim. Assumption~\ref{ass:ARMAuniform} and the stated stochastic boundedness conditions imply
$$
\max_{1\leq\ell\leq T}\left\|\frac{1}{T}\sum_{t=1}^{\ell}(\widehat{\mathbf w}_t\widehat v_t^{\,*}-\mathbf w_tv_t)\right\|_\infty=O_p\left(r_{T,n_T}(1+\|\pmb\alpha^*\|_1)\right)=o_p(\lambda_T).
$$
Finally, by the triangle inequality,
$$
\max_{1\leq\ell\leq T}\left\|\frac{1}{T}\sum_{t=1}^{\ell}\widehat{\mathbf w}_t\widehat v_t^{\,*}\right\|_\infty\leq\max_{1\leq\ell\leq T}\left\|\frac{1}{T}\sum_{t=1}^{\ell}\mathbf w_tv_t\right\|_\infty+o_p(\lambda_T).
$$
The asserted high-probability implication follows because $b>a$.
 \hfill$\square$\\

\subsection{Proof of Theorem~\ref{theo:PrA}}\label{theo:PrAB_LEMMAS_Proof}
We start by introducing some important definitions.

\begin{definition}\label{def1}
    Let $\left( \Omega, \mathcal{F},P\right)$ be a probability space and let $\mathcal{G}$ and $\mathcal{H}$ be sub-$\mathbf{\sigma}$-fields of $\mathcal{F}$. Then
\[\alpha(\mathcal{G},\mathcal{H})=\sup_{G\in\mathcal{G},H\in\mathcal{H}}|\Pr(G\cap H)-\Pr(G)\Pr(H)|\]
is known as the strong mixing coefficient. For a sequence $\left\{X_t\right\}_{-\infty}^{+\infty}$ let $\left\{\mathcal{F}_{-\infty}^{t}\right\}=\mathbf{\sigma}(\dots,X_{t-1},X_t)$ and similarly define $\left\{\mathcal{F}_{t+m}^{\infty}\right\}=\mathbf{\sigma}(X_{t+m},X_{t+m+1},\dots)$. The sequence is said to be $\alpha$-mixing (or strong mixing) if lim$_{m\rightarrow\infty}\alpha_m=0$ where
\[\alpha_m=\sup_t\alpha(\mathcal{F}_{-\infty}^{t},\mathcal{F}_{t+m}^{\infty}).\]
\end{definition}

\begin{definition}\label{def2}
    (Mixingale,~\cite{davidson1994}, ch. 16). The sequence of pairs $\left\{X_t,\mathcal{F}\right\}_{-\infty}^{+\infty}$ in a filtered probability space $\left( \Omega, \mathcal{F},P\right)$ where the $X_t$ are integrable r.v.s is called $L_p$-mixingale if, for $p\geq1$, there exist sequences of non-negative constants $\left\{d_t\right\}_{-\infty}^{+\infty}$ and $\left\{\nu_m\right\}_{0}^{\infty}$ such that $\nu_m\rightarrow0$ as $m\rightarrow\infty$ and 
    \[||E(X_t|\mathcal{F}_{t-m})||_p\leq d_t\nu_m\]
    \[||X_t-E(X_t|\mathcal{F}_{t+m})||_p\leq d_t\nu_{m+1},\]
    hold for all $t$ and $m\geq0$. Furthermore, we say that $\left\{X_t\right\}$ is $L_p$-mixingale of size -$a$ with respect to  $\mathcal{F}_t$ if $\nu_m=O(m^{-a-\epsilon})$ for some $\epsilon>0$.
\end{definition}

\begin{definition}\label{def3}
    (Near-Epoch Dependence,~\cite{davidson1994}, ch. 17). For a possibly vector-valued stochastic sequence $\left\{\mathbf{V}_t\right\}_{-\infty}^{+\infty}$, in a probability space $\left( \Omega, \mathcal{F},P\right)$ let $\mathcal{F}_{t-m}^{t+m}=\mathbf{\sigma\left(\mathbf{V}_{t-m},\dots,\mathbf{V}_{t+m}\right)}$, such that $\left\{\mathcal{F}_{t-m}^{t+m}\right\}_{m=0}^{\infty}$ is a non-decreasing sequence of $\mathbf{\sigma}$-fields. If for $p>0$ a sequence of integrable r.v.s $\left\{X_t\right\}_{-\infty}^{+\infty}$ satisfies 
    \[||X_t-E(X_t|{F}_{t-m}^{t+m})||_p\leq d_t\nu_m,\]
    where $\nu_m\rightarrow0$ and $\left\{d_t\right\}_{-\infty}^{+\infty}$ is a sequence of positive constants, $X_t$ will be said to be near-epoch dependent in $L_p$-norm ($L_p$-NED) on $\left\{\mathbf{V}_t\right\}_{-\infty}^{+\infty}$. Furthermore, we say that $\left\{X_t\right\}$ is $L_p$-NED of size -$a$ on $\mathbf{V}_t$ if $\nu_m=O(m^{-a-\epsilon})$ for some $\epsilon>0$.
\end{definition}

\noindent Note that we use the same notation for the constants $d_t$ and sequence $\nu_m$ as for the near-epoch dependence, since they play the same role in both types of dependence.

\begin{rem}
The error term $v_t$ is modeled as a stationary Near-Epoch-Dependent (NED) process which, under appropriate decay conditions, can be approximated by a strongly mixing sequence. In particular, this includes the case of stationary ARMA processes under suitable innovation regularity, which can be geometrically mixing of the mixing coefficients (see, e.g., \citealt{davidson1994}). For valid asymptotic results, we assume that the error term has finite $q$-th moments for some $q>4$, ensuring tails sufficiently light for our limit theory, although stronger assumptions such as exponential moment bounds could also be accommodated.
\end{rem}

\begin{rem}[NED and mixingale dependence] Assumption~\ref{ass:ScoreMix} can be justified through the NED framework of \cite{davidson1994} when the oracle working- model variables are uniformly NED on an underlying mixing process, with sufficiently fast dependence decay and appropriate moment conditions. In particular, the NED-to-mixingale results in Chapter 17 provide a route to establishing the uniform dependence condition required for the maximal inequality used in Theorem~\ref{theo:PrA}.
\end{rem}

\bigskip

\noindent Throughout the proofs, $C$ and $C_i$ denote generic positive finite constants, independent of $T$ and $n_T$, whose values may change from line to line.

\bigskip

\noindent\textbf{Proof of Theorem~\ref{theo:PrA}}
Set $\xi_{j,t}=w_{j,t}v_t$ for $j\leq n_T+p_y$. Since $\pmb{\beta}^*$ is the population linear-projection coefficient, $E(\xi_{j,t})=0$. By Assumption~\ref{ass:ScoreMix}, these oracle-score processes are uniformly $L_{c_1}$-mixingale with summable dependence coefficients. We first establish the maximal-moment bound used in the original NED/mixingale argument.

Let $D_{j,t,m}=E(\xi_{j,t}\mid\mathcal F_{t-m})-E(\xi_{j,t}\mid\mathcal F_{t-m-1})$, $m\geq0$. For each fixed $m$, $(D_{j,t,m})_t$ is a martingale-difference sequence with respect to $(\mathcal F_{t-m})_t$, and Assumption~\ref{ass:ScoreMix} gives $\|D_{j,t,m}\|_{c_1}\leq a_m+a_{m+1}$. Since $\xi_{j,t}$ is adapted and $\|E(\xi_{j,t}\mid\mathcal F_{t-m})\|_{c_1}\to0$, we have $\xi_{j,t}=\sum_{m=0}^{\infty}D_{j,t,m}$ in $L_{c_1}$. For $c_1>2$, the martingale maximal inequality and Minkowski's inequality therefore yield
$
\left\|\max_{1\leq\ell\leq T}\left|\sum_{t=1}^{\ell}\xi_{j,t}\right|\right\|_{c_1}\leq C\sqrt T\sum_{m=0}^{\infty}(a_m+a_{m+1})\leq C'\sqrt T,
$
uniformly over $j\leq n_T+p_y$. Equivalently,
\begin{equation}\label{eq:ARMAr_oracle_maximal}
\sup_{j\leq n_T+p_y}E\left[\max_{1\leq\ell\leq T}\left|\sum_{t=1}^{\ell}w_{j,t}v_t\right|^{c_1}\right]\leq K T^{c_1/2}.
\end{equation}
This is the maximal inequality underlying the NED/mixingale approach of \cite{davidson1994}; the preceding remark explains how the NED framework motivates Assumption~\ref{ass:ScoreMix}.

Define the oracle-score maximum
$
S_T^o=\max_{1\leq j\leq n_T+p_y,\,1\leq\ell\leq T}\left|\frac1T\sum_{t=1}^{\ell}w_{j,t}v_t\right|.
$
By the union bound, Markov's inequality and~\eqref{eq:ARMAr_oracle_maximal},
$$
P\left(S_T^o>\frac{\lambda_T}{8}\right)\leq\frac{8^{c_1}}{(T\lambda_T)^{c_1}}\sum_{j=1}^{n_T+p_y}E\left[\max_{\ell\leq T}\left|\sum_{t=1}^{\ell}w_{j,t}v_t\right|^{c_1}\right]\leq\frac{K8^{c_1}(n_T+p_y)}{T^{c_1/2}\lambda_T^{c_1}}.
$$
The lower bound on $\lambda_T$ implies
$$
P\left(S_T^o>\frac{\lambda_T}{8}\right)\leq\frac{C(n_T+p_y)}{(2n_T+p_y)\sqrt{\log T}}\leq C(\sqrt{\log T})^{-1}.
$$

Next, define the perturbation due to the estimated ARMA innovations,
$$
R_T=\max_{1\leq\ell\leq T}\left\|\frac1T\sum_{t=1}^{\ell}\left(\widehat{\mathbf w}_t\widehat v_t^{\,*}-\mathbf w_tv_t\right)\right\|_\infty.
$$
By Lemma~\ref{lem:ARMAscore}, Assumption~\ref{ass:ARMAuniform} and the stated rate conditions,
$$
R_T\leq\delta_T\left(\|\mathbf v\|_T+M_T\|\pmb{\alpha}^*\|_1\right)=o_p(\lambda_T),
$$
so that $P(R_T>\lambda_T/8)=o(1)$. Finally, the triangle inequality gives
$$
\max_{1\leq\ell\leq T}\left\|\frac1T\sum_{t=1}^{\ell}\widehat{\mathbf w}_t\widehat v_t^{\,*}\right\|_\infty\leq S_T^o+R_T.
$$
Thus, on the event $\{S_T^o\leq\lambda_T/8\}\cap\{R_T\leq\lambda_T/8\}$, we have $\mathcal{A}_T$. Consequently,
$$
P(\mathcal{A}_T)\geq1-P\left(S_T^o>\frac{\lambda_T}{8}\right)-P\left(R_T>\frac{\lambda_T}{8}\right)\geq1-C(\sqrt{\log T})^{-1}-o(1),
$$
which completes the proof.
 \hfill$\square$\\

\begin{rem}\label{rem:NED}
    The concentration arguments in Theorem~\ref{theo:PrA} are based on the NED-mixingale framework of~\cite{ADAMEK2023}, which allows for very general forms of temporal dependence, random regressors, and only finite moment assumptions. This choice is particularly suited to our projection-based and potentially misspecified regression setting with lagged dependent variables and estimated components. Sharper deviation bounds, and hence more aggressive tuning rates, could be obtained under stronger structural assumptions using, for example, the functional dependence approach of~\cite{WuWu2016} or Fuk-Nagaev type inequalities for $\tau$- or $s$-mixing processes as in~\cite{Babii2024}. These alternatives, however, require additional restrictions on mixing rates, tail behavior, and the design structure. We therefore favor the more general framework of~\cite{ADAMEK2023}, which yields slightly more conservative but broadly valid theoretical guarantees.

\end{rem}

\subsection{Proof of Theorem~\ref{theo:Bounds}}
The argument follows the usual LASSO basic-inequality approach (see~\citealt{Buhlmann2011,ADAMEK2023}), adapting their argument to the feasible design under exact
sparsity. Set $\mathbf d=\widehat{\pmb{\beta}}-\pmb{\beta}^*$, $Q_T=T^{-1}\|\widehat{\mathbf W}'\mathbf d\|_2^2$, $a=\|\mathbf d_{\widetilde S}\|_1$, and $b=\|\mathbf d_{\widetilde S^c}\|_1$. Since $\mathbf y=\widehat{\mathbf W}'\pmb{\beta}^*+\widehat{\mathbf v}^{\,*}$, optimality of $\widehat{\pmb{\beta}}$ yields
\begin{equation}\label{eq:ARMAr_T4_basic_ineq}
\frac12Q_T\leq\frac1T(\widehat{\mathbf W}\widehat{\mathbf v}^{\,*})'\mathbf d+\lambda_T\big(\|\pmb{\beta}^*\|_1-\|\widehat{\pmb{\beta}}\|_1\big).
\end{equation}
On $\mathcal{A}_T$, Theorem~\ref{theo:PrA} implies $T^{-1}\|\widehat{\mathbf W}\widehat{\mathbf v}^{\,*}\|_\infty\leq\lambda_T/4$. By H\"older's inequality and $\|\pmb{\beta}^*\|_1-\|\widehat{\pmb{\beta}}\|_1\leq a-b$, inequality~\eqref{eq:ARMAr_T4_basic_ineq} gives
\begin{equation}\label{eq:ARMAr_T4_cone}
Q_T+\frac{3\lambda_T}{2}b\leq\frac{5\lambda_T}{2}a.
\end{equation}
Thus $b\leq5a/3\leq3a$. The compatibility condition in Assumption~\ref{ass:REC}, applied on this cone, gives
\begin{equation}\label{eq:ARMAr_T4_compat}
Q_T\geq\frac{\gamma_{\widehat w}^2}{\widetilde s}a^2.
\end{equation}
Combining~\eqref{eq:ARMAr_T4_cone} and~\eqref{eq:ARMAr_T4_compat}, and maximizing the resulting quadratic expression over $a\geq0$, we obtain
$$
\|\mathbf d\|_1=a+b\leq\frac83a-\frac{2Q_T}{3\lambda_T}\leq\frac83a-\frac{2\gamma_{\widehat w}^2}{3\widetilde s\lambda_T}a^2\leq\frac{8\widetilde s\lambda_T}{3\gamma_{\widehat w}^2}=C\frac{\widetilde s\lambda_T}{\gamma_{\widehat w}^2},
$$
proving point ({\textit a}).

For the prediction bound, the KKT conditions yield a subgradient $\widehat{\mathbf z}$ with $\|\widehat{\mathbf z}\|_\infty\leq1$ such that
$$
\frac1T\widehat{\mathbf W}\big(\mathbf y-\widehat{\mathbf W}'\widehat{\pmb{\beta}}\big)=\lambda_T\widehat{\mathbf z}.
$$
Using the feasible residual identity, this becomes $T^{-1}\widehat{\mathbf W}\widehat{\mathbf W}'\mathbf d=T^{-1}\widehat{\mathbf W}\widehat{\mathbf v}^{\,*}-\lambda_T\widehat{\mathbf z}$. Taking its inner product with $\mathbf d$ and applying the previous bound yields
$$
Q_T\leq\left(\frac{\lambda_T}{4}+\lambda_T\right)\|\mathbf d\|_1\leq\frac{10\widetilde s\lambda_T^2}{3\gamma_{\widehat w}^2}=C\frac{\widetilde s\lambda_T^2}{\gamma_{\widehat w}^2},
$$
which proves ({\textit b}). Finally, $|\widehat\beta_j-\beta_j^*|\leq\|\mathbf d\|_1$ for every $j$, so the stated beta-min condition excludes false negatives on the same event. Theorem~\ref{theo:PrA} gives $P(\mathcal{A}_T)\geq1-C(\sqrt{\log T})^{-1}-o(1)$, completing the proof.
 \hfill$\square$\\

\subsection{Proof of Corollary~\ref{rem4}}
By Theorem~\ref{theo:Bounds}, with probability tending to one,
$$
\left\|\widehat{\pmb{\beta}}-\pmb{\beta}^*\right\|_1\leq\frac{C\widetilde s\lambda_T}{\gamma_{\widehat w}^2},\qquad
\frac{1}{T}\left\|\widehat{\mathbf W}'(\widehat{\pmb{\beta}}-\pmb{\beta}^*)\right\|_2^2\leq\frac{C\widetilde s\lambda_T^2}{\gamma_{\widehat w}^2}.
$$
Substituting the assumed order of $\lambda_T$ and using $\gamma_{\widehat w}^{-2}=O_p(1)$ yield the stated rates, since $2n_T+p_y\asymp n_T+p_y$. Finally, as $\widetilde s\geq1$, the condition $\widetilde s\lambda_T\to0$ implies $\lambda_T\to0$ and hence $\widetilde s\lambda_T^2\to0$. Both error bounds therefore vanish in probability.
 \hfill$\square$\\

\section{Upper Bound for Minimum Eigenvalue}\label{UppMinPsy}
\sloppy
Here, we would like to point out the  role of $\widehat{\rho}_{ij}^x$ for $\widehat{\psi}_{min}^x$. To this end, we start by recalling an inequality that links off-diagonal elements and eigenvalues of $\widehat{\mathbf{R}}_x$; namely, $\widehat{\psi}_{min}^x\leq1-\underset{i\neq j}{\text{max}}|\widehat{\rho}_{ij}^x|$.
Because of this, for any given $\tau\in[0,1)$ we have
\[
\Pr\left(\widehat{\psi}_{min}^x\leq1-\tau\right)\geq \Pr\left(1-\max_{i\neq j}|\widehat{\rho}_{ij}^x|\leq1-\tau\right) \geq \Pr\left(1-|\widehat{\rho}_{i\neq j}^x|\leq1-\tau\right)=\Pr\left(|\widehat{\rho}_{i\neq j}^x|\geq\tau\right)
\]
which emphasizes how the probability of a generic sample correlation being large in absolute value affects the probability of the minimum eigenvalue being small -- and thus the estimation error bounds of the LASSO, as established by~\cite{Bickel2009}. As the next example shows, point the inequality $\widehat{\psi}_{min}^x\leq1-\underset{i\neq j}{\text{max}}|\widehat{\rho}_{ij}^x|$ can be easily fixed.
\begin{example}
Let $e_i$ and $e_j$ be vectors from the standard basis of $\mathbb{R}^n$, $i,j\in{1,\dots,N}$. Moreover, let $x_{\pm}=2^{-1/2}(e_i\pm e_j)$, satisfying $||x_{\pm}||_2=1$, and let $A$ be a correlation matrix with $a_k$ be the $k$-th column. Then we have
\[x_{\pm}'Ax_{\pm}=\frac{1}{2}(e_i\pm e_j)'(a_i\pm a_j)=\frac{1}{2}(a_{ii}\pm2a_{ij}+a_{jj})=1\pm a_{ij}\ \ .\]
Thus, $\psi_{min}\leq1-|a_{ij}|$ for all $i\neq j$ and so
\[\psi_{min}\leq1-\max_{i\neq j}|a_{ij}|\ \ .\]
\end{example}

\section{Distribution of Sample Correlation Between Serially Correlated Processes: Simulation Experiments}\label{MonteCarlo}
In this Section, we expand the results in Section~\ref{sec:spurious_correlations} in more generic contexts, relaxing the assumption that the covariates are orthogonal Gaussian AR(1) processes. 

\bigskip

\noindent
{\bf The Impact of $Sign(\ddot{\phi})$}
\normalsize

\noindent 
Consider a first order $2$-variate autoregressive process $\mathbf{x}_t=\pmb{\phi}\mathbf{x}_{t-1}+\mathbf{u}_t$, $t=1,\dots,T$, where $\pmb{\phi}$ is the $2\times 2$ diagonal matrix with $\text{diag}(\pmb{\phi})=(\phi_1,\phi_2)$. Lemma 1 shows that the sign of $\ddot\phi$ determines the asymptotic dispersion of Pearson correlations. In particular, when $-1<\ddot{\phi}<0$, an increment on $|\ddot{\phi}|$ makes the density of $\widehat{\rho}_{12}^x$ more concentrated around $0$. 
In order to 
validate this result, 
we 
run simulations with
$T=100$ and 
different values for 
the second element of the diagonal of $\pmb{\phi}$; namely,
$-0.3,-0.6,-0.9,-0.95$. 
In this case, we see that when $Sign(\phi_1)\neq Sign(\phi_2)$ and $|\ddot{\phi}|$ increases,
$d(r)$ 
increases its concentration around $0$.

\bigskip

\noindent{\bf General Case}

\begin{figure}[t]
\graphicspath{{images/}}
\centering
{\includegraphics[scale=0.35]{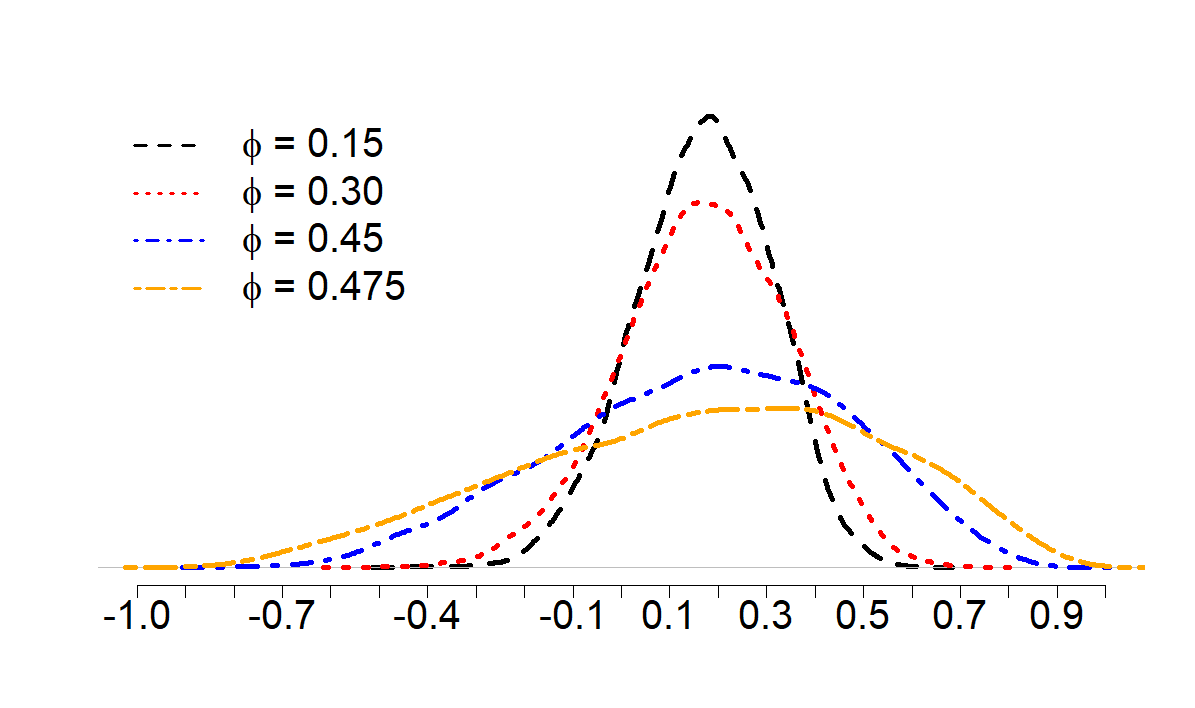}}
\caption{\footnotesize 
Densities for $\widehat{\rho}_{12}^x$
in the case of Laplace weakly correlated AR(3) and ARMA(2,1) processes, for $T=100$ and various 
$\phi$'s.}
\label{fig:GenCas}
\end{figure}

\noindent To generalize 
our findings to the case of non-Gaussian weakly correlated AR and ARMA processes,
we generate covariates according to the following DGPs:
${x}_{1t}=(\phi+0.1){x}_{1t-1}+(\phi+0.1){x}_{1t-2}-0.2{x}_{1t-3}+u_{1t}$, and ${x}_{2t}=\phi x_{2t-1}+\phi x_{2t-2}+{u}_{2t}+0.8{u}_{2t-1}$, where $t=1,
\ldots,100$ and $\phi=0.15, 0.3, 0.45, 0.475$. Moreover, we generate $u_{1t}$ and $u_{2t}$ 
from a bivariate Laplace distribution with means $0$, 
variances $1$, and $c_{12}^u=0.2$.
In these more general cases, we do not know an approximate 
theoretical density 
for $\widehat{\rho}_{12}^u$.
Therefore, we rely entirely on 
simulations to show the effect of serial correlation on $\Pr\left\{|\widehat{\rho}_{12}^x|\geq\tau\right\}$. Figure~\ref{fig:GenCas} shows $d(r)$ obtained 
from 5000 Monte Carlo 
replications 
for the different values of $\phi$. In short, also in 
the more general 
cases where covariates are non-Gaussian, weakly correlated AR(3) and ARMA(2,1) processes, the probability of getting large sample cross-correlations 
depends on the degree of serial correlation. 
More simulation results are provided below.
\begin{figure}[t]
\graphicspath{{images/}}
\centering
  \captionsetup[subfigure]{oneside,margin={0.5cm,0cm}}
{\includegraphics[width=5.4cm]{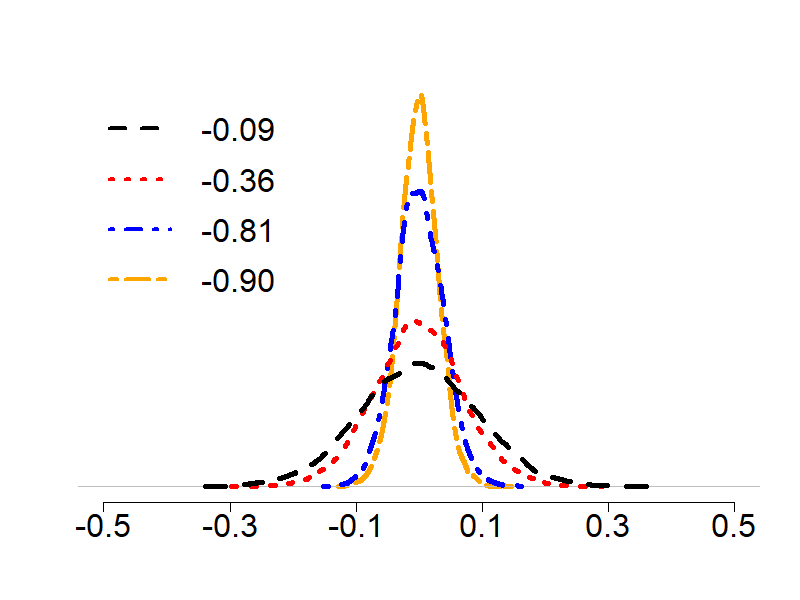}}\hfil
\caption{\footnotesize 
Monte Carlo 
densities for
$\widehat{\rho}_{12}^x$, 
for $T$=100 and various (negative) 
$\ddot{\phi}$'s.}
\label{negSignPhi}
\end{figure}

\subsection{More General Cases}\label{MoreCases}

We study the density of $\widehat{\rho}_{12}^x$ in three different cases: non-Gaussian processes; weakly and high cross-correlated processes; and ARMA processes with different order. Note that for the first two cases the variables are AR(1) processes with $T=100$ and autocorrelation coefficient $\phi=0.3,0.6,0.9,0.95$. 
We rely on the densities obtained on 5000 Monte Carlo replications, i.e. $d(r)$, to show the effect of serial correlation on $\Pr\left\{|\widehat{\rho}_{12}^x|\geq\tau\right\}$.

\vspace{1cm}

\noindent{\bf The Impact of non-Gaussianity}\normalsize

\begin{figure}[t]
\graphicspath{{images/}}
\centering
  \captionsetup[subfigure]{oneside,margin={0.5cm,0cm}}
\subfloat[\footnotesize Laplace]{\includegraphics[width=5cm]{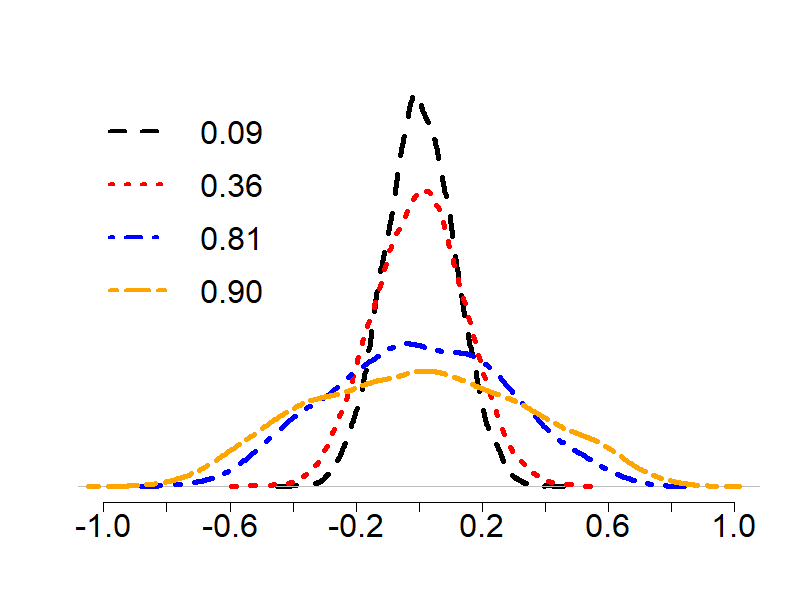}}\hfil
\subfloat[\footnotesize Cauchy]{\includegraphics[width=5cm]{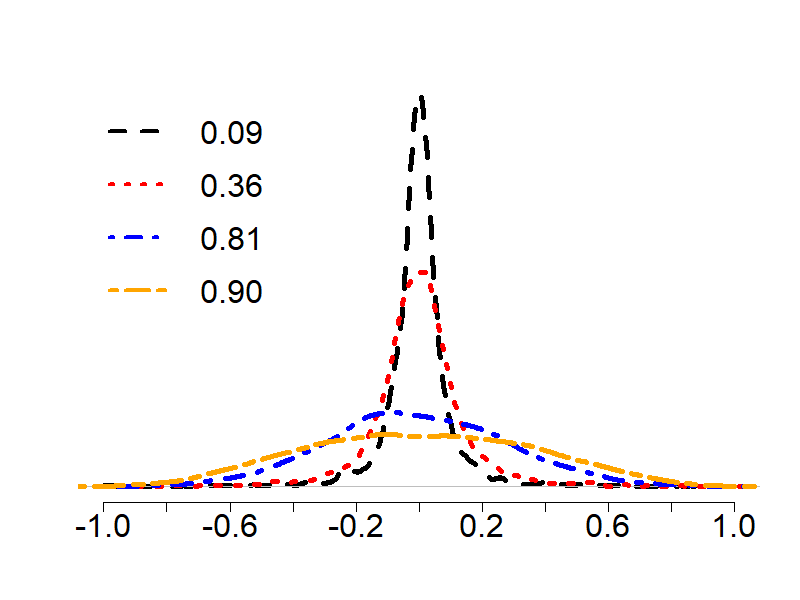}}\hfil 
\subfloat[\footnotesize t-Student]{\includegraphics[width=5cm]{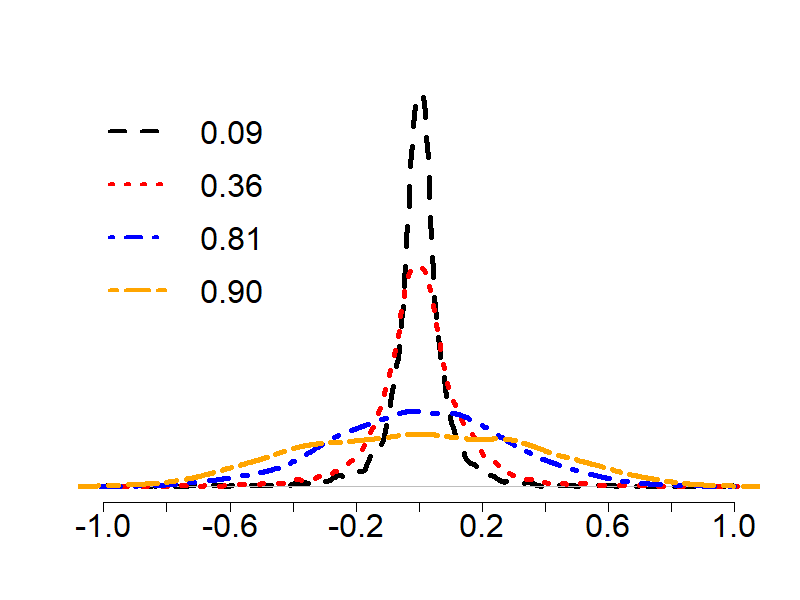}} 
\caption{\footnotesize Simulated density of $\widehat{\rho}_{12}^x$ in the case of non-Gaussian processes, for $T=100$ and various values of $\ddot\phi$.}\label{NoGauss}
\end{figure}

\noindent The theoretical contribution reported in Section~\ref{sec:spurious_correlations} of the main text requires the Gaussianity of $u_1$ and $u_2$. With the following simulation experiments we show that the impact of $\ddot{\phi}$ on the density of $\widehat{\rho}_{12}^x$ is relevant also when $u_{1t}$ and $u_{2t}$ are non-Gaussian  random variables. To this end, we generate $u_{1t}$ and $u_{2t}$ from the following distributions: Laplace with mean 0 and variance 1 (case (a)); Cauchy with location parameter 0 and scale parameter 1 (case (b)); and from a $t$-student with 1 degree of freedom (case (c)). Figure~\ref{NoGauss} reports the results of the simulation experiment. We can state that regardless the distribution of the processes, whenever $Sign(\phi_1)=Sign(\phi_2)$, the probability of large values of $\widehat{\rho}_{12}^x$ increases with $\ddot{\phi}$. As a curiosity, this result is more evident for the case of Laplace variables, whereas for Cauchy and $t$-student the effect of $\ddot{\phi}$ declines. 

\vspace{1cm}

\noindent{\bf The Impact of Population Cross-Correlation}\normalsize

\begin{figure}[t]
\graphicspath{{images/}}
\centering
  \captionsetup[subfigure]{oneside,margin={0.5cm,0cm}}
\subfloat[\footnotesize $c_{12}^u=0.2$]{\includegraphics[width=5cm]{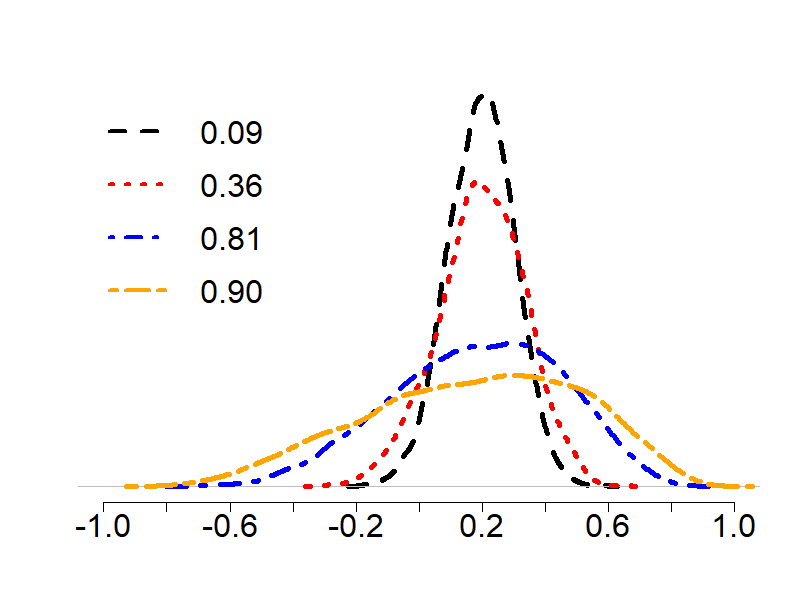}}\hfil
\subfloat[\footnotesize $c_{12}^u=0.8$]{\includegraphics[width=5cm]{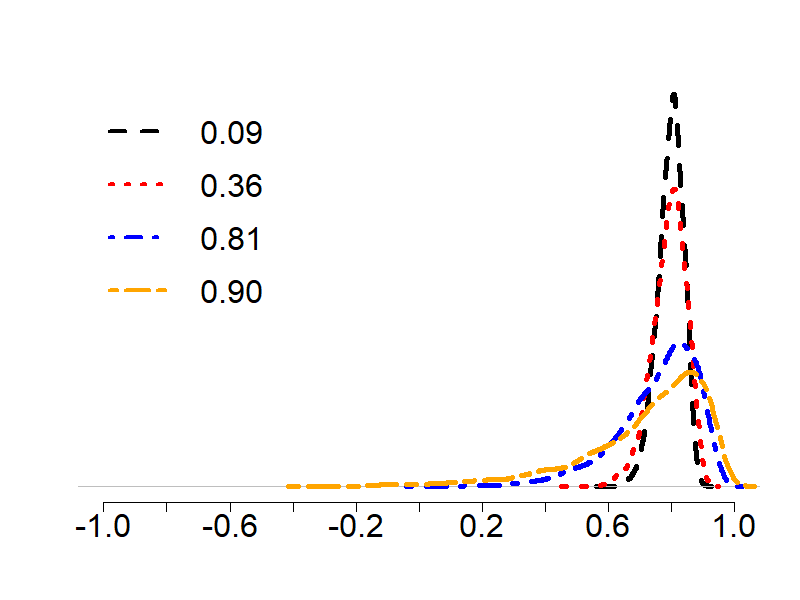}}
\caption{\footnotesize $d(r)$ obtained through simulations in the case of $c_{12}^x=0.2$ (a) and $c_{12}^x=0.8$ (b), for $T=100$ and various values of $\ddot\phi$.}
\label{fig:rho}
\end{figure}

\noindent Since orthogonality is an unrealistic assumption for most economic applications, here we admit population cross-correlation. In Figure~\ref{fig:rho} we report $d(r)$ when the processes are weakly cross-correlated with $c_{12}^u=0.2$, and when the processes are multicollinear with $c_{12}^u=0.8$ (usually we refer to multicollinearity when $c_{12}^u\geq0.7$). We observe that the impact of $\ddot{\phi}$ on $d(r)$ depends on the degree of (population) cross-correlation as follows. In the case of weakly correlated processes, an increase in $\ddot{\phi}$ yields a high probability of observing large sample correlations in absolute value. 
In the case of multicollinear processes, on the other hand, an increase in $\ddot{\phi}$ leads to a high probability of underestimating the true population cross-correlation.

\vspace{1cm}

\noindent{\bf Density of $\widehat{\rho}_{12}^x$ in the case of ARMA$(p_i,q_i)$ processes}\normalsize

\begin{figure}[t]
\begin{center}
\graphicspath{{images/}}
\includegraphics[scale=0.3]{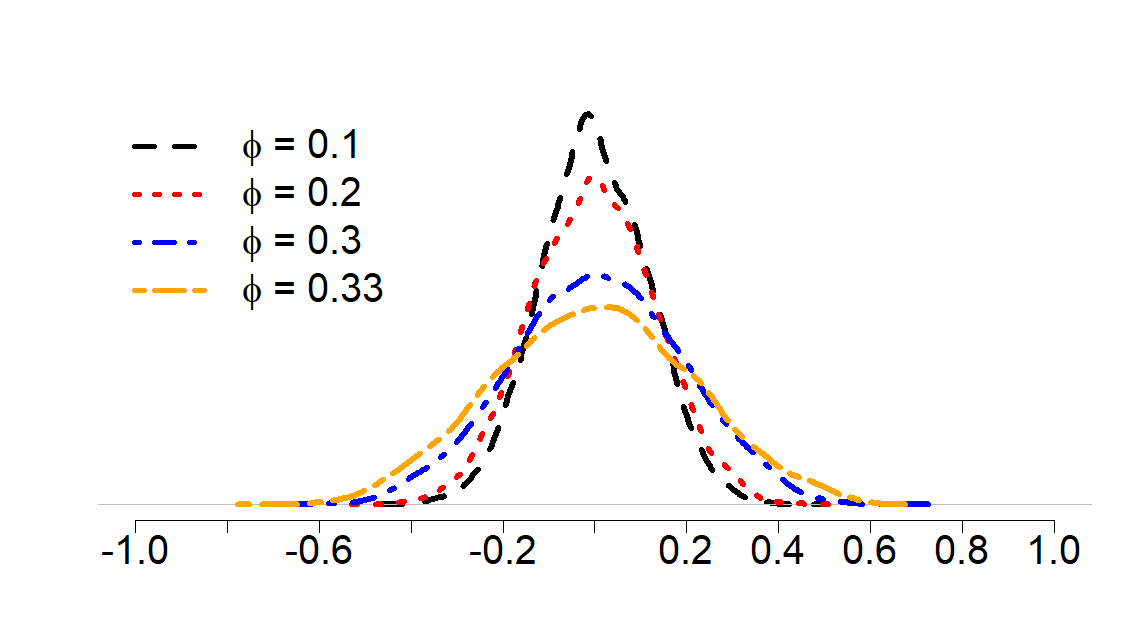}
\end{center}
\caption{\footnotesize Densities of $d(r)$ between two uncorrelated ARMA Gaussian processes, for $T=100$ and various values of $\ddot\phi$.}\label{fig:ARMA}
\end{figure}

\noindent To show the effect of serial correlation on a more general case, we generate $x_1$ and $x_2$ through the following ARMA processes
\begin{align*}
 & {x}_{1t}={\phi}{x}_{1t-1}+\phi{x}_{1t-2}-\phi{x}_{1t-3}{u}_{1t}+0.5u_{1t-1},\\
 & {x}_{2t}={\phi}{x}_{2t-1}+\phi{x}_{2t-2}+{u}_{2t}+0.7{u}_{2t-1}-0.4u_{3t-2},
\end{align*}
where $t=1,\dots,100$ and $u_i\sim \mathcal{N}(0,1)$. In Figure~\ref{fig:ARMA} we report the density of $\widehat{\rho}_{12}^x$ in the case of $T=100$ and $\phi=0.1,0.2,0.3,0.33$. With no loss of generality we can observe that $d(r)$ gets larger as $\phi$ increases, that is $\Pr\left\{|\widehat{\rho}_{12}^x|\geq\tau\right\}$ increases with $|\phi|$.

\section{Comparison with ARDL and GLS Estimators}\label{COMP_ARDLandGLS}
Two natural points of comparison for our proposal are the
AutoRegressive Distributed Lag (ARDL) and the Generalized Least Square (GLS) estimators, which are widely used in the literature to tackle serial correlation.

The ARDL 
consists of regressing the
response on its past realizations -- the autoregressive component -- as well as on current and past values of the 
predictors -- the distributed lag component (see, e.g.,~ \citealp{Panopoulou2004}). Although this method does mitigate 
serial correlation, 
it has the drawback of requiring a very large number of 
coefficients to be estimated. This issue becomes particularly relevant when the sample size is small. 
In contrast, our proposal only requires the addition of a few 
response lags.

The popular Cochrane-Orcutt GLS estimator
approximates the serial correlation structure of the error term while retaining consistent coefficient estimation (see, e.g.,~ \citealp{CocOrc1949}). 
Under suitable exogeneity conditions, GLS can improve efficiency but does not specifically target predictor
persistence. In particular, while the GLS filter may reduce predictors' serial correlations,
it does not remove 
them completely if the AR 
structure of the error term differs from the AR or ARMA 
structures of the 
predictors. 
The GLS-LASSO (\citealp{chronopoulos2023}) can be summarized in the following steps:
\vspace{0.3cm}

\noindent{\it Step 1: Estimation of $\varepsilon_t$.} The estimates of the error term are obtained as $\widehat{\varepsilon}_t=y_t-\mathbf{x}'_t\widetilde{\pmb{\alpha}}$, where $\widetilde{\pmb{\alpha}}$ is the solution to the classical Lasso problem using $\mathbf{X}$ as a design matrix.
\vspace{0.3cm}

\noindent{\it Step 2: Estimation of $\phi_{\varepsilon}$.} The estimates of the parameters of model~\eqref{DGP_varep} of the main text, i.e $\phi_{\varepsilon 1},\dots,\phi_{\varepsilon p_{\varepsilon}}$, is obtained as a solution of the following AR($p_{\varepsilon}$) model
$\widehat{\varepsilon}_{t}=\phi_{1}\widehat{\varepsilon}_{t-1}+\dots+\phi_{\varepsilon p_{\varepsilon}}\widehat{\varepsilon}_{t-p_{\varepsilon}},$
where $\widehat{\varepsilon}_t,\dots,\widehat{\varepsilon}_{t-p_{\varepsilon}}$ are obtained at step 1.
\vspace{0.3cm}

\noindent{\it Step 3: GLS-LASSO.} The LASSO based on the Cochrane-Orcutt GLS filter is

\begin{equation}\label{WM_GLS}
\widehat{\pmb{\alpha}}= \aggregate{argmin}{\pmb{\alpha}\in\mathbb{R}^{n}}\:\left \{\:\frac{1}{2(T)}\:\left|\left|\widetilde{\mathbf{y}}-\widetilde{\mathbf{X}}\pmb{\alpha} \right|\right|_2^2 + \lambda||\pmb{\alpha}||_1\:\right \},
\end{equation}

\noindent where, in scalar representation,
\[\widetilde{y}_t=y_t-\sum_{j=1}^{p_{\varepsilon}}\widehat{\phi}_{\varepsilon j}y_{t-j},\ \ \ \ \ \widetilde{x}_{it}=x_{it}-\sum_{j=1}^{p_{\varepsilon}}\widehat{\phi}_{\varepsilon j}x_{it-j},\ \ \ \ \ t=1,\dots,T,\ \ \ \ \ i=1,\dots,n.\]

The loss function in~\eqref{WM_GLS} corresponds to the $\ell_1$-penalized regression considering the estimates of $\phi_{\varepsilon l}$, $l=1,\dots,p_{\varepsilon}$.~\cite{chronopoulos2023} provide the theoretical properties of this procedure and support them through simulation results. Thus, the working model of GLS-LASSO is

\begin{equation}\label{WorkModel_GLS}
y_t-\sum_{j=1}^{p_{\varepsilon}}\phi_{\varepsilon j}y_{t-j}=\sum_{i=1}^{n}\alpha_{i}^*\left(x_{it}-\sum_{j=1}^{p_{\varepsilon}}\phi_{\varepsilon j}x_{it-j}\right)+\omega_t.
\end{equation}

\noindent Here we compare ARMAr-LASSO and GLM-LASSO in two different cases, namely when the common factor restriction holds and when it does not hold. 

\bigskip

\noindent{\it Common Factor Restriction.}
The common factor restriction holds when predictors and error term are generated by the same AR($p$) process (\citealp{Mizon1995}), as in the Example~\ref{ex:commonAR1} and Remark~\ref{rem1} of the main text. In this case, we can easily observe that the working model of ARMAr-LASSO (see~\eqref{WM} in the main text) and~\eqref{WorkModel_GLS} estimate the true coefficients $\alpha_i^*$ by means of the AR($p$) residuals $u_{it}$. To this end, we consider the simplest case where both predictors and error term are AR(1) processes with autoregression coefficient $\phi$. The GLS filter then yields $\widetilde{x}_{it}=x_{it}-\phi x_{it-1}=u_{it}$. 

However, also in this case two main differences emerge between the procedures. First, GLS-LASSO requires one more estimation step compared to ARMAr-LASSO. In step 1 GLS-LASSO estimates $\varepsilon_t$ by means of classical LASSO applied directly on time series, which we know to be a non-optimal procedure for the LASSO for the problems listed so far. In particular, without removing residuals serial correlation the variance of $T^{-1}\mathbf{x}_i'\pmb{\varepsilon}$ depends on both $\phi$ and $\phi_{\varepsilon}$ also after the standardization of $\mathbf{x}_i$. The sample cross-product $T^{-1}\sum_tx_{i,t}\varepsilon_t$ can fluctuate under serial dependence. Therefore estimates of $\varepsilon_t$ can be problematic in finite samples. Second, GLS-LASSO has poor forecasting performance compared to ARMAr-LASSO. GLS-LASSO reduces the explained variance of $y_t$ compared to ARMAr-LASSO since it does not consider the past of $y_t$. This can be mitigated by considering the term $\widehat{\phi}y_t$ in the forecasting equation, but $\widehat{\phi}$ obtained at step 2 of GLS-LASSO  is affected by estimation issues due to the estimate of $\varepsilon_t$ at step 1.

\bigskip

\noindent{\it Out of the Common Factor Restriction.}
Here we consider the case where $\phi\neq\phi_{\varepsilon}$, namely, all predictors have the same autoregressive coefficient, which differs from that of the error term. Without loss of generality, we note that in this case $\widetilde{x}_{it}=x_{it}-\phi_{\varepsilon}x_{it-1}=(\phi-\phi_{\varepsilon})x_{it-1}+u_{it}$ exhibits the following variance
\[\frac{(1-2\phi\phi_{\varepsilon}+\phi_{\varepsilon}^2)\sigma_{u_i}^2}{1-\phi^2},\]
which corresponds to the variance of an ARMA(1,1) with AR coefficient $\phi$ and MA coefficient $-\phi_{\varepsilon}$. Thus, $\widetilde x_{i,t}\neq u_{i,t}$ when $\phi\neq\phi_\varepsilon$, and GLS filtering generally leaves residual serial dependence in the predictors. Its implications for minimum eigenvalues and LASSO performance depend on the DGP. Numerical comparisons are reported in Supplement~\ref{Sec:MinEig}.

\section{$\widehat{u}$'s Estimation Error}
\begin{figure}[t]
\begin{center}
\graphicspath{{images/}}
\includegraphics[scale=0.6]{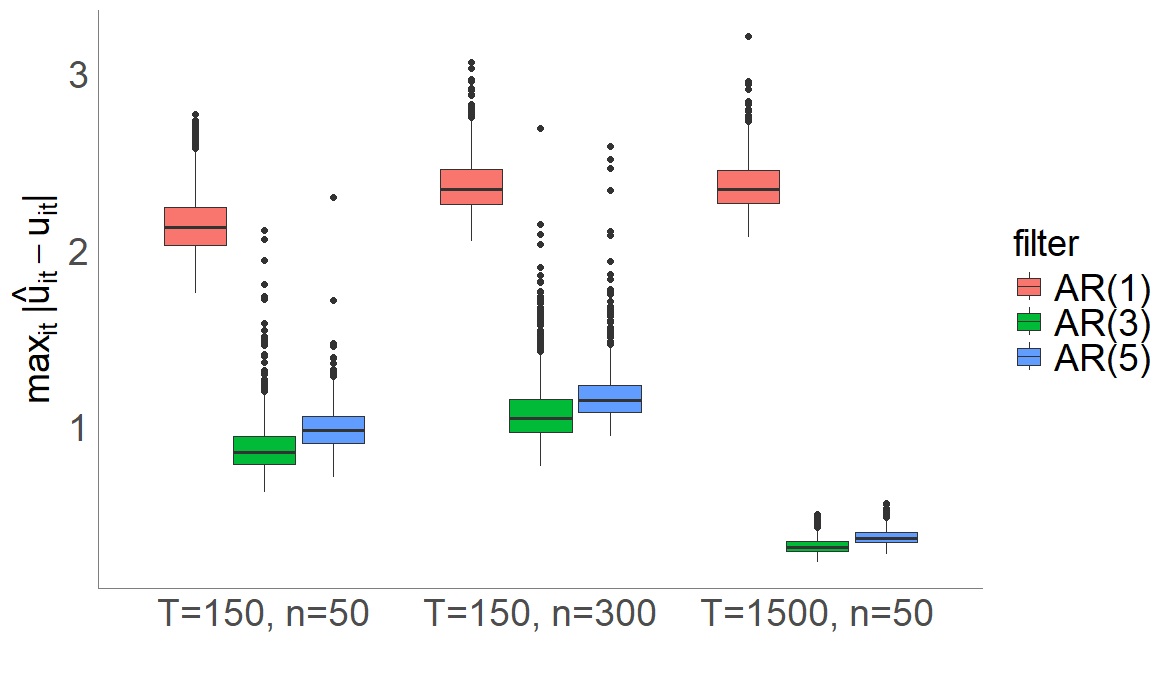}
\end{center}
\caption{\footnotesize Maximum AR residuals estimation error for different $n$, $T$ and filter settings.}
\label{fig:UError}
\end{figure}

We generated n variables from an AR(3) process and applied three different filters. In the first case, we fitted an AR(1) (underspecified order) model to each variable and used the corresponding parameter for filtering. In the second and third cases, we repeated the same exercise fitting AR(3) (correct order) and AR(5) (overspecified order) processes, respectively. For each of the three scenarios, we reported $\aggregate{max}{i\leq i\leq n,1\leq t\leq T}|\widehat{u}_{i,t}-u_{i,t}|$ under three regimes: $n<T$ (classical setting), $n>T$ (high-dimensional setting), and $n\ll T$ (asymptotic setting). Results, reported in Figure~\ref{fig:UError}, are obtained on 1000 Monte Carlo replications. Correctly specified or overfitted filters show decreasing innovation errors in the simulated settings.

\section{More Experiments}

\subsection{Filters Evaluation}\label{sec:FiltersExp}
We consider the following DGP:
    $$x_{i,t}=0.45x_{i,t-1}+0.45x_{i,t-2}+u_{i,t},\ \varepsilon_t=0.45\varepsilon_{t-1}+0.45\varepsilon_{t-2}+\omega_t\ \ .$$
We test several ARMAr-LASSO settings in terms of the number of lags of $y_t$ included in the model ($p_y$), the order of the AR filter used to obtain the estimated $\widehat{u}$'s ($p_i$), and the model selection method for the filter. We consider as reference setting the case where $p_y=3$ and the $\widehat{u}$'s are obtained by filtering each variable with an AR($p_i$), where $p_i$ (max 3) is selected via BIC and the AR parameters are estimated via ML (see Section~\ref{ARMArLASSOMonteCarlo}). The dimensionality is kept fixed at $n = 150$, while the sample size varies; we consider
$T = 75, 150, 300$. Both CoEr and RMSFE are relative to the working model~\ref{WM}, which includes the true $u$’s. Results are summarized in Table~\ref{Tab:FiltersExpBoth}. To assess the impact of including lags of $y_t$, we exclude them from the penalization process and compare the cases with one lag ($y_{t-1}$) and three lags ($y_{t-1}^{t-3}$) treated as fixed predictors in the model. Note that under DGP (A), the number of lags of $y_t$ that 
yields white-noise error terms is $p_y=2$. Using too few lags ($y_{t-1}$) reduces both estimation and forecast accuracy, whereas using an excessive number of lags ($y_{t-1}^{t-3}$) does not provide additional benefits since, especially for large $T$, the performance coincides with that of ARMAr-LS. To evaluate the effect of filter misspecification, instead of selecting $p_i$ using BIC, we filter the variables with AR(1), AR(2), and AR(3) processes, thus under-specifying, correctly-specifying, and over-specifying the true order, respectively. Underestimating the true order worsens coefficient estimation and forecast accuracy, particularly at high SNRs, while over-specifying the true order provides no gains. Finally, we evaluate alternative model selection criteria
other than BIC. In particular, we
consider the Akaike Information Criterion (AIC), its small-sample corrected version (AICc), and 
block cross-validation (blCV), which is a time-series variant of cross-validation~\citep{racine1997}. For blCV, we use a rolling window of size $T-3$, and for each variable, we select the order $p_i$ (max 3) that best predicts one step over the three out-of-sample horizons. All methods yield similar results, except for blCV, which performs worse than the others for high SNRs and large $T$.
\begin{table}
\centering
\caption{\footnotesize DGP (A). CoEr, RMSFE (relative to WM), \%TP and \%FP for various ARMAr-LASSO settings. 
For each $T$ setting the best CoEr and RMSFE are in bold.}
\label{Tab:FiltersExpBoth}
\makebox[\linewidth]{
\scalebox{0.55}{
\begin{tabular}{@{}lcrrrr@{\hspace{0.75cm}}rrrr@{\hspace{0.75cm}}rrrr@{}}
\hline
\multicolumn{13}{c}{\textbf{SNR = 1}}\\
\hline
&$T$&\multicolumn{4}{c}{75}&\multicolumn{4}{c}{150}&\multicolumn{4}{c}{300}\\\hline
& &CoEr&RMSFE& TP& FP
   &CoEr&RMSFE& TP& FP
   &CoEr&RMSFE& TP& FP\\
\hline
ARMAr-LAS       &&1.06 &1.02 &0.49 &0.13  &\bf1.02 &\bf1.00 &0.61 &0.01  &\bf1.00 &\bf1.00 &0.77 &0.01\\
$y_{t-1}$       &&1.30 &1.20 &0.49 &0.15  &1.08    &1.21    &0.61 &0.01  &1.08    &1.24    &0.77 &0.01\\
$y_{t-1}^{(3)}$ &&1.06 &1.02 &0.44 &0.18  &\bf1.02 &\bf1.00 &0.52 &0.01  &\bf1.00 &\bf1.00 &0.69 &0.01\\
AR(1)           &&1.02 &1.00 &0.49 &0.16  &1.03    &1.05    &0.61 &0.01  &1.02    &1.06    &0.75 &0.01\\
AR(2)           &&1.10 &1.02 &0.48 &0.15  &\bf1.02 &\bf1.00 &0.57 &0.02  &\bf1.00 &\bf1.00 &0.72 &0.01\\
AR(3)           &&1.09 &1.00 &0.49 &0.15  &\bf1.02 &\bf1.00 &0.61 &0.01  &\bf1.00 &\bf1.00 &0.77 &0.01\\
AIC             &&1.06 &1.01 &0.49 &0.15  &\bf1.02 &\bf1.00 &0.61 &0.01  &\bf1.00 &\bf1.00 &0.77 &0.01\\
AICc            &&\bf1.05 &\bf0.99 &0.49 &0.14  &\bf1.02 &\bf1.00 &0.61 &0.01  &\bf1.00 &\bf1.00 &0.77 &0.01\\
blCV            &&1.07 &1.00 &0.49 &0.14  &\bf1.02 &1.01    &0.61 &0.01  &1.02    &1.01    &0.77 &0.01\\
\hline\hline
\multicolumn{13}{c}{\textbf{SNR = 10}}\\
\hline
&$T$&\multicolumn{4}{c}{75}&\multicolumn{4}{c}{150}&\multicolumn{4}{c}{300}\\\hline
& &CoEr&RMSFE& TP& FP
   &CoEr&RMSFE& TP& FP
   &CoEr&RMSFE& TP& FP\\
\hline
ARMAr-LAS       &&\bf1.05 &1.03 &0.87 &0.05  &\bf1.03 &\bf1.01 &0.97 &0.02  &\bf1.00 &\bf1.00 &1.00 &0.01\\
$y_{t-1}$       &&1.69    &1.91 &0.86 &0.05  &1.63    &1.91    &0.97 &0.02  &1.81    &1.86    &1.00 &0.01\\
$y_{t-1}^{(3)}$ &&\bf1.05 &1.03 &0.65 &0.09  &\bf1.03 &\bf1.01 &0.80 &0.01  &\bf1.00 &\bf1.00 &0.93 &0.01\\
AR(1)           &&1.13    &1.04 &0.85 &0.06  &1.17    &1.12    &0.96 &0.02  &1.19    &1.10    &1.00 &0.01\\
AR(2)           &&\bf1.05 &\bf1.01 &0.84 &0.07  &\bf1.03 &\bf1.01 &0.95 &0.03  &\bf1.00 &\bf1.00 &0.99 &0.02\\
AR(3)           &&1.07    &\bf1.01 &0.86 &0.05  &\bf1.03 &\bf1.01 &0.97 &0.02  &\bf1.00 &\bf1.00 &1.00 &0.01\\
AIC             &&1.06    &1.02 &0.86 &0.05  &\bf1.03 &1.02    &0.97 &0.02  &\bf1.00 &\bf1.00 &1.00 &0.01\\
AICc            &&\bf1.05 &1.02 &0.85 &0.05  &\bf1.03 &1.02    &0.97 &0.02  &\bf1.00 &\bf1.00 &1.00 &0.01\\
blCV            &&1.07    &1.04 &0.85 &0.05  &1.10    &1.09    &0.97 &0.02  &1.15    &1.07    &1.00 &0.01\\
\hline
\end{tabular}
}}
\end{table}

\subsection{Simulation
Experiments with More DGPs}\label{ARMArLASSOMonteCarlo_supp}
The response variable 
is generated using the model $y_t=\sum_{i=1}^n\alpha_i^*x_{i,t-1}+\varepsilon_t$, and we consider the following data generating processes (DGPs) for 
predictors and error terms:
\setcounter{bean}{0}
\begin{list}
{(\Alph{bean})}{\usecounter{bean}}    \item Common AR(1) Restriction: $x_{i,t}=\phi x_{i,t-1}+u_{i,t},\ \varepsilon_t=\phi\varepsilon_{t-1}+\omega_t$.
    \item Common AR(1) Restriction with Common Factor: $x_{i,t}= f_t+z_{i,t}$, where $f_t=\phi f_{t-1}+\delta_t,$ $z_{i,t}=\phi z_{i,t-1}+\eta_{i,t},$ $\varepsilon_t=\phi\varepsilon_{t-1}+\omega_t$.
\end{list}

\noindent The shocks are generated as follows: $u_{i,t} \sim i.i.d.\ \mathcal{N}(0,1)$ with $\left(\mathbf{C}_{u}\right){ij}=c_{ij}^{u}=\rho^{|i-j|}$, $\delta_t,\ \eta_{i,t} \sim i.i.d.\ \mathcal{N}(0,1)$ with $\left(\mathbf{C}_{\eta}\right){ij}=c_{ij}^{\eta}=\rho^{|i-j|}$, and $\omega_t \sim i.i.d.\ \mathcal{N}(0,\sigma_{\omega}^2)$. For the DGP in (A) we set $\rho=0.8$, while for the DGP in (B) we set $\rho=0.4$ to generate predictors primarily influenced by the common factor, with weakly correlated AR and/or ARMA idiosyncratic components. Finally, we vary the value of $\sigma_{\omega}^2$ to explore different signal-to-noise 
ratios (SNRs).

We compare our ARMAr-LASSO (ARMAr-LAS) with the LASSO-based benchmarks employed in Section~\ref{ARMArLASSOMonteCarlo} of the main text.

For the DGP in (A) 
we set $\rho=0.8$, while for the DGP in (B) 
we set $\rho=0.4$ to generate predictors primarily influenced by the common factor, with weakly correlated AR or ARMA idiosyncratic components. Finally, we vary the value of $\sigma_{\omega}^2$ to explore different signal-to-noise ratios (SNRs).

We compare our ARMAr-LASSO (ARMAr-LAS) with the standard LASSO applied to the observed time series (LAS), LASSO applied to the observed time series plus lags of $y_t$ (LASy), GLS-LASSO as proposed by~\cite{chronopoulos2023} (GLS-LAS), autoregressive distributed lag LASSO (ARDL-LAS), and FarmSelect as proposed by~\cite{Fan2020}, which employs LASSO on factor model residuals (FaSel). The performance of each method is evaluated based on average results from 1000 independent simulations, focusing on the coefficient estimation error (CoEr) obtained as $||\widehat{\pmb\alpha}-\pmb\alpha||_2$, the Root Mean Square Forecast Error (RMSFE), and the percentages of true positives (\%TP) and false positives (\%FP) in selecting relevant predictors. Simulations have varying numbers of predictors (dimensionality), $n=50, 150, 300$, and a fixed sample size, $T=150$. In this way we cover low ($n=50$), intermediate ($n=150$), and high ($n=300$) dimensional scenarios. For all methods, the tuning parameter $\lambda$ is selected using the Bayesian Information Criterion (BIC). Finally, regardless of the choice of $n$, $\pmb\alpha^*$ is always taken to have the first $10$ entries equal to $1$ and all others equal to $0$. In this way, as $n$ varies, we also cover different levels of sparsity. 

For DGP (A), we investigate settings with different $\phi$ ($0.3, 0.6, 0.9, 0.95$) and different SNR ($0.5, 1, 5, 10$).  For GLS-LAS, we estimate an AR(1) model on $\widehat{\varepsilon}_t$ (see Supplement~\ref{COMP_ARDLandGLS}) and use the resulting autoregressive coefficient to filter both response and predictors. For ARDL-LAS, we consider one lag of the response and one lag of each predictor as additional regressors, bringing the number of terms undergoing selection to $n\times2+1$. For the working model underlying ARMAr-LAS, the $\widehat{u}$'s are obtained by filtering each series with an AR(1) process, and we consider $p_y=1$; that is, we take one lag of $y_t$ as additional predictor. Results are presented in Table~\ref{TabDGPA&B} for SNR values of 1 and 10 (complete results are provided in Supplement~\ref{ARMArLASSOMonteCarlo_supp}). For each SNR, CoEr and RMSFE (both expressed in relative terms to the values obtained by LAS), as well as \%TP and \%FP are given for every $n$ and $\phi$ considered (the best CoEr and RMSFE are in bold). Results have ARMAr-LAS as the best performer in terms of CoEr and RMSFE across values of $\phi$, $n$, and SNR, demonstrating superior accuracy in both estimation and forecasting compared to the other LASSO-based methods considered. ARMAr-LAS also shows superior performance in feature selection, with a higher \%TP and a lower \%FP. These gains are more evident when serial correlations are stronger ($\phi=0.6$ or higher).

Notably, under the common AR(1) restriction, the ARMAr-LAS and GLS-LAS estimators should be equivalent (this is the one case where the GLS-LAS estimator removes the serial correlations of the predictors). Nevertheless, GLS-LAS performs on par with ARMAr-LAS only when serial correlations are low; ARMAr-LAS outperforms GLS-LAS when serial correlations are medium/high, likely because the latter requires the estimation of $\widehat{\varepsilon}_t$ (see Supplement~\ref{COMP_ARDLandGLS}). Also, in some instances, ARDL-LAS exhibits a slightly lower \%FP than ARMAr-LAS. However, this metric is calculated on $n\times2+1$ predictors for the former; in terms of the absolute number of false positives, ARDL-LAS has more than ARMAr-LAS (see Supplement~\ref{COMP_ARDLandGLS}).

Finally, we note that the superior performance of ARMAr-LAS in 
DGP (B)
indicates its effectiveness in handling 
factor structures, where 
multicollinearities 
are more complex than for simple AR processes (DGP (A)).
\begin{table}[t]
\centering
\caption{\footnotesize DGPs (A) and (B). CoEr, RMSFE (relative to LAS), \%TP and \%FP for LASSO-based benchmarks and ARMAr-LASSO. For each $n$ and $\phi$ setting the best CoEr and RMSFE are in bold. 
}\label{TabDGPA&B}
    \makebox[\linewidth]{
\scalebox{0.4}{
\begin{tabular}
{@{}cccccrrrr@{\hspace{0.3cm}}rrrrr@{\hspace{0.3cm}}rrrrrr
@{\hspace{0.9cm}}rrrr@{\hspace{0.3cm}}rrrrr@{\hspace{0.3cm}}rrrrrr@{}}
\hline\hline
&&&&& \multicolumn{14}{c}{(A)} && \multicolumn{14}{c}{(B)}\\\hline
&&&n&&\multicolumn{4}{c}{75}&&\multicolumn{4}{c}{150}&&\multicolumn{4}{c}{300}&&\multicolumn{4}{c}{50}&&\multicolumn{4}{c}{150}&&\multicolumn{4}{c}{300}\\\hline
&&&$\phi$& &0.3 &0.6 &0.9 &0.95 &&0.3 &0.6 &0.9 &0.95 &&0.3 &0.6 &0.9 &0.95 & &0.3 &0.6 &0.9 &0.95 &&0.3 &0.6 &0.9 &0.95 &&0.3 &0.6 &0.9 &0.95\\\hline
SNR&&&&&&&&&&&&&&&&&& &&&&&&&&&&&&&&&\\\hline
0.5&&&&&&&&&&&&&&&&&&& &&&&&&&&&&&&&&&\\
&CoEr&&&&&&&&&&&&&&&&&& &&&&&&&&&&&&&&&\\\cline{2-2}
&&    LASSOy       &&&    0.99 &  0.81 &  0.37 &  0.36   &&   0.99 &  0.82 &  0.46 &  0.46   &&   0.99 &  0.83 &  0.53 &  0.53     &&   0.99 &  0.83 &  0.41 &  0.39   &&  0.99 &  0.85 &  0.47 &  0.46  &&   1.00 &  0.86 &  0.53 &  0.51                  \\
&&    GLS-LASSO    &&&    0.96 &  0.82 &  0.64 &  0.70   &&   0.97 &  0.81 &  0.75 &  0.79   &&   0.97 &  0.83 &  0.82 &  0.85     &&   0.96 &  0.79 &  0.64 &  0.73   &&  0.96 &  0.80 &  0.76 &  0.79  &&   0.96 &  0.82 &  0.82 &  0.84                  \\
&&    ARDL-LAS     &&&    0.96 &  0.79 &  0.34 &  0.32   &&   0.97 &  0.80 &  0.42 &  0.42   &&   1.32 &  0.83 &  0.50 &  0.48     &&   0.97 &  0.81 &  0.38 &  0.35   &&  0.99 &  0.83 &  0.44 &  0.42  &&   1.46 &  0.88 &  0.50 &  0.48                  \\
&&    FaSel        &&&    1.01 &  0.99 &  1.04 &  1.00   &&   1.04 &  1.01 &  1.00 &  1.00   &&   2.38 &  1.60 &  1.01 &  0.99     &&   0.91 &  0.86 &  1.01 &  1.02   &&  0.89 &  0.85 &  1.03 &  1.04  &&   5.04 &  2.81 &  1.04 &  1.04                  \\
&&    ARMAr-LAS    &&&    0.98 &  0.82 &  0.34 &  0.31   &&   0.99 &  0.82 &  0.41 &  0.41   &&   0.99 &  0.82 &  0.47 &  0.45     &&   0.99 &  0.82 &  0.37 &  0.34   &&  0.99 &  0.83 &  0.43 &  0.41  &&   1.01 &  0.83 &  0.48 &  0.46                 \\
&RMSFE&&&&&&&&&&&&&&&&&& &&&&&&&&&&&&&&&\\\cline{2-2}
&&    LASSOy      &&&     0.98 &  0.89 &  0.77 &  0.76   &&   0.99 &  0.88 &  0.82 &  0.78   &&   0.98 &  0.90 &  0.84 &  0.81     &&   0.99 &  0.88 &  0.76 &  0.72   &&  1.00 &  0.89 &  0.79 &  0.79  &&   0.99 &  0.89 &  0.85 &  0.82                 \\
&&    GLS-LASSO   &&&     0.97 &  0.83 &  0.84 &  0.86   &&   0.96 &  0.83 &  0.86 &  0.88   &&   0.95 &  0.85 &  0.91 &  0.91     &&   0.96 &  0.83 &  0.82 &  0.85   &&  0.98 &  0.84 &  0.90 &  0.89  &&   0.97 &  0.86 &  0.93 &  0.91                 \\
&&    ARDL-LAS    &&&     1.00 &  0.89 &  0.77 &  0.75   &&   1.00 &  0.89 &  0.82 &  0.78   &&   1.05 &  0.91 &  0.84 &  0.82     &&   0.99 &  0.88 &  0.72 &  0.70   &&  1.00 &  0.89 &  0.77 &  0.77  &&   1.05 &  0.88 &  0.82 &  0.80                 \\
&&    FaSel       &&&     0.96 &  0.97 &  0.96 &  0.96   &&   0.96 &  0.99 &  0.97 &  0.94   &&   0.90 &  0.92 &  0.97 &  0.97     &&   0.97 &  1.00 &  0.99 &  0.99   &&  0.98 &  1.00 &  0.97 &  0.98  &&   1.29 &  1.07 &  0.99 &  0.94                 \\
&&    ARMAr-LAS   &&&     0.98 &  0.83 &  0.69 &  0.68   &&   0.98 &  0.84 &  0.73 &  0.72   &&   0.98 &  0.86 &  0.75 &  0.74     &&   0.97 &  0.83 &  0.67 &  0.64   &&  0.99 &  0.83 &  0.72 &  0.72  &&   1.01 &  0.86 &  0.77 &  0.77                 \\
& TP&&&&&&&&&&&&&&&&&& &&&&&&&&&&&&&&&\\\cline{2-2}
&&    LASSO      &&&     44.80 & 41.40 & 47.50 & 48.60   &&  42.60 & 38.60 & 43.70 & 43.80   &&  40.60 & 37.10 & 42.30 & 40.40     &&  34.50 & 30.40 & 42.00 & 43.50   && 28.50 & 23.30 & 33.00 & 30.00  &&  20.90 & 18.70 & 24.10 & 23.40               \\
&&    LASSOy     &&&     45.40 & 42.10 & 32.40 & 31.90   &&  43.70 & 39.00 & 31.50 & 31.00   &&  41.40 & 37.80 & 30.40 & 29.40     &&  35.30 & 32.00 & 21.70 & 22.50   && 28.90 & 25.00 & 17.70 & 16.40  &&  21.50 & 20.00 & 13.20 & 13.20               \\
&&    GLS-LASSO  &&&     45.70 & 46.20 & 43.10 & 43.80   &&  44.20 & 43.50 & 41.40 & 41.50   &&  41.90 & 42.10 & 40.50 & 39.00     &&  36.20 & 35.10 & 32.00 & 34.80   && 28.60 & 27.70 & 28.00 & 26.50  &&  21.50 & 21.20 & 21.70 & 20.80               \\
&&    ARDL-LAS   &&&     43.20 & 40.20 & 33.20 & 33.30   &&  41.70 & 37.30 & 31.80 & 30.90   &&  41.90 & 36.20 & 28.80 & 27.70     &&  33.90 & 31.60 & 25.30 & 24.60   && 28.00 & 24.80 & 20.30 & 18.60  &&  28.30 & 20.80 & 14.10 & 14.60               \\
&&    FaSel      &&&     11.80 & 18.00 & 39.80 & 42.00   &&  19.30 & 24.70 & 39.30 & 42.10   &&  39.90 & 35.00 & 38.90 & 40.30     &&   2.60 &  8.10 & 38.60 & 40.50   &&  2.70 &  7.10 & 30.30 & 29.70  &&  56.50 & 38.40 & 22.90 & 22.30               \\
&&    ARMAr-LAS  &&&     48.10 & 49.90 & 49.10 & 48.90   &&  46.20 & 47.50 & 47.80 & 48.10   &&  42.80 & 45.10 & 45.20 & 44.20     &&  38.30 & 40.80 & 40.50 & 40.80   && 30.30 & 33.90 & 33.80 & 33.60  &&  22.60 & 26.20 & 27.00 & 27.80               \\
& FP&&&&&&&&&&&&&&&&&& &&&&&&&&&&&&&&&\\\cline{2-2}
&&    LASSO      &&&     1.50 &  4.90 & 31.40 & 34.20   &&   0.70 &  3.00 & 17.90 & 17.60   &&   0.30 &  1.90 &  9.90 &  9.80     &&   6.30 &  9.30 & 32.50 & 34.20   &&  3.80 &  5.50 & 18.90 & 18.30  &&   2.20 &  3.50 & 10.40 &  9.90                 \\
&&    LASSOy     &&&     1.60 &  2.80 &  7.80 & 10.60   &&   0.80 &  1.50 &  5.10 &  5.90   &&   0.40 &  0.80 &  3.00 &  3.60     &&   6.70 &  8.00 & 10.20 & 12.20   &&  3.90 &  4.80 &  6.30 &  6.70  &&   2.30 &  2.90 &  3.80 &  3.90                 \\
&&    GLS-LASSO  &&&     1.00 &  1.30 & 16.90 & 22.00   &&   0.50 &  0.60 & 11.10 & 12.10   &&   0.20 &  0.50 &  6.90 &  7.30     &&   5.90 &  6.20 & 17.50 & 22.80   &&  3.40 &  3.50 & 12.20 & 13.00  &&   2.00 &  2.30 &  7.40 &  7.50                 \\
&&    ARDL-LAS   &&&     0.80 &  1.40 &  5.00 &  6.00   &&   0.50 &  0.70 &  2.90 &  3.50   &&   4.60 &  0.80 &  1.60 &  1.90     &&   3.00 &  4.40 &  6.90 &  7.50   &&  2.30 &  2.60 &  3.90 &  4.10  &&   8.10 &  2.40 &  2.20 &  2.30                 \\
&&    FaSel      &&&     0.90 &  3.80 & 29.20 & 31.50   &&   0.60 &  2.70 & 16.60 & 17.10   &&  10.70 &  7.50 &  9.50 &  9.60     &&   0.30 &  3.10 & 31.80 & 33.80   &&  0.20 &  1.80 & 19.10 & 19.80  &&  48.30 & 29.80 & 11.30 & 10.80                 \\
&&    ARMAr-LAS  &&&     1.80 &  1.80 &  2.00 &  2.00   &&   0.80 &  0.90 &  0.90 &  1.00   &&   0.40 &  0.50 &  0.50 &  0.50     &&   7.20 &  8.20 &  8.00 &  8.20   &&  4.10 &  5.00 &  5.00 &  4.90  &&   2.60 &  3.00 &  2.90 &  2.90                 \\\cline{2-34}

1&&&&&&&&&&&&&&&&&&& &&&&&&&&&&&&&&&\\
&CoEr&&&&&&&&&&&&&&&&&& &&&&&&&&&&&&&&&\\\cline{2-2}
&&    LASSOy       &&&   0.99 &  0.85 &  0.51 &  0.50   &&   0.99 &  0.86 &  0.61 &  0.59   &&   0.99 &  0.87 &  0.66 &  0.66     &&   0.99 &  0.87 &  0.56 &  0.54   &&   0.99 &  0.89 &  0.63 &  0.62   &&  0.99 &  0.90 &  0.69 &  0.66                \\
&&    GLS-LASSO    &&&   0.97 &  0.81 &  0.70 &  0.74   &&   0.97 &  0.83 &  0.80 &  0.83   &&   0.97 &  0.83 &  0.86 &  0.88     &&   0.96 &  0.81 &  0.71 &  0.75   &&   0.96 &  0.82 &  0.80 &  0.82   &&  0.96 &  0.84 &  0.85 &  0.85                \\
&&    ARDL-LAS     &&&   0.97 &  0.83 &  0.44 &  0.41   &&   0.98 &  0.84 &  0.55 &  0.54   &&   1.01 &  0.87 &  0.61 &  0.60     &&   0.98 &  0.85 &  0.49 &  0.45   &&   0.98 &  0.86 &  0.56 &  0.55   &&  1.07 &  0.89 &  0.63 &  0.60                \\
&&    FaSel        &&&   1.10 &  1.05 &  1.04 &  1.00   &&   1.08 &  1.04 &  1.02 &  1.02   &&   1.93 &  1.34 &  1.00 &  0.98     &&   1.01 &  0.92 &  1.03 &  1.02   &&   0.98 &  0.90 &  1.04 &  1.06   &&  3.86 &  2.29 &  1.04 &  1.05                \\
&&    ARMAr-LAS    &&&   0.98 &  0.82 &  0.41 &  0.38   &&   0.98 &  0.83 &  0.50 &  0.49   &&   0.98 &  0.83 &  0.55 &  0.53     &&   0.98 &  0.81 &  0.44 &  0.41   &&   0.99 &  0.83 &  0.50 &  0.50   &&  1.00 &  0.84 &  0.57 &  0.54                \\
&RMSFE&&&&&&&&&&&&&&&&&& &&&&&&&&&&&&&&&\\\cline{2-2}
&&    LASSOy      &&&    0.98 &  0.90 &  0.81 &  0.76   &&   0.99 &  0.92 &  0.87 &  0.79   &&   1.00 &  0.90 &  0.89 &  0.85     &&   0.99 &  0.90 &  0.81 &  0.83   &&   0.99 &  0.93 &  0.83 &  0.83   &&  1.00 &  0.92 &  0.88 &  0.85                 \\
&&    GLS-LASSO   &&&    0.95 &  0.82 &  0.84 &  0.86   &&   0.95 &  0.85 &  0.88 &  0.87   &&   0.96 &  0.82 &  0.94 &  0.92     &&   0.96 &  0.85 &  0.82 &  0.91   &&   0.96 &  0.87 &  0.89 &  0.89   &&  0.95 &  0.86 &  0.92 &  0.90                 \\
&&    ARDL-LAS    &&&    0.99 &  0.91 &  0.79 &  0.74   &&   1.00 &  0.92 &  0.86 &  0.79   &&   1.01 &  0.91 &  0.88 &  0.83     &&   1.00 &  0.89 &  0.76 &  0.78   &&   0.99 &  0.91 &  0.77 &  0.80   &&  1.01 &  0.92 &  0.84 &  0.81                  \\
&&    FaSel       &&&    0.96 &  0.98 &  0.97 &  0.95   &&   0.97 &  0.98 &  0.96 &  0.95   &&   0.90 &  0.95 &  0.96 &  0.96     &&   1.00 &  1.03 &  0.99 &  0.99   &&   0.99 &  1.02 &  0.98 &  0.94   &&  1.16 &  1.05 &  0.99 &  0.96                 \\
&&    ARMAr-LAS   &&&    0.96 &  0.82 &  0.67 &  0.65   &&   0.97 &  0.85 &  0.72 &  0.71   &&   0.98 &  0.84 &  0.76 &  0.73     &&   0.97 &  0.86 &  0.68 &  0.71   &&   0.98 &  0.85 &  0.72 &  0.75   &&  0.99 &  0.85 &  0.77 &  0.76                 \\
& TP&&&&&&&&&&&&&&&&&& &&&&&&&&&&&&&&&\\\cline{2-2}
&&    LASSO      &&&     60.10 & 54.00 & 55.80 & 56.30   &&  58.40 & 52.10 & 53.50 & 51.70   &&  56.70 & 51.10 & 52.30 & 50.20     &&  53.90 & 45.30 & 49.00 & 50.80   &&  47.60 & 39.10 & 41.10 & 38.60   && 40.30 & 33.70 & 33.70 & 32.50                \\
&&    LASSOy     &&&     60.50 & 54.90 & 46.10 & 46.10   &&  58.70 & 52.80 & 46.00 & 43.90   &&  56.80 & 51.90 & 44.80 & 42.50     &&  54.20 & 46.60 & 35.40 & 36.40   &&  47.80 & 40.30 & 30.50 & 28.00   && 40.60 & 34.50 & 25.10 & 24.20                \\
&&    GLS-LASSO  &&&     62.10 & 61.90 & 54.10 & 53.40   &&  60.20 & 58.80 & 52.10 & 50.60   &&  58.80 & 58.60 & 51.50 & 49.40     &&  56.30 & 55.20 & 43.90 & 44.50   &&  49.40 & 48.50 & 38.80 & 35.20   && 41.90 & 42.30 & 31.90 & 30.90                \\
&&    ARDL-LAS   &&&     59.10 & 53.60 & 49.20 & 48.00   &&  57.30 & 51.50 & 45.70 & 44.20   &&  55.70 & 50.30 & 44.20 & 41.70     &&  53.20 & 47.30 & 39.30 & 38.50   &&  47.00 & 41.10 & 33.40 & 30.00   && 40.90 & 34.70 & 26.50 & 26.20                \\
&&    FaSel      &&&     19.10 & 26.20 & 45.60 & 49.00   &&  30.70 & 36.30 & 48.00 & 50.60   &&  51.90 & 47.80 & 50.20 & 50.80     &&   7.40 & 15.40 & 45.30 & 47.20   &&   6.50 & 13.00 & 38.50 & 37.30   && 60.60 & 44.90 & 32.30 & 32.30                \\
&&    ARMAr-LAS  &&&     63.40 & 64.00 & 63.40 & 64.00   &&  61.30 & 60.80 & 61.90 & 61.70   &&  59.90 & 60.00 & 60.20 & 59.70     &&  57.80 & 57.60 & 58.50 & 57.40   &&  50.50 & 51.20 & 51.30 & 50.20   && 43.20 & 44.80 & 44.50 & 45.00                \\
& FP&&&&&&&&&&&&&&&&&& &&&&&&&&&&&&&&&\\\cline{2-2}
&&    LASSO      &&&     1.80 &  5.80 & 31.70 & 33.60   &&   0.90 &  3.40 & 17.40 & 16.80   &&   0.40 &  2.10 &  9.80 &  9.40     &&   8.30 & 11.10 & 32.50 & 34.40   &&   5.20 &  7.10 & 18.30 & 17.60   &&  3.30 &  4.40 & 10.20 &  9.70               \\
&&    LASSOy     &&&     1.80 &  3.60 & 12.80 & 15.90   &&   0.90 &  1.80 &  7.70 &  8.10   &&   0.50 &  1.10 &  4.40 &  4.90     &&   8.40 &  9.80 & 15.50 & 17.60   &&   5.30 &  6.20 &  9.10 &  9.30   &&  3.30 &  3.80 &  5.40 &  5.30               \\
&&    GLS-LASSO  &&&     1.30 &  1.70 & 17.90 & 22.10   &&   0.60 &  0.70 & 11.00 & 11.70   &&   0.30 &  0.50 &  6.80 &  7.20     &&   8.00 &  8.50 & 19.20 & 22.90   &&   5.00 &  5.30 & 12.50 & 12.50   &&  3.10 &  3.40 &  7.50 &  7.20               \\
&&    ARDL-LAS   &&&     0.80 &  2.00 &  7.80 &  8.50   &&   0.40 &  0.90 &  4.60 &  4.80   &&   0.70 &  0.80 &  2.50 &  2.80     &&   3.90 &  5.70 &  9.20 &  9.30   &&   2.60 &  3.40 &  5.20 &  5.30   &&  3.50 &  2.30 &  3.00 &  3.00               \\
&&    FaSel      &&&     1.40 &  5.20 & 29.40 & 30.70   &&   0.90 &  3.30 & 16.60 & 16.90   &&   9.30 &  6.20 &  9.30 &  9.20     &&   0.40 &  4.10 & 33.10 & 34.30   &&   0.20 &  2.30 & 19.80 & 19.80   && 46.40 & 28.60 & 11.40 & 11.00               \\
&&    ARMAr-LAS  &&&     2.10 &  2.30 &  2.20 &  2.40   &&   1.00 &  1.10 &  1.20 &  1.20   &&   0.50 &  0.60 &  0.60 &  0.60     &&   9.40 &  9.60 &  9.80 &  9.50   &&   5.80 &  6.20 &  6.00 &  6.10   &&  3.70 &  3.70 &  3.80 &  3.60               \\\cline{2-34}

5&&&&&&&&&&&&&&&&&&& &&&&&&&&&&&&&&&\\
&CoEr&&&&&&&&&&&&&&&&&& &&&&&&&&&&&&&&&\\\cline{2-2}
&&    LASSOy       &&&     1.00 &  0.95 &  0.83 &  0.83  &&   1.00 &  0.96 &  0.88 &  0.87   &&   1.00 &  0.96 &  0.91 &  0.89   &&   1.00 &  0.96 &  0.87 &  0.85  &&   1.00 &  0.97 &  0.91 &  0.89   &&  1.00  &  0.98 &  0.93 &  0.92                 \\
&&    GLS-LASSO    &&&     0.95 &  0.81 &  0.79 &  0.82  &&   0.95 &  0.82 &  0.87 &  0.88   &&   0.96 &  0.83 &  0.92 &  0.92   &&   0.95 &  0.81 &  0.82 &  0.84  &&   0.95 &  0.83 &  0.87 &  0.86   &&  0.96  &  0.85 &  0.91 &  0.90                 \\
&&    ARDL-LAS     &&&     1.00 &  0.93 &  0.68 &  0.68  &&   1.00 &  0.95 &  0.81 &  0.81   &&   1.00 &  0.95 &  0.85 &  0.84   &&   0.99 &  0.93 &  0.70 &  0.70  &&   0.99 &  0.94 &  0.81 &  0.81   &&  0.99  &  0.95 &  0.85 &  0.85                 \\
&&    FaSel        &&&     1.45 &  1.27 &  1.06 &  1.07  &&   1.27 &  1.14 &  1.04 &  1.01   &&   1.42 &  1.15 &  0.99 &  0.98   &&   1.18 &  1.09 &  1.05 &  1.06  &&   1.14 &  1.06 &  1.07 &  1.06   &&  2.08  &  1.46 &  1.06 &  1.05                 \\
&&    ARMAr-LAS    &&&     0.95 &  0.80 &  0.50 &  0.49  &&   0.96 &  0.81 &  0.60 &  0.58   &&   0.97 &  0.81 &  0.64 &  0.61   &&   0.95 &  0.78 &  0.52 &  0.50  &&   0.96 &  0.79 &  0.60 &  0.58   &&  0.97  &  0.81 &  0.65 &  0.63                 \\
&RMSFE&&&&&&&&&&&&&&&&&& &&&&&&&&&&&&&&&\\\cline{2-2}
&&    LASSOy      &&&      1.00 &  0.97 &  0.93 &  0.92  &&   1.00 &  0.97 &  0.95 &  0.91   &&   1.00 &  0.98 &  0.96 &  0.95   &&   1.00 &  0.98 &  0.94 &  0.90  &&   1.00 &  0.97 &  0.96 &  0.91   &&  1.00  &  0.99 &  0.94 &  0.95                  \\
&&    GLS-LASSO   &&&      0.95 &  0.83 &  0.86 &  0.86  &&   0.96 &  0.84 &  0.88 &  0.83   &&   0.94 &  0.85 &  0.94 &  0.90   &&   0.95 &  0.84 &  0.88 &  0.87  &&   0.97 &  0.87 &  0.88 &  0.83   &&  0.96  &  0.91 &  0.90 &  0.87                  \\
&&    ARDL-LAS    &&&      1.00 &  0.98 &  0.94 &  0.93  &&   1.01 &  0.98 &  0.95 &  0.90   &&   1.01 &  0.99 &  0.99 &  0.96   &&   1.02 &  0.96 &  0.86 &  0.84  &&   1.01 &  0.97 &  0.90 &  0.84   &&  1.01  &  1.00 &  0.90 &  0.89                  \\
&&    FaSel       &&&      0.99 &  0.99 &  0.98 &  0.93  &&   0.99 &  0.96 &  0.95 &  0.93   &&   0.92 &  0.95 &  0.97 &  0.96   &&   1.05 &  1.03 &  0.97 &  0.96  &&   1.04 &  1.03 &  0.97 &  0.93   &&  1.13  &  1.09 &  0.97 &  0.91                  \\
&&    ARMAr-LAS   &&&      0.97 &  0.84 &  0.70 &  0.64  &&   0.97 &  0.84 &  0.70 &  0.68   &&   0.98 &  0.84 &  0.77 &  0.74   &&   0.97 &  0.82 &  0.71 &  0.65  &&   0.97 &  0.86 &  0.74 &  0.66   &&  0.96  &  0.89 &  0.74 &  0.69                  \\
& TP&&&&&&&&&&&&&&&&&& &&&&&&&&&&&&&&&\\\cline{2-2}
&&    LASSO      &&&      90.10 & 84.10 & 78.30 & 78.10  &&  89.10 & 82.80 & 79.60 & 75.70   &&  88.80 & 82.10 & 78.70 & 75.90   &&  90.20 & 83.30 & 76.30 & 75.80  &&  89.00 & 79.80 & 71.30 & 65.70   && 84.80  & 77.00 & 67.80 & 63.20                 \\
&&    LASSOy     &&&      90.10 & 84.30 & 76.70 & 77.00  &&  89.20 & 83.10 & 78.20 & 74.40   &&  88.80 & 82.60 & 77.60 & 74.10   &&  90.10 & 83.60 & 74.00 & 72.90  &&  88.90 & 79.90 & 69.70 & 63.60   && 84.90  & 77.10 & 65.70 & 60.70                 \\
&&    GLS-LASSO  &&&      91.70 & 90.50 & 79.60 & 78.50  &&  91.00 & 89.70 & 79.90 & 77.90   &&  90.40 & 89.00 & 78.70 & 76.50   &&  91.90 & 91.00 & 77.90 & 75.60  &&  90.70 & 87.90 & 74.10 & 70.10   && 87.00  & 84.70 & 70.50 & 66.90                 \\
&&    ARDL-LAS   &&&      89.20 & 83.40 & 76.50 & 74.40  &&  88.70 & 82.30 & 77.30 & 74.00   &&  88.00 & 81.80 & 76.80 & 73.20   &&  89.70 & 84.00 & 73.70 & 71.10  &&  88.50 & 80.10 & 71.00 & 65.90   && 84.50  & 77.00 & 66.60 & 62.40                 \\
&&    FaSel      &&&      40.60 & 50.40 & 66.10 & 71.50  &&  60.90 & 69.10 & 75.90 & 76.60   &&  81.90 & 78.70 & 78.30 & 77.30   &&  66.40 & 64.70 & 71.20 & 70.80  &&  65.40 & 61.40 & 69.90 & 67.90   && 77.30  & 70.50 & 67.20 & 65.60                 \\
&&    ARMAr-LAS  &&&      92.10 & 91.30 & 90.40 & 91.00  &&  91.50 & 90.70 & 90.50 & 90.20   &&  90.30 & 90.50 & 89.40 & 89.00   &&  92.60 & 92.20 & 90.90 & 90.30  &&  91.00 & 90.10 & 89.00 & 87.20   && 87.30  & 87.20 & 85.90 & 84.80                 \\
& FP&&&&&&&&&&&&&&&&&& &&&&&&&&&&&&&&&\\\cline{2-2}
&&    LASSO      &&&       2.50 &  7.30 & 30.90 & 32.00  &&   1.30 &  3.80 & 15.20 & 13.70   &&   0.60 &  2.50 &  8.50 &  8.30   &&  11.10 & 14.00 & 33.50 & 34.40  &&   7.00 &  9.50 & 16.40 & 14.50   &&  4.50  &  6.10 &  9.50 &  8.60              \\
&&    LASSOy     &&&       2.40 &  6.20 & 24.70 & 27.00  &&   1.20 &  3.30 & 12.50 & 11.30   &&   0.60 &  2.00 &  7.40 &  6.90   &&  10.90 & 13.40 & 28.60 & 29.00  &&   7.00 &  9.00 & 14.30 & 12.30   &&  4.50  &  5.80 &  8.30 &  7.30              \\
&&    GLS-LASSO  &&&       2.00 &  2.40 & 17.80 & 21.00  &&   1.00 &  1.20 &  9.90 &  9.80   &&   0.50 &  0.70 &  5.90 &  6.30   &&  10.50 & 11.20 & 23.20 & 24.90  &&   6.70 &  7.50 & 12.50 & 11.90   &&  4.30  &  4.90 &  7.80 &  7.40              \\
&&    ARDL-LAS   &&&       1.10 &  3.10 & 12.30 & 11.00  &&   0.60 &  1.70 &  7.20 &  6.60   &&   0.30 &  1.10 &  4.20 &  4.00   &&   5.20 &  7.50 & 12.50 & 11.20  &&   3.40 &  4.70 &  7.20 &  6.50   &&  2.20  &  3.10 &  4.20 &  3.90              \\
&&    FaSel      &&&       4.10 &  9.20 & 29.40 & 32.10  &&   2.00 &  5.20 & 15.90 & 15.20   &&   6.70 &  4.90 &  8.40 &  8.70   &&   4.00 & 10.90 & 34.80 & 36.40  &&   1.90 &  6.20 & 20.60 & 19.90   && 28.20  & 19.20 & 12.30 & 11.60              \\
&&    ARMAr-LAS  &&&       2.80 &  3.00 &  2.90 &  3.10  &&   1.50 &  1.30 &  1.50 &  1.50   &&   0.70 &  0.70 &  0.70 &  0.80   &&  11.70 & 11.90 & 12.10 & 11.60  &&   7.70 &  7.80 &  7.60 &  7.70   &&  4.90  &  4.90 &  4.90 &  4.80              \\\cline{2-34}

10&&&&&&&&&&&&&&&&&&& &&&&&&&&&&&&&&&\\
&CoEr&&&&&&&&&&&&&&&&&& &&&&&&&&&&&&&&&\\\cline{2-2}
&&    LASSOy       &&&     1.00 &  0.97 &  0.92 &  0.92  &&  1.00 &  0.98 &  0.93 &  0.92   &&   1.00 &  0.98 &  0.95 &  0.94    &&   1.00 &  0.98 &  0.94 &  0.93   &&  1.00 &  0.98 &  0.96 &  0.95  &&  1.00 &  0.99 &  0.97 &  0.96                  \\
&&    GLS-LASSO    &&&     0.94 &  0.78 &  0.82 &  0.85  &&  0.94 &  0.79 &  0.87 &  0.87   &&   0.94 &  0.81 &  0.92 &  0.90    &&   0.94 &  0.79 &  0.85 &  0.86   &&  0.95 &  0.81 &  0.87 &  0.85  &&  0.95 &  0.85 &  0.90 &  0.88                  \\
&&    ARDL-LAS     &&&     1.01 &  0.96 &  0.78 &  0.80  &&  1.01 &  0.97 &  0.87 &  0.88   &&   1.01 &  0.97 &  0.90 &  0.90    &&   0.99 &  0.95 &  0.80 &  0.81   &&  1.00 &  0.96 &  0.88 &  0.89  &&  0.99 &  0.97 &  0.91 &  0.92                  \\
&&    FaSel        &&&     1.77 &  1.42 &  1.12 &  1.08  &&  1.42 &  1.20 &  1.04 &  1.00   &&   1.36 &  1.13 &  1.00 &  0.96    &&   1.19 &  1.13 &  1.07 &  1.07   &&  1.16 &  1.09 &  1.07 &  1.04  &&  1.52 &  1.22 &  1.04 &  1.02                  \\
&&    ARMAr-LAS    &&&     0.94 &  0.77 &  0.51 &  0.50  &&  0.95 &  0.77 &  0.59 &  0.57   &&   0.95 &  0.79 &  0.63 &  0.60    &&   0.94 &  0.75 &  0.52 &  0.50   &&  0.95 &  0.75 &  0.59 &  0.56  &&  0.96 &  0.78 &  0.63 &  0.60                  \\
&RMSFE&&&&&&&&&&&&&&&&&& &&&&&&&&&&&&&&&\\\cline{2-2}
&&    LASSOy      &&&      1.00 &  0.98 &  0.96 &  0.96  &&  1.00 &  0.98 &  0.96 &  0.96   &&   1.00 &  0.99 &  0.98 &  0.94    &&   1.00 &  0.99 &  0.94 &  0.96   &&  1.00 &  0.99 &  0.97 &  0.95  &&  1.00 &  0.99 &  0.97 &  0.97                \\
&&    GLS-LASSO   &&&      0.95 &  0.83 &  0.87 &  0.89  &&  0.95 &  0.83 &  0.87 &  0.82   &&   0.96 &  0.85 &  0.89 &  0.86    &&   0.94 &  0.86 &  0.87 &  0.87   &&  0.96 &  0.88 &  0.84 &  0.81  &&  0.96 &  0.88 &  0.88 &  0.85                \\
&&    ARDL-LAS    &&&      1.01 &  0.99 &  0.99 &  1.05  &&  1.01 &  0.99 &  0.98 &  0.96   &&   1.01 &  1.01 &  0.97 &  0.97    &&   1.00 &  0.98 &  0.91 &  0.96   &&  1.01 &  0.99 &  0.93 &  0.89  &&  1.00 &  0.99 &  0.94 &  0.93                \\
&&    FaSel       &&&      0.98 &  0.98 &  0.97 &  0.96  &&  1.00 &  0.97 &  0.96 &  0.93   &&   0.94 &  0.97 &  0.95 &  0.91    &&   1.04 &  1.03 &  0.99 &  0.98   &&  1.06 &  1.04 &  0.93 &  0.88  &&  1.08 &  1.05 &  0.92 &  0.90                \\
&&    ARMAr-LAS   &&&      0.97 &  0.83 &  0.70 &  0.69  &&  0.99 &  0.83 &  0.70 &  0.65   &&   0.98 &  0.84 &  0.74 &  0.69    &&   0.96 &  0.85 &  0.68 &  0.64   &&  0.97 &  0.86 &  0.66 &  0.65  &&  0.96 &  0.86 &  0.70 &  0.69                \\
& TP&&&&&&&&&&&&&&&&&& &&&&&&&&&&&&&&&\\\cline{2-2}
&&    LASSO      &&&      97.20 & 93.10 & 88.00 & 87.70  && 97.10 & 93.60 & 87.30 & 85.30   &&  96.70 & 93.00 & 88.30 & 84.90    &&  98.20 & 94.70 & 88.40 & 87.20   && 97.40 & 93.00 & 84.40 & 78.70  && 96.60 & 91.80 & 82.80 & 77.50              \\
&&    LASSOy     &&&      97.20 & 93.20 & 87.90 & 87.60  && 97.10 & 93.40 & 87.00 & 84.70   &&  96.70 & 93.20 & 87.90 & 84.60    &&  98.20 & 94.80 & 87.70 & 86.50   && 97.40 & 93.10 & 84.10 & 77.20  && 96.60 & 91.70 & 82.20 & 76.50              \\
&&    GLS-LASSO  &&&      98.00 & 97.30 & 89.80 & 88.10  && 97.80 & 97.50 & 89.20 & 87.90   &&  97.60 & 97.30 & 89.20 & 87.30    &&  98.70 & 98.10 & 90.10 & 88.50   && 98.00 & 97.00 & 88.80 & 85.60  && 97.40 & 96.00 & 86.60 & 83.70              \\
&&    ARDL-LAS   &&&      96.80 & 92.70 & 85.70 & 83.10  && 96.70 & 93.10 & 86.70 & 83.80   &&  96.50 & 92.70 & 87.40 & 83.50    &&  98.10 & 94.60 & 86.30 & 82.30   && 97.20 & 93.00 & 84.50 & 78.80  && 96.50 & 91.40 & 82.60 & 77.80              \\
&&    FaSel      &&&      47.50 & 61.00 & 76.90 & 80.80  && 70.60 & 80.30 & 85.60 & 85.60   &&  91.00 & 90.00 & 87.80 & 87.20    &&  91.90 & 86.80 & 84.80 & 84.30   && 89.30 & 85.00 & 84.60 & 82.70  && 90.30 & 84.90 & 84.40 & 81.60              \\
&&    ARMAr-LAS  &&&      98.20 & 97.80 & 97.20 & 96.60  && 97.70 & 97.90 & 96.90 & 96.30   &&  97.50 & 97.60 & 96.60 & 95.90    &&  98.80 & 98.70 & 97.60 & 97.10   && 98.00 & 97.90 & 96.90 & 96.40  && 97.30 & 97.20 & 95.80 & 95.00              \\
& FP&&&&&&&&&&&&&&&&&& &&&&&&&&&&&&&&&\\\cline{2-2}
&&    LASSO      &&&      2.60 &  7.10 & 30.20 & 31.60  &&  1.30 &  4.30 & 13.70 & 12.30   &&   0.70 &  2.60 &  8.00 &  7.30    &&  11.20 & 15.00 & 33.50 & 34.30   &&  7.60 & 10.00 & 15.20 & 12.90  &&  4.70 &  6.20 &  8.80 &  7.70              \\
&&    LASSOy     &&&      2.70 &  7.00 & 28.00 & 29.50  &&  1.30 &  4.00 & 12.20 & 10.70   &&   0.70 &  2.40 &  7.50 &  6.60    &&  11.10 & 14.40 & 31.10 & 31.60   &&  7.50 &  9.60 & 14.10 & 11.70  &&  4.70 &  6.10 &  8.20 &  7.10              \\
&&    GLS-LASSO  &&&      2.10 &  2.70 & 18.60 & 21.30  &&  1.00 &  1.40 &  8.90 &  9.30   &&   0.50 &  0.80 &  5.80 &  5.90    &&  11.00 & 11.40 & 24.10 & 26.00   &&  7.30 &  7.80 & 12.50 & 11.90  &&  4.50 &  5.10 &  7.70 &  7.30              \\
&&    ARDL-LAS   &&&      1.30 &  3.20 & 11.40 & 10.00  &&  0.60 &  2.10 &  7.10 &  6.20   &&   0.40 &  1.30 &  4.30 &  3.70    &&   5.30 &  7.30 & 11.90 & 10.40   &&  3.70 &  4.90 &  7.20 &  6.30  &&  2.30 &  3.10 &  4.30 &  3.80              \\
&&    FaSel      &&&      7.10 & 12.30 & 30.70 & 30.70  &&  2.80 &  6.20 & 15.10 & 14.10   &&   5.30 &  4.70 &  8.00 &  8.20    &&   6.50 & 14.70 & 36.10 & 38.00   &&  3.00 &  8.90 & 21.40 & 20.40  && 13.50 & 11.40 & 12.80 & 11.90              \\
&&    ARMAr-LAS  &&&      2.90 &  3.20 &  3.30 &  3.10  &&  1.40 &  1.50 &  1.50 &  1.60   &&   0.70 &  0.80 &  0.80 &  0.80    &&  12.10 & 12.00 & 12.50 & 12.30   &&  8.10 &  7.90 &  7.80 &  7.70  &&  5.00 &  5.00 &  5.00 &  4.90              \\\cline{2-34}

\hline\hline
\end{tabular}
}
}
\end{table}

\subsection{Analysis of the minimum eigenvalues}\label{Sec:MinEig}
\begin{figure}[t]
\graphicspath{{images/}}
\centering
  \captionsetup[subfigure]{oneside,margin={0.5cm,0cm}}
\subfloat[\scriptsize DGP (A)]{\includegraphics[width=6cm]{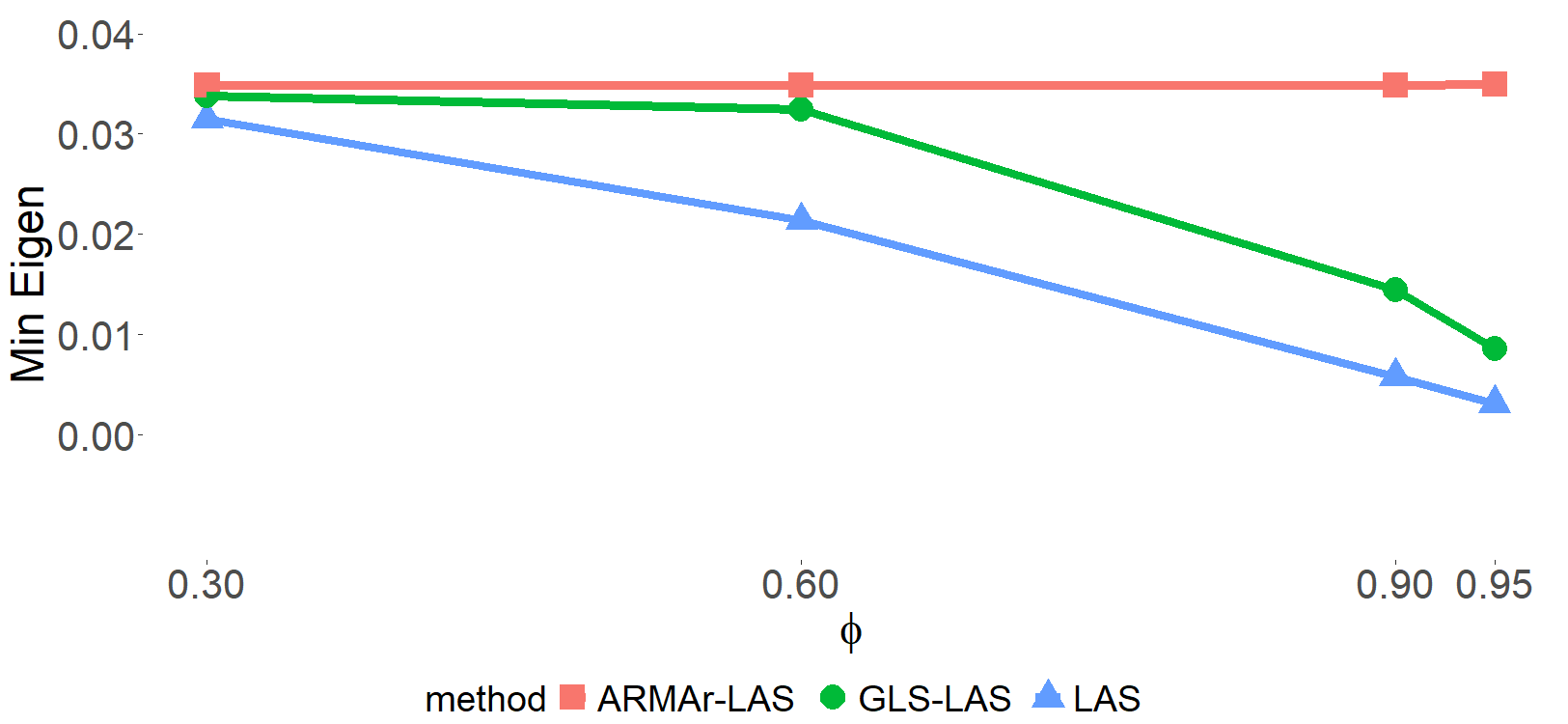}}\hfil
\subfloat[\scriptsize DGP (B)]{\includegraphics[width=6cm]{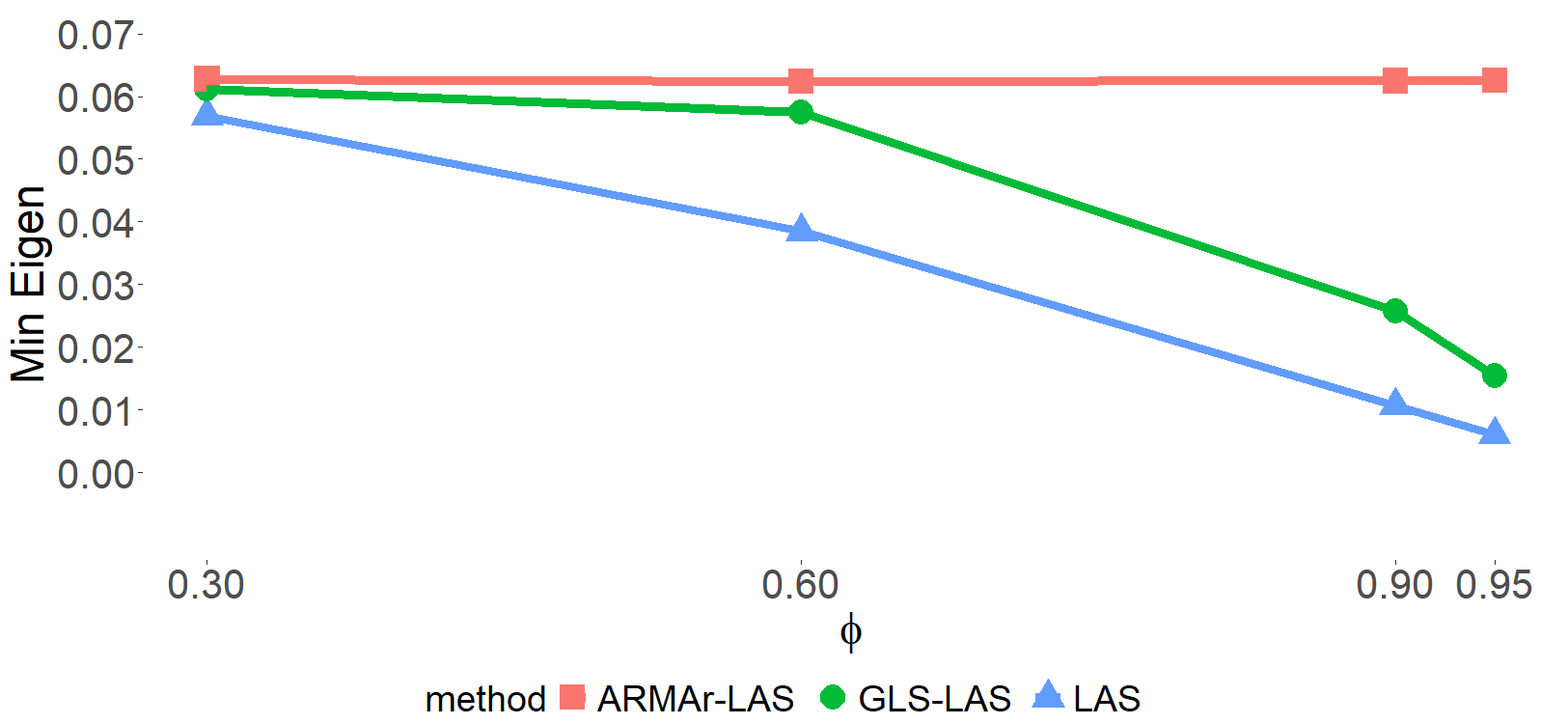}}
\caption{\footnotesize Minimum eigenvalues for the design matrix of LAS, GLS-LAS, and ARMAr-LAS, for various degrees of serial correlation ($\phi$) under DGPs (A) and (B).}
\label{fig:MinEigAB}
\end{figure}

In this section, we compare the minimum eigenvalues of the design matrix of LAS, GLS-LAS, and ARMAr-LAS in the case of $n=50$ and $SNR=10$. Figure~\ref{fig:MinEigAB} shows the average of the minimum eigenvalues obtained in the experiments presented in Section~\ref{ARMArLASSOMonteCarlo_supp}. Both LAS and GLS-LAS reduce their minimum eigenvalues as $\phi$ increases. This does not happen for ARMAr-LAS, which maintains the same value regardless of the degree of serial correlation. Figure~\ref{fig:MinEigCD} shows the same results but for the experiments presented in Supplement~\ref{ARMArLASSOMonteCarlo_supp}. In this case, we compare the minimum eigenvalues for the two DGPs reported as 0 for no common factor ($q=0$) and 1 for common factor ($q=1$). ARMAr-LASSO also exhibits larger average minimum eigenvalues in these experiments. This numerical comparison does not establish a corresponding ordering of compatibility constants.

\begin{figure}[t]
\begin{center}
\graphicspath{{images/}}
\includegraphics[scale=0.2]{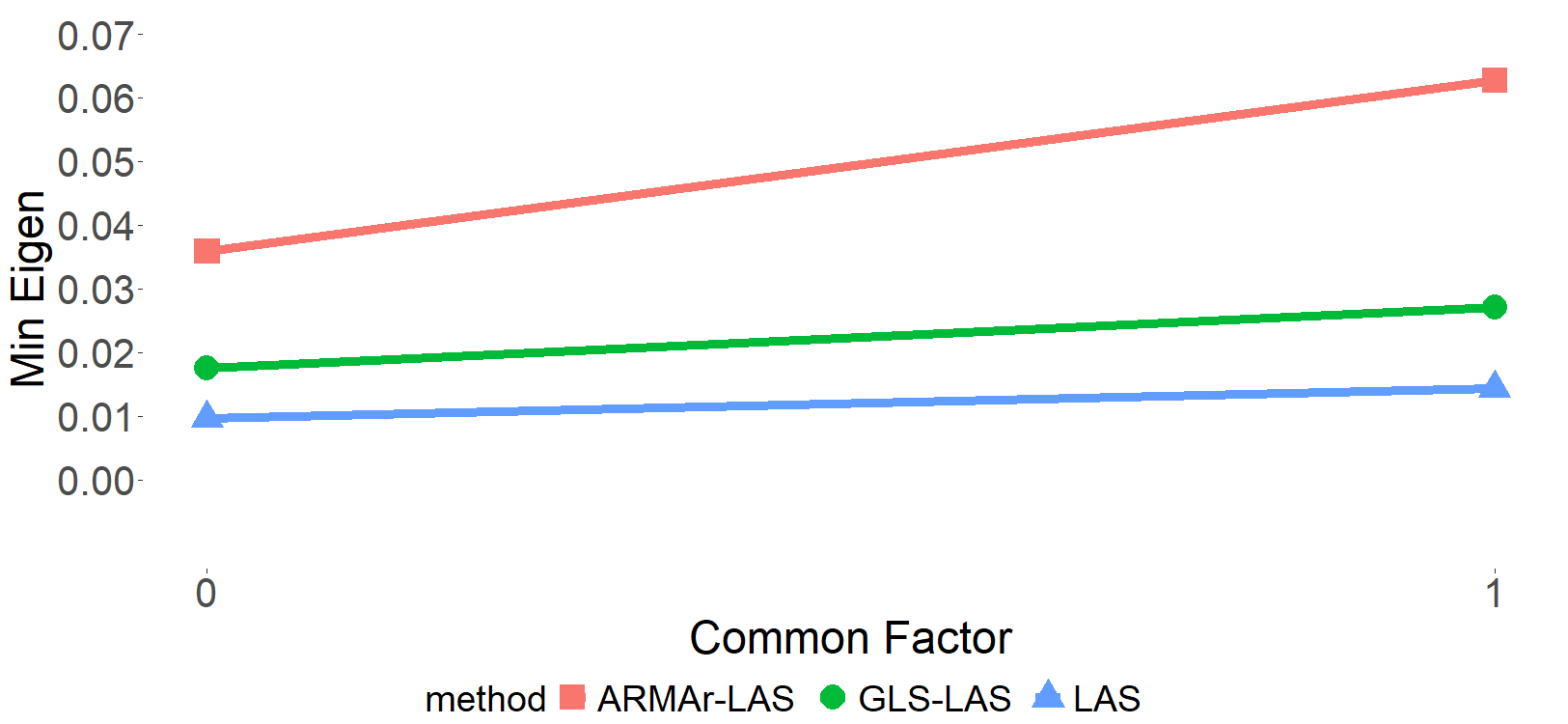}
\end{center}
\caption{\footnotesize Minimum eigenvalues for the design matrix of LAS, GLS-LAS, and ARMAr-LAS, for DGP in Supplement~\ref{ARMArLASSOMonteCarlo_supp}. In this case, we compare the minimum eigenvalues for the two DGPs reported as 0 for no common factor ($q=0$) and 1 for common factor ($q=1$).}
\label{fig:MinEigCD}
\end{figure}

\subsection{Performance in a Large $T$ Regime}\label{sec:largeT}
Here we compare our ARMAr-LAS with the employed LASSO-based benchmarks in the case of DGP (A) (see Section~\ref{ARMArLASSOMonteCarlo_supp}) with $T=1500$, $n=50$, and SNR=10. This section aims to evaluate the performances of ARMAr-LAS in a large sample size regime. Results in Table~\ref{TabDGPAT1500} show that ARMAr-LAS performs as GLS-LAS. This result is expected since under DGP (A) these two estimators coincide. Further, both outperform the other LASSO-based methods providing more accurate coefficient estimates and forecasts, as well as a perfect variable selection accuracy.   
\begin{table}[t]
\centering
\caption{\footnotesize DGPs (A). CoEr, RMSFE (relative to LAS), \%TP and \%FP for LASSO-based benchmarks and ARMAr-LASSO, under 4 values of $\phi$ with $T=1500$ and $n=50$. 
}\label{TabDGPAT1500}
\scalebox{0.75}{
\begin{tabular}
{@{}cccrrrr@{}}
\hline\hline
&& &0.3 &0.6 &0.9 &0.95\\\hline
CoEr&&&&&& \\
&    LASSOy     &   &  1.00 &  1.01 &  0.95 &  0.93     \\
&    GLS-LAS  &      &  0.91 &  0.68 &  0.22 &  0.16     \\
&    ARDL-LAS  &     &  1.01 &  0.83 &  0.28 &  0.24     \\
&    FaSel      &    & 13.13 &  9.86 &  3.04 &  2.02     \\
&    ARMAr-LAS  &    &  0.91 &  0.68 &  0.21 &  0.15     \\\hline
RMSFE&&&&&& \\
&    LASSOy     &     &  1.00 &  0.99 &  0.97 &  0.96    \\
&    GLS-LAS  &       &  0.95 &  0.81 &  0.46 &  0.35    \\
&    ARDL-LAS  &      &  1.00 &  0.83 &  0.49 &  0.39    \\
&    FaSel      &     &  1.01 &  0.99 &  0.99 &  1.00    \\
&    ARMAr-LAS  &     &  0.95 &  0.82 &  0.46 &  0.34    \\\hline
 TP&&&&&& \\
&    LASSO      &     &100.00 &100.00 &100.00 &100.00   \\
&    LASSOy     &     &100.00 &100.00 &100.00 &100.00   \\
&    GLS-LAS  &       &100.00 &100.00 &100.00 &100.00    \\
&    ARDL-LAS  &      &100.00 &100.00 &100.00 &100.00    \\
&    FaSel      &     & 65.80 & 67.90 & 80.90 & 86.90     \\
&    ARMAr-LAS  &     &100.00 &100.00 &100.00 &100.00     \\\hline
 FP&&&&&& \\
&    LASSO      &    &  0.10 &  0.10 &  1.70 &  2.20    \\
&    LASSOy     &    &  0.10 &  0.10 &  1.60 &  2.20     \\
&    GLS-LAS  &      &  0.00 &  0.00 &  0.00 &  0.10     \\
&    ARDL-LAS  &     &  0.00 &  1.40 &  1.40 &  1.40      \\
&    FaSel      &    &  1.10 &  1.10 &  1.80 &  2.10      \\
&    ARMAr-LAS  &    &  0.10 &  0.10 &  0.00 &  0.00      \\\cline{2-7}  \hline\hline
\end{tabular}
}
\end{table}

\subsection{Performance with Misspecified Autoregressive Structure}\label{sec:misspecifiedAR}
In this section, we compare our ARMAr-LASSO with the LASSO-based benchmarks in the case where the former misspecifies the autoregressive model of predictors. In particular, we generated both predictors and error terms from an AR(2) model with autoregressive coefficients equal to 1.2 and -0.4, but ARMAr-LASSO filters predictors through an AR(1) model. We consider $T=150$, SNR=10, and $n=50,150,300$. Results are reported in Table~\ref{TabMisspecifiedAR}. Also in this case where the predictors are filtered with a misspecified autoregressive model, ARMAr-LAS outperforms LASSO-based benchmarks. This is because, despite the misspecification, the filter can remove the majority of serial correlation with the proper estimation of a single autoregressive coefficient. This is corroborated by the averages of the minimum eigenvalues of the correlation matrices for LAS, GLS-LAS, and ARMAr-LAS which are 0.00680, 0.01824, and 0.02861, respectively.

\begin{table}[t]
\centering
\caption{\footnotesize  CoEr, RMSFE (relative to LAS), \%TP and \%FP for LASSO-based benchmarks and ARMAr-LASSO, under 3 values of $n$.
}\label{TabMisspecifiedAR}
\scalebox{0.75}{
\begin{tabular}
{@{}cccrrr@{}}
\hline
&&50&150&300\\\hline
CoEr&&&&\\
&    LASSOy        &0.97   &0.98   &0.99     \\
&    GLS-LAS       &0.59   &0.72   &0.80   \\
&    ARDL-LAS      &0.76   &0.94   &0.94    \\
&    FaSel         &1.67   &1.33   &1.28     \\
&    ARMAr-LAS     &0.43   &0.60   &0.64    \\\hline
RMSFE&&&&\\
&    LASSOy        &0.99   &0.99   &0.98     \\
&    GLS-LAS       &0.70   &0.70   &0.76   \\
&    ARDL-LAS      &0.94   &0.98   &0.98    \\
&    FaSel         &1.00   &0.91   &0.80     \\
&    ARMAr-LAS     &0.59   &0.61   &0.66        \\\hline
 TP&&&&\\
&    LASSO        & 99.70 & 99.70 & 99.80        \\
&    LASSOy       & 99.70 & 99.70 & 99.80        \\
&    GLS-LAS      & 99.90 &100.00 &100.00      \\
&    ARDL-LAS     & 99.80 & 99.80 & 99.80       \\
&    FaSel        & 78.90 & 95.80 & 98.30        \\
&    ARMAr-LAS    &100.00 &100.00 &100.00        \\\hline
 FP&&&&\\
&    LASSO           &51.20  &12.90  & 8.40     \\
&    LASSOy          &50.10  &12.40  & 8.10     \\
&    GLS-LAS         &17.70  & 4.80  & 4.30   \\
&    ARDL-LAS        &35.50  & 7.50  & 4.80    \\
&    FaSel           &54.00  &19.30  &12.80     \\
&    ARMAr-LAS       & 7.00  & 2.50  & 1.40     \\
\hline
\end{tabular}
}
\end{table}

\section{List of Time Series 
in the Euro Area Data}\label{Stat_EA}
We report the list of series for the Euro Area dataset adopted in the forecasting exercise (obtained from~\cite{Ale2021}). As for the FRED data, the column tcode denotes the data transformation for a given series $x_t$: (1) no transformation; (2) $\Delta x_t$; (3)$\Delta^2 x_t$; (4) $ \text{log}(x_t)$; (5) $\Delta \text{log}(x_t)$; (6) $\Delta^2 \text{log}(x_t)$. (7) $\Delta (x_t/x_{t-} - 1.0)$.

\noindent The acronyms for the sectors refer to:
\setcounter{bean}{0}
\begin{list}
{(\alph{bean})}{\usecounter{bean}}  \item ICS: Industry \& Construction Survey
  \item CCI: Consumer Confidence Indicators
  \item M\&IR: Money \& Interest Rates
  \item IP: Industrial Production
  \item HCPI: Harm. Consumer Price Index
  \item PPI: Producer Price Index
  \item T\/O: Turnover \& Retail Sale
  \item HUR: Harm. Unemployment rate
  \item SI: Service Svy.
\end{list}

As mentioned in the main text, for the first variable of each group we report in brackets its autocorrelation function to show that predictors are serially correlated.

\scriptsize
\begin{longtable}{cllcc}
\caption{\footnotesize Euro Area macroeconomic variables from \cite{Ale2021}} \label{tab:long} \\

\hline \multicolumn{1}{c}{\textbf{ID}} &  \multicolumn{1}{c}{\textbf{Description}} & \multicolumn{1}{c}{\textbf{Area}} & \multicolumn{1}{c}{\textbf{Sector}} & \multicolumn{1}{c}{\textbf{Tcode}} \\ \hline
\endfirsthead

\multicolumn{5}{c}%
{{\bfseries \tablename\ \thetable{} -- continued from previous page}} \\
\hline \multicolumn{1}{c}{\textbf{ID}} &  \multicolumn{1}{c}{\textbf{Description}} & \multicolumn{1}{c}{\textbf{Area}} & \multicolumn{1}{c}{\textbf{Sector}} & \multicolumn{1}{c}{\textbf{Tcode}} \\ \hline
\endhead

\hline \multicolumn{5}{r}{{Continued on next page}} \\
\endfoot

\hline \hline
\endlastfoot

    1     & Ind Svy: Employment Expectations ($\pmb{acf:\  0.97}$)  & EA    & ICS   & 1 \\
    2     & Ind Svy: Export Order-Book Levels  & EA    & ICS   & 1 \\
    3     & Ind Svy: Order-Book Levels  & EA    & ICS   & 1 \\
    4     & Ind Svy: Mfg - Selling Price Expectations  & EA    & ICS   & 1 \\
    5     & Ind Svy: Production Expectations  & EA    & ICS   & 1 \\
    6     & Ind Svy: Production Trend  & EA    & ICS   & 1 \\
    7     & Ind Svy: Mfg - Stocks Of Finished Products  & EA    & ICS   & 1 \\
    8     & Constr. Svy: Price Expectations & EA    & ICS   & 1 \\
    9     & Ind Svy: Export Order Book Position  & EA    & ICS   & 1 \\
    10    & Ind Svy: Production Trends In Recent Mth. & EA    & ICS   & 1 \\
    11    & Ind Svy: Selling Prc. Expect. Mth. Ahead  & EA    & ICS   & 1 \\
    12    & Ret. Svy: Employment & EA    & ICS   & 1 \\
    13    & Ret. Svy: Orders Placed With Suppliers & EA    & ICS   & 1 \\
    14    & Constr. Svy: Synthetic Bus. Indicator & FR    & ICS   & 1 \\
    15    & Bus. Svy: Constr. Sector - Capacity Utilisation Rate & FR    & ICS   & 1 \\
    16    & Constr. Svy: Activity Expectations & FR    & ICS   & 1 \\
    17    & Constr. Svy: Price Expectations & FR    & ICS   & 1 \\
    18    & Constr. Svy: Unable To Increase Capacity & FR    & ICS   & 1 \\
    19    & Constr. Svy: Workforce Changes & FR    & ICS   & 1 \\
    20    & Constr. Svy: Workforce Forecast Changes & FR    & ICS   & 1 \\
    21    & Svy: Mfg Output - Order Book \& Demand & FR    & ICS   & 1 \\
    22    & Svy: Mfg Output - Order Book \& Foreign Demand & FR    & ICS   & 1 \\
    23    & Svy: Mfg Output - Personal Outlook & FR    & ICS   & 1 \\
    24    & Svy: Auto Ind - Order Book \& Demand & FR    & ICS   & 1 \\
    25    & Svy: Auto Ind - Personal Outlook & FR    & ICS   & 1 \\
    26    & Svy: Basic \& Fab Met Pdt Ex Mach \& Eq - Personal Outlook & FR    & ICS   & 1 \\
    27    & Svy: Ele \& Elec Eq, Mach Eq - Order Book \& Demand & FR    & ICS   & 1 \\
    28    & Svy: Ele \& Elec Eq, Mach Eq - Order Book \& Foreign Demand & FR    & ICS   & 1 \\
    29    & Svy: Ele \& Elec Eq, Mach Eq - Personal Outlook & FR    & ICS   & 1 \\
    30    & Svy: Mfg Output - Price Outlook & FR    & ICS   & 1 \\
    31    & Svy: Mfg Of Chemicals \& Chemical Pdt - Order Book \& Demand & FR    & ICS   & 1 \\
    32    & Svy: Mfg Of Chemicals \& Chemical Pdt - Personal Outlook & FR    & ICS   & 1 \\
    33    & Svy: Mfg Of Food Pr \& Beverages - Order Book \& Demand & FR    & ICS   & 1 \\
    34    & Svy: Mfg Of Food Pr \& Beverages - Order Book \& Foreign Demand & FR    & ICS   & 1 \\
    35    & Svy: Mfg Of Trsp Eq - Finished Goods Inventories & FR    & ICS   & 1 \\
    36    & Svy: Mfg Of Trsp Eq - Order Book \& Demand & FR    & ICS   & 1 \\
    37    & Svy: Mfg Of Trsp Eq - Order Book \& Foreign Demand & FR    & ICS   & 1 \\
    38    & Svy: Mfg Of Trsp Eq - Personal Outlook & FR    & ICS   & 1 \\
    39    & Svy: Oth Mfg, Mach \& Eq Rpr \& Instal - Ord Book \& Demand & FR    & ICS   & 1 \\
    40    & Svy: Oth Mfg, Mach \& Eq Rpr \& Instal - Ord Book \& Fgn Demand & FR    & ICS   & 1 \\
    41    & Svy: Oth Mfg, Mach \& Eq Rpr \& Instal - Personal Outlook & FR    & ICS   & 1 \\
    42    & Svy: Other Mfg - Order Book \& Demand & FR    & ICS   & 1 \\
    43    & Svy: Rubber, Plastic \& Non Met Pdt - Order Book \& Demand & FR    & ICS   & 1 \\
    44    & Svy: Rubber, Plastic \& Non Met Pdt - Order Book \& Fgn Demand & FR    & ICS   & 1 \\
    45    & Svy: Rubber, Plastic \& Non Met Pdt - Personal Outlook & FR    & ICS   & 1 \\
    46    & Svy: Total Ind - Order Book \& Demand & FR    & ICS   & 1 \\
    47    & Svy: Total Ind - Order Book \& Foreign Demand & FR    & ICS   & 1 \\
    48    & Svy: Total Ind - Personal Outlook & FR    & ICS   & 1 \\
    49    & Svy: Total Ind - Price Outlook & FR    & ICS   & 1 \\
    50    & Svy: Wood \& Paper, Print \& Media - Ord Book \& Fgn Demand & FR    & ICS   & 1 \\
    51    & Trd. \& Ind: Bus Sit & DE    & ICS   & 1 \\
    52    & Trd. \& Ind: Bus Expect In 6Mo & DE    & ICS   & 1 \\
    53    & Trd. \& Ind: Bus Sit & DE    & ICS   & 1 \\
    54    & Trd. \& Ind: Bus Climate & DE    & ICS   & 1 \\
    55    & Cnstr Ind: Bus Climate & DE    & ICS   & 1 \\
    56    & Mfg: Bus Climate & DE    & ICS   & 1 \\
    57    & Mfg: Bus Climate & DE    & ICS   & 1 \\
    58    & Mfg Cons Gds: Bus Climate & DE    & ICS   & 1 \\
    59    & Mfg (Excl Fbt): Bus Climate & DE    & ICS   & 1 \\
    60    & Whsle (Incl Mv): Bus Climate & DE    & ICS   & 1 \\
    61    & Mfg: Bus Sit & DE    & ICS   & 1 \\
    62    & Mfg: Bus Sit & DE    & ICS   & 1 \\
    63    & Mfg (Excl Fbt): Bus Sit & DE    & ICS   & 1 \\
    64    & Mfg (Excl Fbt): Bus Sit & DE    & ICS   & 1 \\
    65    & Cnstr Ind: Bus Expect In 6Mo & DE    & ICS   & 1 \\
    66    & Cnstr Ind: Bus Expect In 6Mo & DE    & ICS   & 1 \\
    67    & Mfg: Bus Expect In 6Mo & DE    & ICS   & 1 \\
    68    & Mfg: Bus Expect In 6Mo & DE    & ICS   & 1 \\
    69    & Mfg Cons Gds: Bus Expect In 6Mo & DE    & ICS   & 1 \\
    70    & Mfg (Excl Fbt): Bus Expect In 6Mo & DE    & ICS   & 1 \\
    71    & Mfg (Excl Fbt): Bus Expect In 6Mo & DE    & ICS   & 1 \\
    72    & Rt (Incl Mv): Bus Expect In 6Mo & DE    & ICS   & 1 \\
    73    & Whsle (Incl Mv): Bus Expect In 6Mo & DE    & ICS   & 1 \\
    74    & Bus. Conf. Indicator & IT    & ICS   & 1 \\
    75    & Order Book Level: Ind & ES    & ICS   & 1 \\
    76    & Order Book Level: Foreign - Ind & ES    & ICS   & 1 \\
    77    & Order Book Level: Investment Goods & ES    & ICS   & 1 \\
    78    & Order Book Level: Int. Goods & ES    & ICS   & 1 \\
    79    & Production Level - Ind & ES    & ICS   & 1 \\
    80    & Cons. Confidence Indicator ($\pmb{acf:\  0.98}$)  & EA    & CCI   & 1 \\
    81    & Cons. Svy: Economic Situation Last 12 Mth. - Emu 11/12 & EA    & CCI   & 1 \\
    82    & Cons. Svy: Possible Savings Opinion & FR    & CCI   & 1 \\
    83    & Cons. Svy: Future Financial Situation & FR    & CCI   & 1 \\
    84    & Svy - Households, Economic Situation Next 12M & FR    & CCI   & 1 \\
    85    & Cons. Confidence Indicator - DE & DE    & CCI   & 1 \\
    86    & Cons. Confidence Index & DE    & CCI   & 5 \\
    87    & Gfk Cons. Climate Svy - Bus. Cycle Expectations & DE    & CCI   & 1 \\
    88    & Cons.S Confidence Index & DE    & CCI   & 5 \\
    89    & Cons. Confidence Climate (Balance) & DE    & CCI   & 1 \\
    90    & Cons. Svy: Economic Climate Index (N.West It) & IT    & CCI   & 5 \\
    91    & Cons. Svy: Economic Climate Index (Southern It) & IT    & CCI   & 5 \\
    92    & Cons. Svy: General Economic Situation (Balance) & IT    & CCI   & 1 \\
    93    & Cons. Svy: Prices In Next 12 Mths. - Lower & IT    & CCI   & 5 \\
    94    & Cons. Svy: Unemployment Expectations (Balance) & IT    & CCI   & 1 \\
    95    & Cons. Svy: Unemployment Expectations - Approx. Same & IT    & CCI   & 5 \\
    96    & Cons. Svy: Unemployment Expectations - Large Increase & IT    & CCI   & 5 \\
    97    & Cons. Svy: Unemployment Expectations - Small Increase & IT    & CCI   & 5 \\
    98    & Cons. Svy: General Economic Situation (Balance) & IT    & CCI   & 1 \\
    99    & Cons. Svy: Household Budget - Deposits To/Withdrawals & ES    & CCI   & 5 \\
    100   & Cons. Svy: Household Economy (Cpy) - Much Worse & FR    & CCI   & 5 \\
    101   & Cons. Svy: Italian Econ.In Next 12 Mths.- Much Worse & FR    & CCI   & 5 \\
    102   & Cons. Svy: Major Purchase Intentions - Balance & FR    & CCI   & 1 \\
    103   & Cons. Svy: Major Purchase Intentions - Much Less & FR    & CCI   & 5 \\
    104   & Cons. Svy: Households Fin Situation - Balance & FR    & CCI   & 1 \\
    105   & Indl. Prod. - Excluding Constr. ($\pmb{acf:\  -0.21}$)   & EA    & IP    & 5 \\
    106   & Indl. Prod. - Cap. Goods  & EA    & IP    & 5 \\
    107   & Indl. Prod. - Cons. Non-Durables  & EA    & IP    & 5 \\
    108   & Indl. Prod. - Cons. Durables  & EA    & IP    & 5 \\
    109   & Indl. Prod. - Cons. Goods  & EA    & IP    & 5 \\
    110   & Indl. Prod. & FR    & IP    & 5 \\
    111   & Indl. Prod. - Mfg & FR    & IP    & 5 \\
    112   & Indl. Prod. - Mfg (2010=100) & FR    & IP    & 5 \\
    113   & Indl. Prod. - Manuf. Of Motor Vehicles, Trailers, Semitrailers & FR    & IP    & 5 \\
    114   & Indl. Prod. - Int. Goods & FR    & IP    & 5 \\
    115   & Indl. Prod. - Indl. Prod. - Constr. & FR    & IP    & 5 \\
    116   & Indl. Prod. - Manuf. Of Wood And Paper Products & FR    & IP    & 5 \\
    117   & Indl. Prod. - Manuf. Of Computer, Electronic And Optical Prod & FR    & IP    & 5 \\
    118   & Indl. Prod. - Manuf. Of Electrical Equipment & FR    & IP    & 5 \\
    119   & Indl. Prod. - Manuf. Of Machinery And Equipment & FR    & IP    & 5 \\
    120   & Indl. Prod. - Manuf. Of Transport Equipment & FR    & IP    & 5 \\
    121   & Indl. Prod. - Other Mfg & FR    & IP    & 5 \\
    122   & Indl. Prod. - Manuf. Of Chemicals And Chemical Products & FR    & IP    & 5 \\
    123   & Indl. Prod. - Manuf. Of Rubber And Plastics Products & FR    & IP    & 5 \\
    124   & Indl. Prod. - Investment Goods & IT    & IP    & 5 \\
    125   & Indl. Prod. & IT    & IP    & 5 \\
    126   & Indl. Prod.  & IT    & IP    & 5 \\
    127   & Indl. Prod. - Cons. Goods - Durable & IT    & IP    & 5 \\
    128   & Indl. Prod. - Investment Goods & IT    & IP    & 5 \\
    129   & Indl. Prod. - Int. Goods & IT    & IP    & 5 \\
    130   & Indl. Prod. - Chemical Products \& Synthetic Fibres & IT    & IP    & 5 \\
    131   & Indl. Prod. - Machines \& Mechanical Apparatus & IT    & IP    & 5 \\
    132   & Indl. Prod. - Means Of Transport & IT    & IP    & 5 \\
    133   & Indl. Prod. - Metal \& Metal Products & IT    & IP    & 5 \\
    134   & Indl. Prod. - Rubber Items \& Plastic Materials & IT    & IP    & 5 \\
    135   & Indl. Prod. - Wood \& Wood Products & IT    & IP    & 5 \\
    136   & Indl. Prod. & IT    & IP    & 5 \\
    137   & Indl. Prod. - Computer, Electronic And Optical Products & IT    & IP    & 5 \\
    138   & Indl. Prod. - Basic Pharmaceutical Products & IT    & IP    & 5 \\
    139   & Indl. Prod. - Constr. & DE    & IP    & 5 \\
    140   & Indl. Prod. - Ind Incl Cnstr & DE    & IP    & 5 \\
    141   & Indl. Prod. - Mfg & DE    & IP    & 5 \\
    142   & Indl. Prod. - Rebased To 1975=100 & DE    & IP    & 5 \\
    143   & Indl. Prod. - Chems \& Chem Prds & DE    & IP    & 5 \\
    144   & Indl. Prod. - Ind Excl Cnstr & DE    & IP    & 5 \\
    145   & Indl. Prod. - Ind Excl Energy \& Cnstr & DE    & IP    & 5 \\
    146   & Indl. Prod. - Mining \& Quar & DE    & IP    & 5 \\
    147   & Indl. Prod. - Cmptr, Eleccl \& Opt Prds, Elecl Eqp & DE    & IP    & 5 \\
    148   & Indl. Prod. - Interm Goods & DE    & IP    & 5 \\
    149   & Indl. Prod. - Cap. Goods & DE    & IP    & 5 \\
    150   & Indl. Prod. - Durable Cons Goods & DE    & IP    & 5 \\
    151   & Indl. Prod. - Tex \& Wearing Apparel & DE    & IP    & 5 \\
    152   & Indl. Prod. - Pulp, Paper\&Prds, Pubshg\&Print & DE    & IP    & 5 \\
    153   & Indl. Prod. - Chem Prds & DE    & IP    & 5 \\
    154   & Indl. Prod. - Rub\&Plast Prds & DE    & IP    & 5 \\
    155   & Indl. Prod. - Basic Mtls & DE    & IP    & 5 \\
    156   & Indl. Prod. - Cmptr, Eleccl \& Opt Prds, Elecl Eqp & DE    & IP    & 5 \\
    157   & Indl. Prod. - Motor Vehicles, Trailers\&Semi Trail & DE    & IP    & 5 \\
    158   & Indl. Prod. - Tex \& Wearing Apparel & DE    & IP    & 5 \\
    159   & Indl. Prod. - Paper \& Prds, Print, Reprod Of Recrd Media & DE    & IP    & 5 \\
    160   & Indl. Prod. - Chems \& Chem Prds & DE    & IP    & 5 \\
    161   & Indl. Prod. - Basic Mtls, Fab Mtl Prds, Excl Mach\&Eqp & DE    & IP    & 5 \\
    162   & Indl. Prod. - Repair \& Install Of Mach \& Eqp & DE    & IP    & 5 \\
    163   & Indl. Prod. - Mfg Excl Cnstr \& Fbt & DE    & IP    & 5 \\
    164   & Indl. Prod. - Mining \& Ind Excl Fbt & DE    & IP    & 5 \\
    165   & Indl. Prod. - Ind Excl Fbt & DE    & IP    & 5 \\
    166   & Indl. Prod. - Interm \& Cap. Goods & DE    & IP    & 5 \\
    167   & Indl. Prod. - Fab Mtl Prds Excl Mach \& Eqp & ES    & IP    & 5 \\
    168   & Indl. Prod.  & ES    & IP    & 5 \\
    169   & Indl. Prod. - Cons. Goods  & ES    & IP    & 5 \\
    170   & Indl. Prod. - Cap. Goods  & ES    & IP    & 5 \\
    171   & Indl. Prod. - Int. Goods  & ES    & IP    & 5 \\
    172   & Indl. Prod. - Energy  & ES    & IP    & 5 \\
    173   & Indl. Prod. - Cons. Goods, Non-Durables & ES    & IP    & 5 \\
    174   & Indl. Prod. - Mining  & ES    & IP    & 5 \\
    175   & Indl. Prod. - Mfg Ind  & ES    & IP    & 5 \\
    176   & Indl. Prod. - Other Mining \& Quarrying  & ES    & IP    & 5 \\
    177   & Indl. Prod. - Textile  & ES    & IP    & 5 \\
    178   & Indl. Prod. - Chemicals \& Chemical Products  & ES    & IP    & 5 \\
    179   & Indl. Prod. - Plastic \& Rubber Products  & ES    & IP    & 5 \\
    180   & Indl. Prod. - Other Non-Metal Mineral Products  & ES    & IP    & 5 \\
    181   & Indl. Prod. - Metal Processing Ind  & ES    & IP    & 5 \\
    182   & Indl. Prod. - Metal Products Excl. Machinery  & ES    & IP    & 5 \\
    183   & Indl. Prod. - Electrical Equipment  & ES    & IP    & 5 \\
    184   & Indl. Prod. - Automobile  & ES    & IP    & 5 \\
    185   & Euro Interbank Offered Rate - 3-Month (Mean) ($\pmb{acf:\  0.67}$) & EA    & M\&IR & 5 \\
    186   & Money Supply: Loans To Other Ea Residents Excl. Govt. & EA    & M\&IR & 5 \\
    187   & Money Supply: M3  & EA    & M\&IR & 5 \\
    188   & Euro Short Term Repo Rate & FR    & M\&IR & 5 \\
    189   & Datastream Euro Share Price Index (Mth. Avg.) & FR    & M\&IR & 1 \\
    190   & Euribor: 3-Month (Mth. Avg.) & FR    & M\&IR & 5 \\
    191   & Mfi Loans To Resident Private Sector & FR    & M\&IR & 5 \\
    192   & Money Supply - M1 & FR    & M\&IR & 5 \\
    193   & Money Supply - M3 & FR    & M\&IR & 5 \\
    194   & Share Price Index - Sbf 250 & DE    & M\&IR & 1 \\
    195   & Fibor - 3 Month (Mth.Avg.) & DE    & M\&IR & 5 \\
    196   & Money Supply - M3 & DE    & M\&IR & 5 \\
    197   & Money Supply - M2 & DE    & M\&IR & 5 \\
    198   & Bank Prime Lending Rate / Ecb Marginal Lending Facility & DE    & M\&IR & 5 \\
    199   & Dax Share Price Index, Ep & IT    & M\&IR & 1 \\
    200   & Interbank Deposit Rate-Average On 3-Months Deposits & IT    & M\&IR & 5 \\
    201   & Official Reserve Assets & ES    & M\&IR & 5 \\
    202   & Money Supply: M3 - Spanish  & ES    & M\&IR & 5 \\
    203   & Madrid S.E - General Index & ES    & M\&IR & 5 \\
    204   & Hicp - Overall Index ($\pmb{acf:\  -0.54}$) & EA    & HCPI  & 6 \\
    205   & Hicp - All-Items Excluding Energy, Index & EA    & HCPI  & 6 \\
    206   & Hicp - Food Incl. Alcohol And Tobacco, Index & EA    & HCPI  & 6 \\
    207   & Hicp - Processed Food Incl. Alcohol And Tobacco, Index & EA    & HCPI  & 6 \\
    208   & Hicp - Unprocessed Food, Index & EA    & HCPI  & 6 \\
    209   & Hicp - Goods, Index & EA    & HCPI  & 6 \\
    210   & Hicp - Industrial Goods, Index & EA    & HCPI  & 6 \\
    211   & Hicp - Industrial Goods Excluding Energy, Index & EA    & HCPI  & 6 \\
    212   & Hicp - Services, Index & EA    & HCPI  & 6 \\
    213   & Hicp - All-Items Excluding Tobacco, Index & EA    & HCPI  & 6 \\
    214   & Hicp - All-Items Excluding Energy And Food, Index & EA    & HCPI  & 6 \\
    215   & Hicp - All-Items Excluding Energy And Unprocessed Food, Index & EA    & HCPI  & 6 \\
    216   & All-Items Hicp & DE    & HCPI  & 6 \\
    217   & All-Items Hicp & ES    & HCPI  & 6 \\
    218   & All-Items Hicp & FR    & HCPI  & 6 \\
    219   & All-Items Hicp & IT    & HCPI  & 6 \\
    220   & Goods (Overall Index Excluding Services) & DE    & HCPI  & 6 \\
    221   & Goods (Overall Index Excluding Services) & FR    & HCPI  & 6 \\
    222   & Processed Food Including Alcohol And Tobacco & DE    & HCPI  & 6 \\
    223   & Processed Food Including Alcohol And Tobacco & ES    & HCPI  & 6 \\
    224   & Processed Food Including Alcohol And Tobacco & FR    & HCPI  & 6 \\
    225   & Processed Food Including Alcohol And Tobacco & IT    & HCPI  & 6 \\
    226   & Unprocessed Food & DE    & HCPI  & 6 \\
    227   & Unprocessed Food & ES    & HCPI  & 6 \\
    228   & Unprocessed Food & FR    & HCPI  & 6 \\
    229   & Unprocessed Food & IT    & HCPI  & 6 \\
    230   & Non-Energy Industrial Goods & DE    & HCPI  & 6 \\
    231   & Non-Energy Industrial Goods & FR    & HCPI  & 6 \\
    232   & Services (Overall Index Excluding Goods) & DE    & HCPI  & 6 \\
    233   & Services (Overall Index Excluding Goods) & FR    & HCPI  & 6 \\
    234   & Overall Index Excluding Tobacco & DE    & HCPI  & 6 \\
    235   & Overall Index Excluding Tobacco & FR    & HCPI  & 6 \\
    236   & Overall Index Excluding Energy & DE    & HCPI  & 6 \\
    237   & Overall Index Excluding Energy & FR    & HCPI  & 6 \\
    238   & Overall Index Excluding Energy And Unprocessed Food & DE    & HCPI  & 6 \\
    239   & Overall Index Excluding Energy And Unprocessed Food & FR    & HCPI  & 6 \\
    240   & Ppi: Ind Excluding Constr. ($\pmb{acf:\  -0.62}$) \& Energy  & EA    & PPI   & 6 \\
    241   & Ppi: Cap. Goods  & EA    & PPI   & 6 \\
    242   & Ppi: Non-Durable Cons. Goods  & EA    & PPI   & 6 \\
    243   & Ppi: Int. Goods  & EA    & PPI   & 6 \\
    244   & Ppi: Non Dom. - Mining, Mfg \& Quarrying  & EA    & PPI   & 6 \\
    245   & Ppi: Non Dom. Mfg  & DE    & PPI   & 6 \\
    246   & Ppi: Int. Goods Excluding Energy & DE    & PPI   & 6 \\
    247   & Ppi: Cap. Goods & DE    & PPI   & 6 \\
    248   & Ppi: Cons. Goods & DE    & PPI   & 6 \\
    249   & Ppi: Fuel & DE    & PPI   & 6 \\
    250   & Ppi: Indl. Products (Excl. Energy) & DE    & PPI   & 6 \\
    251   & Ppi: Machinery & DE    & PPI   & 6 \\
    252   & Deflated T/O: Ret. Sale In Non-Spcld Str With Food, Bev \& Tob ($\pmb{acf:\  -0.47}$)& DE    & T/O   & 5 \\
    253   & Deflated T/O: Oth Ret. Sale In Non-Spcld Str & DE    & T/O   & 5 \\
    254   & Deflated T/O: Sale Of Motor Vehicle Pts \& Acces & DE    & T/O   & 5 \\
    255   & Deflated T/O: Wholesale Of Agl Raw Matls \& Live Animals & DE    & T/O   & 5 \\
    256   & Deflated T/O: Wholesale Of Household Goods & IT    & T/O   & 5 \\
    257   & T/O: Ret. Trd, Exc Of Mv , Motorcyles \& Fuel & ES    & T/O   & 5 \\
    258   & T/O: Ret. Sale Of Clth \& Leath Gds In Spcld Str & ES    & T/O   & 5 \\
    259   & T/O: Ret. Sale Of Non-Food Prds (Exc Fuel) & ES    & T/O   & 5 \\
    260   & T/O: Ret. Sale Of Info, Househld \& Rec Eqp In Spcld Str & ES    & T/O   & 5 \\
    261   & Ek Unemployment: All ($\pmb{acf:\  0.76}$) & EA    & HUR   & 5 \\
    262   & Ek Unemployment: Persons Over 25 Years Old  & EA    & HUR   & 5 \\
    263   & Ek Unemployment: Women Under 25 Years Old  & EA    & HUR   & 5 \\
    264   & Ek Unemployment: Women Over 25 Years Old  & EA    & HUR   & 5 \\
    265   & Ek Unemployment: Men Over 25 Years Old  & EA    & HUR   & 5 \\
    266   & Fr Hur All Persons (All Ages)  & FR    & HUR   & 5 \\
    267   & Fr Hur Femmes (Ages 15-24)  & FR    & HUR   & 5 \\
    268   & Fr Hur Femmes (All Ages)  & FR    & HUR   & 5 \\
    269   & Fr Hur Hommes (Ages 15-24)  & FR    & HUR   & 5 \\
    270   & Fr Hur Hommes (All Ages)  & FR    & HUR   & 5 \\
    271   & Fr Hur All Persons (Ages 15-24)  & FR    & HUR   & 5 \\
    272   & Fr Hurall Persons(Ages 25 And Over)  & FR    & HUR   & 5 \\
    273   & Fr Hur Females (Ages 25 And Over)  & FR    & HUR   & 5 \\
    274   & Fr Hur Males (Ages 25 And Over)  & FR    & HUR   & 5 \\
    275   & Bd Hur All Persons (All Ages)  & DE    & HUR   & 5 \\
    276   & Bd Hur Femmes (Ages 15-24)  & DE    & HUR   & 5 \\
    277   & Bd Hur Femmes (All Ages)  & DE    & HUR   & 5 \\
    278   & Bd Hur Hommes (Ages 15-24)  & DE    & HUR   & 5 \\
    279   & Bd Hur Hommes (All Ages)  & DE    & HUR   & 5 \\
    280   & Bd Hur All Persons (Ages 15-24)  & DE    & HUR   & 5 \\
    281   & Bd Hurall Persons(Ages 25 And Over)  & DE    & HUR   & 5 \\
    282   & Bd Hur Females (Ages 25 And Over)  & DE    & HUR   & 5 \\
    283   & Bd Hur Males (Ages 25 And Over)  & DE    & HUR   & 5 \\
    284   & It Hur All Persons (All Ages)  & IT    & HUR   & 5 \\
    285   & It Hur Femmes (All Ages)  & IT    & HUR   & 5 \\
    286   & It Hur Hommes (All Ages)  & IT    & HUR   & 5 \\
    287   & It Hur All Persons (Ages 15-24)  & IT    & HUR   & 5 \\
    288   & It Hurall Persons(Ages 25 And Over)  & IT    & HUR   & 5 \\
    289   & Es Hur All Persons (All Ages)  & ES    & HUR   & 5 \\
    290   & Es Hur Femmes (Ages 16-24)  & ES    & HUR   & 5 \\
    291   & Es Hur Femmes (All Ages)  & ES    & HUR   & 5 \\
    292   & Es Hur Hommes (Ages 16-24)  & ES    & HUR   & 5 \\
    293   & Es Hur Hommes (All Ages)  & ES    & HUR   & 5 \\
    294   & Es Hur All Persons (Ages 16-24)  & ES    & HUR   & 5 \\
    295   & Es Hurall Persons(Ages 25 And Over)  & ES    & HUR   & 5 \\
    296   & Es Hur Females (Ages 25 And Over)  & ES    & HUR   & 5 \\
    297   & Es Hur Males (Ages 25 And Over)  & ES    & HUR   & 5 \\
    298   & De - Service  Confidence Indicator ($\pmb{acf:\  0.96}$) & DE    & SI    & 1 \\
    299   & De Services - Buss. Dev. Past 3 Months & DE    & SI    & 1 \\
    300   & De Services - Evol. Demand Past 3 Months & DE    & SI    & 1 \\
    301   & De Services - Exp. Demand Next 3 Months & DE    & SI    & 1 \\
    302   & De Services - Evol. Employ. Past 3 Months & DE    & SI    & 1 \\
    303   & Fr - Service  Confidence Indicator & FR    & SI    & 1 \\
    304   & Fr Services - Buss. Dev. Past 3 Months & FR    & SI    & 1 \\
    305   & Fr Services - Evol. Demand Past 3 Months & FR    & SI    & 1 \\
    306   & Fr Services - Exp. Demand Next 3 Months & FR    & SI    & 1 \\
    307   & Fr Services - Evol. Employ. Past 3 Months & FR    & SI    & 1 \\
    308   & Fr Services - Exp. Employ. Next 3 Months & FR    & SI    & 1 \\
    309   & Fr Services - Exp. Prices Next 
    3 Months & FR    & SI    & 1 \\

\end{longtable}

\normalsize

\section{Five most frequently selected predictors}\label{MostSelPred}
Table~\ref{Tab:VarSel} reports the list of the top 5 predictors in terms of selection frequency across forecasting samples obtained from the empirical application in Section~\ref{EmpApp} of the main text. Regardless of the forecasting horizon $h$, the top predictor for ARMAr-LAS is the Goods Index. The other top predictors, also in the HCPI domain, include EA measurements (e.g.,~Services Index), or are specific to France and Germany (e.g.,~All-Items).
\begin{table}[t]
\centering
\caption{\footnotesize 
Five most frequently selected predictors. Selection percentages are ratios between the times a predictor appears in a forecast and the total number of forecasts (120 for $h$=12 and 96 for $h=24$). 
}\label{Tab:VarSel}
\scalebox{0.75}{
\begin{tabular}{@{}cllll@{}}\hline
Rank & & \multicolumn{3}{c}{Selected Variables} \\\cline{1-1} \cline{3-5}
& & $h$=12 && $h=24$ \\\cline{3-3} \cline{5-5} 
I°  & &Goods, Index & &Goods, Index \\
& &85.8\% & &85.4\%  \\
II° & &Industrial Goods, Index & &Services, Index \\
& &47.5\% & &43.8\%  \\
III°& &Services, Index & &All-Items (De) \\
& &40.8\%  & &35.4\%  \\
IV° & &All-Items Excluding Tobacco, Index & &All-Items Excluding Tobacco, Index\\
& &32.5\% & &32.3\% \\
V°  & &All-Items (Fr) & &Industrial Goods, Index \\
& &24.2\% & &30.2\% \\
\hline
\end{tabular}
}
\end{table}

\end{appendices}

\bibliographystyle{Chicago} 

\bibliography{Bibliography-MM-MC}

@article{SeW2002a,
 ISSN = {01621459},
 URL = {http://www.jstor.org/stable/3085839},
 author = {James H. Stock and Mark W. Watson},
 journal = {Journal of the American Statistical Association},
 number = {460},
 pages = {1167--1179},
 publisher = {[American Statistical Association, Taylor \& Francis, Ltd.]},
 title = {Forecasting Using Principal Components from a Large Number of Predictors},
 volume = {97},
 year = {2002}
}

@article{SeW2002b,
author = {James H Stock and Mark W Watson},
title = {Macroeconomic Forecasting Using Diffusion Indexes},
journal = {Journal of Business \& Economic Statistics},
volume = {20},
number = {2},
pages = {147-162},
year  = {2002},
publisher = {Taylor \& Francis}
}

@article{Lippi2000,
title = {The Generalized Dynamic-Factor Model: Identification And Estimation},
author = {Forni, Mario and Hallin, Marc and Lippi, Marco and Reichlin, Lucrezia},
year = {2000},
journal = {The Review of Economics and Statistics},
volume = {82},
number = {4},
pages = {540-554}
}

@article{Lippi2005,
  author = {Mario Forni and Marc Hallin and Marco Lippi and Lucrezia Reichlin},
title = {The Generalized Dynamic Factor Model: One-Sided Estimation and Forecasting},
journal = {Journal of the American Statistical Association},
volume = {100},
number = {471},
pages = {830-840},
year  = {2005},
publisher = {Taylor & Francis},
doi = {10.1198/016214504000002050},
URL = { 
        https://doi.org/10.1198/016214504000002050
},
eprint = { 
        https://doi.org/10.1198/016214504000002050  
}
}

@article{tibshirani96,
  author = {Tibshirani, R.},
  journal = {Journal of the Royal Statistical Society {Series B}},
  pages = {267-288},
  title = {Regression Shrinkage and Selection via the Lasso},
  volume = {58},
  year = {1996}
}

@ARTICLE{demol2008,
title = {Forecasting using a large number of predictors: Is Bayesian shrinkage a valid alternative to principal components?},
author = {De Mol, Christine and Giannone, Domenico and Reichlin, Lucrezia},
year = {2008},
journal = {Journal of Econometrics},
volume = {146},
number = {2},
pages = {318-328},
url = {https://EconPapers.repec.org/RePEc:eee:econom:v:146:y:2008:i:2:p:318-328}
}

@article{zhaoyu2006,
 author = {Zhao, Peng and Yu, Bin},
 title = {On Model Selection Consistency of Lasso},
 journal = {The Journal of Machine Learning Research},
 issue_date = {12/1/2006},
 volume = {7},
 month = dec,
 year = {2006},
 issn = {1532-4435},
 pages = {2541--2563},
 numpages = {23},
 url = {http://dl.acm.org/citation.cfm?id=1248547.1248637},
 acmid = {1248637},
 publisher = {JMLR.org}
}

@article{Fan2020,
title = "Factor-adjusted regularized model selection",
journal = "Journal of Econometrics",
volume = "216",
number = "1",
pages = "71 - 85",
year = "2020",
note = "Annals Issue in honor of George Tiao: Statistical Learning for Dependent Data",
issn = "0304-4076",
doi = "https://doi.org/10.1016/j.jeconom.2020.01.006",
url = "http://www.sciencedirect.com/science/article/pii/S0304407620300117",
author = "Jianqing Fan and Yuan Ke and Kaizheng Wang"
}

@article{Fan13,
      doi = {10.1093/nsr/nwt032},
      url = {https://doi.org/10.1093/nsr/nwt032},
    year = 2014,
	month = {feb},
	publisher = {Oxford University Press ({OUP})},
	volume = {1},
	number = {2},
	pages = {293--314},
	author = {Jianqing Fan and Fang Han and Han Liu},
	title = {Challenges of Big Data analysis},
	journal = {National Science Review}
}

@article{Fan16,
  author  = {Jianqing Fan and Wen-Xin Zhou},
  title   = {Guarding against Spurious Discoveries in High Dimensions},
  journal = {Journal of Machine Learning Research},
  year    = {2016},
  volume  = {17},
  number  = {203},
  pages   = {1-34},
  url     = {http://jmlr.org/papers/v17/16-068.html}
}

@article{AhnHorenstein2013,
author = {Ahn, Seung C. and Horenstein, Alex R.},
title = {Eigenvalue Ratio Test for the Number of Factors},
journal = {Econometrica},
volume = {81},
number = {3},
pages = {1203-1227},
doi = {https://doi.org/10.3982/ECTA8968},
url = {https://onlinelibrary.wiley.com/doi/abs/10.3982/ECTA8968},
eprint = {https://onlinelibrary.wiley.com/doi/pdf/10.3982/ECTA8968},
year = {2013}
}

@article{Bickel2009,
   title={Simultaneous analysis of Lasso and Dantzig selector},
   volume={37},
   ISSN={0090-5364},
   url={http://dx.doi.org/10.1214/08-AOS620},
   DOI={10.1214/08-aos620},
   number={4},
   journal={The Annals of Statistics},
   publisher={Institute of Mathematical Statistics},
   author={Bickel, Peter J. and Ritov, Ya’acov and Tsybakov, Alexandre B.},
   year={2009},
   month={Aug},
   pages={1705–1732}
}

@article{Negahban2012,
   title={A Unified Framework for High-Dimensional Analysis of $M$-Estimators with Decomposable Regularizers},
   volume={27},
   ISSN={0883-4237},
   url={http://dx.doi.org/10.1214/12-STS400},
   DOI={10.1214/12-sts400},
   number={4},
   journal={Statistical Science},
   publisher={Institute of Mathematical Statistics},
   author={Negahban, Sahand N. and Ravikumar, Pradeep and Wainwright, Martin J. and Yu, Bin},
   year={2012},
   month={Nov}
}

@article{James2013,
  author = {James, Gareth and Witten, Daniela and Hastie, Trevor and Tibshirani, Robert},
  publisher = {Springer},
  title = {An Introduction to Statistical Learning: with Applications in R },
  url = {https://faculty.marshall.usc.edu/gareth-james/ISL/},
  year = 2013
}

@article{Ale2021,
title = {Nowcasting monthly GDP with big data: A model averaging approach},
author = {Proietti, Tommaso and Giovannelli, Alessandro},
year = {2021},
journal = {Journal of the Royal Statistical Society Series A},
volume = {184},
number = {2},
pages = {683-706},
url = {https://EconPapers.repec.org/RePEc:bla:jorssa:v:184:y:2021:i:2:p:683-706}
}

@Article{Bartlett35,
  author =       "M. S. Bartlett",
  title =        "Some Aspects of the Time-Correlation Problem in Regard to Tests of Significance",
  journal =      "j-J-R-STAT-SOC-SUPPL",
  volume =       "98",
  number =       "3",
  pages =        "536--543",
  year =         "1935",
  journal-URL =  "https://www.jstor.org/stable/2342284"
}

@ARTICLE{Uematsu2019,
title = {High‐dimensional macroeconomic forecasting and variable selection via penalized regression},
author = {Uematsu, Yoshimasa and Tanaka, Shinya},
year = {2019},
journal = {Econometrics Journal},
volume = {22},
number = {1},
pages = {34-56},
url = {https://EconPapers.repec.org/RePEc:oup:emjrnl:v:22:y:2019:i:1:p:34-56.}
}

@article{Mizon1995,
  title={A simple message for autocorrelation correctors: Don't},
  author={Mizon, Grayham E},
  journal={Journal of Econometrics},
  volume={69},
  number={1},
  pages={267--288},
  year={1995},
  publisher={Elsevier}
}

@article{CocOrc1949,
 ISSN = {01621459},
 URL = {http://www.jstor.org/stable/2280349},
 author = {D. Cochrane and G. H. Orcutt},
 journal = {Journal of the American Statistical Association},
 number = {245},
 pages = {32--61},
 publisher = {[American Statistical Association, Taylor & Francis, Ltd.]},
 title = {Application of Least Squares Regression to Relationships Containing Auto- Correlated Error Terms},
 urldate = {2022-07-22},
 volume = {44},
 year = {1949}
}

@article{zhang2012,
  title={A General Theory of Concave Regularization for High-Dimensional Sparse Estimation Problems},
  author={Zhang, Cun-Hui and Zhang, Tong},
  journal={Statistical Science},
  pages={576--593},
  year={2012},
  publisher={JSTOR}
}

@book{Buhlmann2011,
  author = {B{\"u}hlmann, Peter and van de Geer, Sara},
  doi = {10.1007/978-3-642-20192-9},
  mrclass = {62-07 (60E15 62G05 62G20 62G30 62H30 62J07)},
  mrnumber = {2807761 (2012e:62006)},
  mrreviewer = {Pierre Alquier},
  note = {Methods, theory and applications},
  pages = {xviii+556},
  publisher = {Springer, Heidelberg},
  series = {Springer Series in Statistics},
  title = {Statistics for high-dimensional data},
  url = {http://dx.doi.org/10.1007/978-3-642-20192-9},
  year = {2011}
}

@book{davidson1994,
  title={Stochastic limit theory: An introduction for econometricians},
  author={Davidson, James},
  year={1994},
  publisher={OUP Oxford}
}

@article{ADAMEK2023,
title = {Lasso inference for high-dimensional time series},
journal = {Journal of Econometrics},
volume = {235},
number = {2},
pages = {1114-1143},
year = {2023},
issn = {0304-4076},
doi = {https://doi.org/10.1016/j.jeconom.2022.08.008},
url = {https://www.sciencedirect.com/science/article/pii/S0304407622001804},
author = {Robert Adamek and Stephan Smeekes and Ines Wilms}
}

@article{nardi2011,
title = {Autoregressive process modeling via the Lasso procedure},
journal = {Journal of Multivariate Analysis},
volume = {102},
number = {3},
pages = {528-549},
year = {2011},
issn = {0047-259X},
doi = {https://doi.org/10.1016/j.jmva.2010.10.012},
url = {https://www.sciencedirect.com/science/article/pii/S0047259X10002186},
author = {Y. Nardi and A. Rinaldo}
}

@article{hansen2019, 
title={THE FACTOR-LASSO AND K-STEP BOOTSTRAP APPROACH FOR INFERENCE IN HIGH-DIMENSIONAL ECONOMIC APPLICATIONS}, 
volume={35}, 
DOI={10.1017/S0266466618000245}, 
number={3}, 
journal={Econometric Theory}, 
publisher={Cambridge University Press}, 
author={Hansen, Christian and Liao, Yuan}, 
year={2019}, 
pages={465–509}
}

@misc{chronopoulos2023,
      title={High Dimensional Generalised Penalised Least Squares}, 
      author={Ilias Chronopoulos and Katerina Chrysikou and George Kapetanios},
      year={2023},
      eprint={2207.07055},
      archivePrefix={arXiv},
      primaryClass={econ.EM},
      url={https://arxiv.org/abs/2207.07055} 
}

@article{Panopoulou2004,
    author = {Panopoulou, Ekaterini and Pittis, Nikitas},
    title = {A comparison of autoregressive distributed lag and dynamic OLS cointegration estimators in the case of a serially correlated cointegration error},
    journal = {The Econometrics Journal},
    volume = {7},
    number = {2},
    pages = {585-617},
    year = {2004},
    month = {11},
    issn = {1368-4221},
    doi = {10.1111/j.1368-423X.2004.00145.x}
}

@article{Babii2024,
    author = {Babii, Andrii and Ghysels, Eric and Striaukas, Jonas},
    title = {High-Dimensional Granger Causality Tests with an Application to VIX and News*},
    journal = {Journal of Financial Econometrics},
    volume = {22},
    number = {3},
    pages = {605-635},
    year = {2022},
    month = {07},
    issn = {1479-8409},
    doi = {10.1093/jjfinec/nbac023},
    url = {https://doi.org/10.1093/jjfinec/nbac023},
    eprint = {https://academic.oup.com/jfec/article-pdf/22/3/605/58229995/nbac023.pdf},
}

@article{McCracken2016,
author = {Michael W. McCracken and Serena Ng},
title = {FRED-MD: A Monthly Database for Macroeconomic Research},
journal = {Journal of Business \& Economic Statistics},
volume = {34},
number = {4},
pages = {574--589},
year = {2016},
publisher = {ASA Website},
doi = {10.1080/07350015.2015.1086655},
URL = { 
        https://doi.org/10.1080/07350015.2015.1086655
},
eprint = { 
        https://doi.org/10.1080/07350015.2015.1086655
}
}

@article{Chernozhukov2021,
author = {Victor Chernozhukov and Wolfgang Karl H{\"a}rdle and Chen Huang and Weining Wang},
title = {{LASSO-driven inference in time and space}},
volume = {49},
journal = {The Annals of Statistics},
number = {3},
publisher = {Institute of Mathematical Statistics},
pages = {1702 -- 1735},
year = {2021},
doi = {10.1214/20-AOS2019},
URL = {https://doi.org/10.1214/20-AOS2019}
}

@article{medeiros2017ADA,
  title={Adaptive LASSO estimation for ARDL models with GARCH innovations},
  author={Medeiros, Marcelo C and Mendes, Eduardo F},
  journal={Econometric Reviews},
  volume={36},
  number={6-9},
  pages={622--637},
  year={2017},
  publisher={Taylor \& Francis}
}

@article{hannan1980,
author = {Hannan, E. J.},
title = {The Estimation of the Order of an {ARMA} Process},
journal = {The Annals of Statistics},
year = {1980},
volume = {8},
number = {5},
pages = {1071--1081},
doi = {10.1214/aos/1176345079},
url = {https://doi.org/10.1214/aos/1176345079}
}

@book{brockwell2016,
  author = {Brockwell, Peter J. and Davis, Richard A.},
  doi = {10.1007/978-3-319-29854-2},
  publisher = {Springer International Publishing},
  series = {Springer Texts in Statistics},
  title = {Introduction to Time Series and Forecasting},
  url = {https://doi.org/10.1007/978-3-319-29854-2},
  year = 2016
}

@article{racine1997,
  author  = {Racine, Jeffrey},
  title   = {Consistent Cross-Validatory Model-Selection for Dependent Data: hv-Block Cross-Validation},
  journal = {Journal of Econometrics},
  volume  = {76},
  pages   = {79--99},
  year    = {1997}
}

@article{WuWu2016,
  author  = {Wu, Wei-Biao and Wu, Ying Nian},
  title   = {Performance bounds for parameter estimates of high-dimensional linear models with correlated errors},
  journal = {Electronic Journal of Statistics},
  volume  = {10},
  number  = {1},
  pages   = {352--379},
  year    = {2016},
  doi     = {10.1214/16-EJS1108}
}

@article{Mcgregor1965,
    author = {Mcgregor, J. R. and Bielenstein, U. M.},
    title = {The approximate distribution of the correlation between two stationary linear Markov series.II†},
    journal = {Biometrika},
    volume = {52},
    number = {1-2},
    pages = {301-302},
    year = {1965},
    month = {06},
    issn = {0006-3444},
    doi = {10.1093/biomet/52.1-2.301},
    url = {https://doi.org/10.1093/biomet/52.1-2.301},
    eprint = {https://academic.oup.com/biomet/article-pdf/52/1-2/301/1029252/52-1-2-301.pdf},
}

@article{BUHLMANN2013,
title = {Correlated variables in regression: Clustering and sparse estimation},
journal = {Journal of Statistical Planning and Inference},
volume = {143},
number = {11},
pages = {1835-1858},
year = {2013},
issn = {0378-3758},
doi = {https://doi.org/10.1016/j.jspi.2013.05.019},
url = {https://www.sciencedirect.com/science/article/pii/S0378375813001225},
author = {Peter Bühlmann and Philipp Rütimann and Sara {van de Geer} and Cun-Hui Zhang}
}

\end{document}